\documentclass[11pt]{article}
\usepackage{pxfonts}
\usepackage{yfonts}
\usepackage{dsfont}
\usepackage{marvosym}
\usepackage{graphicx}
\usepackage{relsize}
\usepackage{amsfonts}
\usepackage{bbm}

\def\bel{\begin{equation}\label}
\def\eeq{\end{equation}}
\def\ds{\displaystyle}
\def\endproof{\hphantom{MM}
\hfill\llap{$\square$}\goodbreak}

\def\mt{\longrightarrow}
\def\v{\vskip 1em}
\def\vsk{\vskip 40em}
\def\ve{\varepsilon}
\def\R{\mathds R}
\def\Z{\mathds Z}
\def\C{\mathfrak{B}}

\def\Cx{\mathds C}
\def\N{{\bf N}}

\def\Re{{\bf Re}}
\def\Im {{\bf Im}}
\def\S{{\bf S}}

\def\O{{\bf O}}
\def\P{{\bf P}}
\def\Q{{\bf Q}}

\def\D{{\bf D}}

\def\J{{\bf J}}
\def\A{{\bf A}}
\def\B{{\bf B}}
\def\H{{\bf H}}
\def\L{{\bf L}}
\def\U{{\bf U}}
\def\V{{\bf V}}

\def\T{{\bf T}}

\def\p{{\partial}}

\def\i{{\bf i}}
\def\Tilde{\widetilde}
\def\Hat{\widehat}

\def\bar{\overline}
\def\supp{{\bf supp}}
\def\vol{{\bf vol}}

\def\I{{\bf I}}
\def\II{{\bf II}}
\def\III{{\bf III}}

\def\q{\mathfrak{q}}

\def\alpha{\alphaup}
\def\beta{\betaup}
\def\xi{{\xiup}}
\def\eta{{\etaup}}
\def\tau{{\tauup}}
\def\rho{{\rhoup}}
\def\phi{{\phiup}}
\def\psi{{\psiup}}
\def\omega{\omegaup}
\def\varphi{{\varphiup}}
\def\lambda{{\lambdaup}}
\def\c{{\bf c}}

\def\m{{\bf m}}

\def\a{{\bf a}}
\def\b{{\bf b}}

\def\p{\partial}

\def\q{{\deltaup}}

\def\epsilon{\sigma}

\newtheorem{lemma}{Lemma}[section]

\newtheorem{remark}{Remark}[section]

\begin{document}
\[\hbox{\LARGE{\bf $\L^p$-boundedness of the Bochner-Riesz operator}}\]

  \[\hbox{Zipeng Wang}\]
   \begin{abstract}
 In this paper, we give a new approach to the Bochner-Riesz summability.
 As a result, we show that the Bochner-Riesz operator $\S^\q, 0<\Re\q<{1\over 2}$ is bounded on $\L^p(\R^n)$ for ${n-1\over 2n}\leq {1\over p}\leq{n+1\over 2n}$.

\[\hbox{\small{\bf Contents}}\]
\[\begin{array}{lr}\ds
\hbox{1 Introduction}~~~~~~~~~~~~~~~~~~~~~~~~~~~~~~~~~~~~~~~~~~~~~~~~~~~~~~~~~~~~~~~~~~~~~~~~~~~~~~~~~~~~~~~~~~~~~~~~~~~~~~~~~~~~~1
\\ \ds
~~~~~\hbox{1.1 Cone multiplier of negative orders}~~~~~~~~~~~~~~~~~~~~~~~~~~~~~~~~~~~~~~~~~~~~~~~~~~~~~~~~~~~~~~3
\\ \ds
~~~~~\hbox{1.2 A pairing formulation}~~~~~~~~~~~~~~~~~~~~~~~~~~~~~~~~~~~~~~~~~~~~~~~~~~~~~~~~~~~~~~~~~~~~~~~~~~~~~~~~~~~~6
\\ \ds
\hbox{2 Theorem Two implies Theorem One}~~~~~~~~~~~~~~~~~~~~~~~~~~~~~~~~~~~~~~~~~~~~~~~~~~~~~~~~~~~~~~~~~~~10
\\ \ds
\hbox{3 End-point estimates regarding $\II^{\alpha}_{j~m}$ and $\III^{\alpha}_{j~m}$}~~~~~~~~~~~~~~~~~~~~~~~~~~~~~~~~~~~~~~~~~~~~~~~~~~~12
\\ \ds~~~~~
\hbox{3.1 ${^\sharp}\U^{\alpha~\beta}_{j~m}$ and ${^\sharp}\V^{\alpha~\beta}_{j~m}$} ~~~~~~~~~~~~~~~~~~~~~~~~~~~~~~~~~~~~~~~~~~~~~~~~~~~~~~~~~~~~~~~~~~~~~~~~~~~~~~~~~~~~~~~~~~~~~~~~17
\\ \ds~~~~~
\hbox{3.2 ${^\flat}\U^{\alpha~\beta}_{j~m}$ and ${^\flat}\V^{\alpha~\beta}_{j~m}$} ~~~~~~~~~~~~~~~~~~~~~~~~~~~~~~~~~~~~~~~~~~~~~~~~~~~~~~~~~~~~~~~~~~~~~~~~~~~~~~~~~~~~~~~~~~~~~~~~18
\\ \ds
\hbox{4 Proof of Theorem Two}~~~~~~~~~~~~~~~~~~~~~~~~~~~~~~~~~~~~~~~~~~~~~~~~~~~~~~~~~~~~~~~~~~~~~~~~~~~~~~~~~~~~~~~~~~~19
\\ \ds
\hbox{5 Proof of Proposition One}~~~~~~~~~~~~~~~~~~~~~~~~~~~~~~~~~~~~~~~~~~~~~~~~~~~~~~~~~~~~~~~~~~~~~~~~~~~~~~~~~~~~~~~25
\\ \ds~~~~~
\hbox{5.1 Proof of Lemma One}~~~~~~~~~~~~~~~~~~~~~~~~~~~~~~~~~~~~~~~~~~~~~~~~~~~~~~~~~~~~~~~~~~~~~~~~~~~~~~~~~~~~~~26
\\ \ds~~~~~
\hbox{5.2 Size of ${^\sharp}\Hat{\U}^{\alpha~\beta}_{j~m}$ and ${^\sharp}\Hat{\V}^{\alpha~\beta}_{j~m}$}~~~~~~~~~~~~~~~~~~~~~~~~~~~~~~~~~~~~~~~~~~~~~~~~~~~~~~~~~~~~~~~~~~~~~~~~~~~~~~~~~~~~29
\\ \ds
\hbox{6 On the $\L^1$-norm of ${^\flat}\Hat{\U}^{\alpha~\beta}_{j~m}$ and ${^\flat}\Hat{\V}^{\alpha~\beta}_{j~m}$}~~~~~~~~~~~~~~~~~~~~~~~~~~~~~~~~~~~~~~~~~~~~~~~~~~~~~~~~~~~~~~~~~~~~~~~35
\\ \ds
\hbox{7 Proof of Proposition Two}~~~~~~~~~~~~~~~~~~~~~~~~~~~~~~~~~~~~~~~~~~~~~~~~~~~~~~~~~~~~~~~~~~~~~~~~~~~~~~~~~~~~~~~38
\\ \ds~~~~~
\hbox{7.1 Proof of Lemma Two}~~~~~~~~~~~~~~~~~~~~~~~~~~~~~~~~~~~~~~~~~~~~~~~~~~~~~~~~~~~~~~~~~~~~~~~~~~~~~~~~~~~~~~39
\\ \ds
\hbox{A Some estimates regarding Bessel functions}~~~~~~~~~~~~~~~~~~~~~~~~~~~~~~~~~~~~~~~~~~~~~~~~~~~~~~~44
\\ \ds
\hbox{B Fourier transform of $\Lambda^\alpha$}~~~~~~~~~~~~~~~~~~~~~~~~~~~~~~~~~~~~~~~~~~~~~~~~~~~~~~~~~~~~~~~~~~~~~~~~~~~~~~~~~~~~~~~~~47
\\ \ds
\hbox{References}~~~~~~~~~~~~~~~~~~~~~~~~~~~~~~~~~~~~~~~~~~~~~~~~~~~~~~~~~~~~~~~~~~~~~~~~~~~~~~~~~~~~~~~~~~~~~~~~~~~~~~~~~~~~~~~~~~~51
\end{array}\]
\end{abstract}  
  
 \section{Introduction}
 \setcounter{equation}{0} 
 A classical problem in harmonic analysis is to make precise the sense of Fourier inversion formulae, for given $f\in\L^p(\R^n)$. One natural way of formulating such identities is  via the summability method due to Bochner and Riesz. This assertion leads us to study for  the $\L^p$-boundedness of Bochner-Riesz operator, defined as
  \bel{S}
  \S^\q f(x)~=~\int_{\R^n}e^{2\pi\i x\cdot\xi} \Hat{f}(\xi) \left( 1-|\xi|^2\right)^\q_+ d\xi,\qquad \hbox{\small{$\Re\q\ge0$}}.
  \eeq
At $\q=0$, we revisit on the famous ball multiplier problem.  $\S^0$ is bounded only on $\L^2(\R^n)$ proved by Fefferman \cite{Fefferman}.

$\diamond$  {\small Throughout,  $\C>0$ is regarded as a generic constant whose value depends on its sub-indices. }  

$\diamond$ {\small We  always write $\c>0$ for some  large constant.}

$\S^\q f$ in (\ref{S})
 can be written as a convolution 
 \bel{S kernel}
 \begin{array}{cc}\ds
  \S^\q f(x)
  ~=~\pi^{-\q}\Gamma(1+\q)\int_{\R^n}f(x-u) \Omega^\q(u)du,
  \\\\ \ds
  \Omega^\q(u)~=~\left({1\over|u|}\right)^{{n\over 2}+\q} \J_{{n\over 2}+\q}\Big(2\pi|u|\Big)  
\end{array}
\eeq
where $\Gamma$ and $\J$ denote Gamma and Bessel functions. See (\ref{Omega^z})-(\ref{Omega^z Transform}).

Moreover,  the estimate in (\ref{J norm}) implies
\[\left|\Omega^\q(u)\right|~\leq~
\C_{\Re\q}~e^{\c|\Im\q|}~\left({1\over 1+|u|}\right)^{{n+1\over 2}+\Re\q}.\]
Observe that for $\Re\q>{n-1\over 2}$, the kernel of $\S^\q$ is an $\L^1$-function in $\R^n$. We thus have
\bel{Result n-1/2}
 \left\| \S^\q f\right\|_{\L^p(\R^n)}~\leq~\C_{\Re\q}~e^{\c|\Im\q|}~\left\| f\right\|_{\L^p(\R^n)},\qquad 1\leq p\leq \infty.
 \eeq
When $0<\q\leq{n-1\over 2}$,  a  counter example given by Herz \cite{Herz} shows that $\S^\q$ is not bounded on $\L^p(\R^n)$   for either $p\leq {2n\over n+1+2\q}$ or $p\ge {2n\over n-1-2\q}$.  
\v
{\bf Bochner-Riesz Conjecture} ~~{\it Let $\S^\q$ defined in (\ref{S}) for $0<\q\leq {n-1\over 2}$. We have
 \bel{CONJECTURE ONE}
 \left\| \S^\q f\right\|_{\L^p(\R^n)}~\leq~\C_{\q~p}~\left\| f\right\|_{\L^p(\R^n)}
 \eeq
 if and only if
 \bel{p RANGE}
 {n-1-2\q\over 2n}~<~{1\over p}~<~{n+1+2\q\over 2n}. 
 \eeq}\\
 In 1972, Carlson and Sj\"{o}lin \cite{Carlson-Sjolin} proved the conjecture for $n=2$. 
A year later, Fefferman \cite{Fefferman'} observes  the restriction estimate: 
\bel{Restriction Estimate}
\left\{\int_{\mathds{S}^{n-1}} \left| \Hat{f}(\xi)\right|^2 d\sigma(\xi)\right\}^{1\over 2}~\leq~\C_p~\left\| f\right\|_{\L^p(\R^n)},\qquad 1<p<\infty
\eeq
 implying $\S^\q, \q>0$ bounded on $\L^{p\over p-1}(\R^n)$. Here, $d\sigma$ denotes the Lebesgue measure on $\mathds{S}^{n-1}$.

Consequently,  the Tomas-Stein restriction theorem implies the $\L^p$-norm inequality in (\ref{CONJECTURE ONE})
for $p$  belonging  to (\ref{p RANGE}) with an extra condition of $p\ge{2n+2\over n-1}$  or $p\leq {2n+2\over n+3}$. See Tomas \cite{Tomas} or Chapter IX of Stein \cite{Stein}. 
Today, the connection between Fourier restriction theorem and Bochner-Riesz summability is extensively explored. For instance, Tao \cite{Tao} proves (\ref{CONJECTURE ONE})-(\ref{p RANGE}) implying the Restriction conjecture on $\mathds{S}^{n-1}$. On the other hand, Carbery \cite{Carbery} shows that such implication can be reversed if the unit sphere is replaced by a paraboloid.

A number of remarkable results have been accomplished in this broad area, for example by Bourgain \cite{Bourgain}-\cite{Bourgain'}, Bourgain and Guth \cite{Bourgain-Guth},  Guth \cite{Guth}-\cite{Guth'}, Tao and Vargas \cite{Tao-Vargas}, Guo-et-al \cite{Guo-et-al}, Guo, Wang and Zhang \cite{Guo-Wang-Zhang}, Gressman \cite{Gressman}-\cite{Gressman'}, Lee \cite{Lee}   and Wolff \cite{Wolff}. 

Our main result is stated in below.
\v
{\bf Theorem One}~~{\it Let $\S^\q$ defined in (\ref{S}) for $0<\Re\q<{1\over 2}$. We have
 \bel{RESULT ONE}
 \left\| \S^\q f\right\|_{\L^p(\R^n)}~\leq~\C_{\Re\q}~e^{\c|\Im\q|}~\left\| f\right\|_{\L^p(\R^n)},\qquad  {n-1\over 2n}~\leq~{1\over p}~\leq~{n+1\over 2n}. 
 \eeq} \\
From (\ref{Result n-1/2}) and (\ref{RESULT ONE}),  by applying Stein interpolation theorem \cite{Stein'}, we obtain (\ref{CONJECTURE ONE})-(\ref{p RANGE}).

In order to prove {\bf Theorem One},  we next introduce another convolution operator, denoted by $\I^\alpha$ for $0<\Re\alpha<1$ whose Fourier transform is given in terms of  cone multipliers having negative orders. The $\L^p$-boundedness this operator implies {\bf Theorem One}. 

Furthermore, we develop two different types of dyadic decomposition respectively on the frequency space and the physical space. Every consisting partial operator  will be split into two after a {\bf pairing formulation} between these dyadic decompositions. Each one of them satisfies the desired $\L^p$-estimates.

\subsection{Cone multiplier of negative orders}
Let $\S^\q, \Re \q>0$ defined in (\ref{S}).  Our key observation  is 
\[
 \left(1-|\xi|^2\right)^\q_+~=~2\q\int_0^1 \left({1\over \tau^2-|\xi|^2}\right)^{1-\q}_+\tau d\tau,\qquad \hbox{\small{$\Re\q>0$}}.
\]
Let $0<\Re\q<{1\over 2}$. The cone multiplier $\left({1\over \tau^2-|\xi|^2}\right)^{1-\q}_+$  having a negative order has been previously investigated by Gelfand and Shilov \cite{Gelfand-Shilov}. 
 
 Denote $1-\q=\alpha$ for $\alpha\in\Cx$.
 $\Lambda^\alpha$ is a distribution defined in $\R^{n+1}$ by analytic continuation from 
\[
\begin{array}{cc}\ds
\hbox{\small{$\Re\alpha>{n-1\over2}$}},
\qquad
 \Lambda^\alpha(x,t)~=~ \hbox{\small{$\pi^{\alpha-{n+1\over2}}\Gamma^{-1}\left(\alpha-{n-1\over 2}\right)$}}
  \left({1\over t^2-|x|^2}\right)^{{n+1\over 2}-\alpha}_+ .
  \end{array}
\]
For ${1\over 2}<\Re\alpha<1$, the Fourier transform of $\Lambda^\alpha$ agrees with the function
\bel{Lambda Transform}
\begin{array}{lr}
  \Hat{\Lambda}^\alpha(\xi,\tau)~=~
\pi^{{n-1\over 2}-2\alpha}\Gamma\left(\alpha\right)
 \left\{ ~{\ds \left({1\over \tau^2-|\xi|^2}\right)^\alpha_-} - \sin\pi\left(\alpha-{1\over 2}\right){\ds \left( {1\over \tau^2- |\xi|^2}\right)^\alpha_+}~\right\}
  \end{array}
 \eeq
whenever $|\tau|\neq|\xi|$.  See appendix B.

Let $\varphi\in\mathcal{C}^\infty_o(\R)$ be a smooth cut-off function such that $\varphi(t)=1$ if $|t|\leq1$ and $\varphi(t)=0$ for $|t|>2$. Define
\bel{hat phi}
\Hat{\phi}(\xi)~=~\varphi\left({2\over 3}|\xi|\right)-\varphi\left(3|\xi|\right).
\eeq
Note that
$\supp\Hat{\phi}\subset\left\{\xi\in\R^n\colon {1\over 3}<|\xi|\leq3\right\}$ and 
$\Hat{\phi}(\xi)=1$ for ${2\over 3}<|\xi|<{3\over 2}$.
\v

In order to prove {\bf Theorem One}, we consider
 \bel{I alpha}
 \begin{array}{lr}\ds
\I^\alpha f(x)
~=~\int_{\R^n} e^{2\pi\i x\cdot \xi} \Hat{f}(\xi) \Hat{\phi}(\xi)\left\{\int_0^1 \Hat{\Lambda}^\alpha(\xi,\tau) \tau d\tau\right\}d\xi.
\end{array}
\eeq

{\bf Theorem Two} ~~{\it Let $\I^\alpha$ defined in (\ref{Lambda Transform})-(\ref{I alpha}) for ${1\over 2}<\Re\alpha<1$.   We have
 \bel{Result Two} 
 \begin{array}{cc}\ds
  \left\| \I^\alpha f\right\|_{\L^p(\R^n)} ~\leq~\C_{\Re\alpha}~e^{\c|\Im\alpha|}~\left\| f\right\|_{\L^p(\R^n)}, \qquad  {n-1\over 2n}~\leq ~{1\over p}~\leq~ {n+1\over 2n}.  
 \end{array}
  \eeq
  }
\begin{remark}
 {\bf Theorem Two} implies {\bf Theorem One}. \end{remark}

As investigated by Strichartz \cite{Strichartz}, $\Omega^\alpha$ is a distribution defined in $\R^n$ by analytic continuation from
\[
\begin{array}{cc}\ds
\hbox{\small{$\Re\alpha>{n-1\over 2}$}},\qquad \Omega^{\alpha}(x)~=~\hbox{\small{$\pi^{\alpha-{n+1\over 2}}\Gamma^{-1}\left(\alpha-{n-1\over 2}\right)$}} 
 \left({1 \over1-|x|^2}\right)^{{n+1\over 2}-\alpha}_+.
 \end{array}
\]
Equivalently,   it can be  defined by
\[
\begin{array}{lr}\ds
\Hat{\Omega}^{\alpha}(\xi)~=~\left({1\over|\xi|}\right)^{{n\over 2}-\big[{n+1\over 2}-\alpha\big]} \J_{{n\over 2}-\big[{n+1\over 2}-\alpha\big]}\Big(2\pi|\xi|\Big)
\\\\ \ds~~~~~~~~~~
~=~
\left({1\over|\xi|}\right)^{\alpha-{1\over 2}} \J_{\alpha-{1\over 2}}\Big(2\pi|\xi|\Big),\qquad \alpha\in\Cx.
\end{array}
\]
In section 3, we show that $\I^\alpha$ defined in (\ref{I alpha}) for ${1\over 2}<\Re\alpha<1$ can be expressed as
\bel{I alpha new express}
 \begin{array}{cc}\ds
 \I^\alpha f(x)
~=~\int_{\R^n} e^{2\pi\i x\cdot \xi} \Hat{f}(\xi) \Hat{\phi}(\xi) \left\{\int_\R e^{-2\pi\i r} \Hat{\Omega}^\alpha(r\xi) \omega(r) |r|^{2\alpha-1} dr\right\} d\xi,
\\\\ \ds
\omega(r)~=~e^{2\pi\i r}\int_0^1 e^{-2\pi\i \tau r} \tau d\tau~=~{-1\over 2\pi\i}{1\over r}-{1\over 4\pi^2r^2} \left[1-e^{-2\pi\i r}\right].
\end{array}
\eeq
Moreover, define
\bel{I<}
\begin{array}{cc}\ds
 \I^\alpha_< f(x)
~=~\int_{\R^n} e^{2\pi\i x\cdot \xi} \Hat{f}(\xi) \Hat{\P}^\alpha_<(\xi) d\xi,
\\\\ \ds
\Hat{\P}^\alpha_<(\xi)~=~\Hat{\phi}(\xi) \int_{-1}^1 e^{-2\pi\i r} \Hat{\Omega}^\alpha(r\xi) \omega(r) |r|^{2\alpha-1} dr
\end{array}
\eeq
and
 \bel{I j}
 \begin{array}{cc}\ds
  \I^\alpha_j f(x)
~=~\int_{\R^n} e^{2\pi\i x\cdot \xi} \Hat{f}(\xi) \Hat{\P}^\alpha_j(\xi)d\xi,
\\\\ \ds
\Hat{\P}^\alpha_j(\xi)~=~
\Hat{\phi}(\xi) \int_{2^{j-1}\leq|r|<2^j} e^{-2\pi\i r} \Hat{\Omega}^\alpha(r\xi) \omega(r) |r|^{2\alpha-1} dr,\qquad j>0.
\end{array}
\eeq
Later, we shall see $\P^\alpha_<\in\L^1(\R^n)$. Hence that $\I^\alpha_<$ is bounded on $\L^p(\R^n)$ for $1\leq p\leq\infty$.

Our aim is to show
 \bel{Result Two j} 
 \begin{array}{cc}\ds
  \left\| \I^\alpha_j f\right\|_{\L^p(\R^n)} ~\leq~\C_{\Re\alpha}~e^{\c|\Im\alpha|}~2^{-\ve j}~\left\| f\right\|_{\L^p(\R^n)}, \qquad  {n-1\over 2n}~\leq ~{1\over p}~\leq~ {n+1\over 2n}  
 \end{array}
  \eeq
for every $j>0$ and some $\ve=\ve(\Re\alpha)>0$.

Given ${1\over2}<\Re\alpha<1$, $\I^{\alpha~z}_j, 0\leq\Re z\leq1$ is a family of analytic operators  which will be defined explicitly in section 4.   In particular, we have 
$\I^{\alpha~{1\over n}}_j f(x)=\I^\alpha_j f(x)$.

Ideally, we hope to establish  (\ref{Result Two j}) by using complex interpolation as follows. Suppose
\bel{end-point.0 fake}
 \begin{array}{cc}\ds
  \left\| \I^{\alpha~0+\i\Im z}_j f\right\|_{\L^2(\R^n)} ~\leq~\C_{\Re\alpha}~e^{\c|\Im\alpha|}~2^{-\ve j}~\left\| f\right\|_{\L^2(\R^n)} 
 \end{array}
  \eeq
and
\bel{end-point.1}
 \begin{array}{cc}\ds
  \left\| \I^{\alpha~1+\i\Im z}_j f\right\|_{\L^p(\R^n)} ~\leq~\C_{\Re\alpha}~e^{\c|\Im\alpha|}~2^{-\ve j}~\left\| f\right\|_{\L^p(\R^n)},\qquad 1\leq p\leq\infty 
 \end{array}
  \eeq
for every $j>0$ and some $\ve=\ve(\Re\alpha)>0$.

We obtain (\ref{Result Two j}) by applying Stein interpolation theorem \cite{Stein'} and taking into account for ${n-1\over 2n}=0\left({1\over n}\right)+{1\over 2}\left({n-1\over n}\right)$ and ${n+1\over 2n}=1\left({1\over n}\right)+{1\over 2}\left({n-1\over n}\right)$.

Unfortunately, the $\L^2$-estimate in (\ref{end-point.0 fake}) is NOT true. The best result can be concluded is 
\bel{end-point.0 real}
 \begin{array}{cc}\ds
  \left\| \I^{\alpha~0+\i\Im z}_j f\right\|_{\L^2(\R^n)} ~\leq~\C_{\Re\alpha}~e^{\c|\Im\alpha|}~2^{j/2}2^{-\ve j}~\left\| f\right\|_{\L^2(\R^n)},\qquad j>0. 
 \end{array}
  \eeq
 This is consequence of which the  Fourier multiplier regarding every $\I^{\alpha~0+\i\Im z}_j$  has its norm bounded by 
 \[\C_{\Re\alpha}~e^{\c|\Im\alpha|}~\left|{1\over 1-|\xi|}\right|^{{1\over 2}-\ve},\qquad\hbox{\small{$|\xi|\approx 1\pm2^{-j}$}}.\]
To obtain the $\L^p$-estimate in (\ref{Result Two j}), we first write  $\I^\alpha_j$ into a group of partial operators.
 
Given ${1\over 2}<\Re\alpha<1$, denote  $0<\sigma=\sigma(\Re\alpha)<{1\over 2}$ for some constant which can be chosen sufficiently small.  In particular, we shall find $\sigma\mt0$ as $\Re\alpha\mt1$.

For every $j>0$, assert $\lambda_m\in[2^{j-1},2^j]$ where $\lambda_0=2^{j-1}$, $\lambda_M=2^j$ and 
\bel{lambda distance}
 2^{\epsilon j-1}~\leq~\lambda_m-\lambda_{m-1}~<~2^{\epsilon j}.
 \eeq
 \begin{remark} There are at most a constant multiple of $2^{(1-\epsilon)j}$ many $\lambda_m$'s inside $[2^{j-1},2^j]$.
 \end{remark}
Let $\I^\alpha_j$ defined in (\ref{I j}). Consider
\bel{I j m}
 \begin{array}{cc}\ds
  \I^\alpha_j f(x)~=~\sum_{m=1}^M \I^\alpha_{j~m} f(x),\qquad 
  \I^\alpha_{j~m} f(x)
~=~\int_{\R^n} e^{2\pi\i x\cdot \xi} \Hat{f}(\xi) \Hat{\P}^\alpha_{j~m}(\xi)d\xi,
\\\\ \ds
\Hat{\P}^\alpha_{j~m}(\xi)~=~
\Hat{\phi}(\xi) \int_{\lambda_{m-1}\leq|r|<\lambda_m} e^{-2\pi\i r} \Hat{\Omega}^\alpha(r\xi) \omega(r) |r|^{2\alpha-1} dr.
\end{array}
\eeq
We claim
 \bel{Result Two j m} 
 \begin{array}{cc}\ds
  \left\| \I^\alpha_{j~m} f\right\|_{\L^p(\R^n)} ~\leq~\C_{\Re\alpha}~e^{\c|\Im\alpha|}~2^{-(1-\epsilon)j}~2^{-\ve j}~\left\| f\right\|_{\L^p(\R^n)}, \qquad  {n-1\over 2n}~\leq ~{1\over p}~\leq~ {n+1\over 2n}  
 \end{array}
  \eeq
for every $j>0$, $m=1,\ldots,M$ and some $\ve=\ve(\Re\alpha)>0$. 

Clearly, (\ref{Result Two j m}) and {\bf Remark 1.2} together with Minkowski inequality imply  (\ref{Result Two j}).

\subsection{A pairing formulation}
First, we introduce an ingenious dyadic decomposition  initially constructed by Fefferman \cite{Fefferman'} and later refined by Seeger, Sogge and Stein \cite{Seeger-Sogge-Stein} in their study of Fourier integral operators. 
Given $j>0$,  $\left\{ \xi^\nu_j\right\}_\nu$ is a collection of points  almost equally distributed on  $\mathds{S}^{n-1}$ having a grid length between $2^{-j/2-1}$ and $2^{-j/2}$:   

{\bf (1)} $|\xi^\mu_j-\xi^\nu_j|\ge 2^{-j/2-1}$ for every $\xi^\mu_j\neq\xi^\nu_j$. 

{\bf (2)} For any $u\in \mathds{S}^{n-1}$, there is a $\xi^\nu_j$ in the open set $\{\xi\in\mathds{S}^{n-1}\colon |\xi-u|<2^{-j/2+1}\}$.

Let $\varphi\in\mathcal{C}^\infty_o(\R)$  be a smooth cut-off function such that $\varphi(t)=1$ if $|t|\leq1$ and $\varphi(t)=0$ for $|t|>2$. Define
\bel{phi^v_j intro}
\varphi^\nu_j(\xi)~=~{\ds\varphi\Bigg[2^{j/2}\left|{\xi\over |\xi|}-\xi^{\nu}_{j}\right|\Bigg]\over \ds \sum_{\nu\colon\xi^\nu_j\in\mathds{S}^{n-1}} \varphi\Bigg[2^{j/2}\left|{\xi\over |\xi|}-\xi^{\nu}_{j}\right|\Bigg]}
\eeq
whose support is contained in the cone
\bel{Gamma_j}
\Gamma^\nu_j~=~\Bigg\{\xi\in\R^n\colon \left|{\xi\over |\xi|}-\xi^{\nu}_{j}\right|<2^{-j/2+1}\Bigg\}.
\eeq
Observe that by the construction in {\bf (1)}-{\bf(2)}, the support of $\sum_{\nu\colon\xi^\nu_j\in\mathds{S}^{n-1}} \varphi\Bigg[2^{j/2}\left|{\xi\over |\xi|}-\xi^{\nu}_{j}\right|\Bigg]$ is  $\R^n$.  
\begin{remark} There are at most a constant multiple of $2^{\big({n-1\over 2}\big)j}$ many $\xi^\nu_j$ s.
\end{remark}
\begin{figure}[h]
\centering
\includegraphics[scale=0.40]{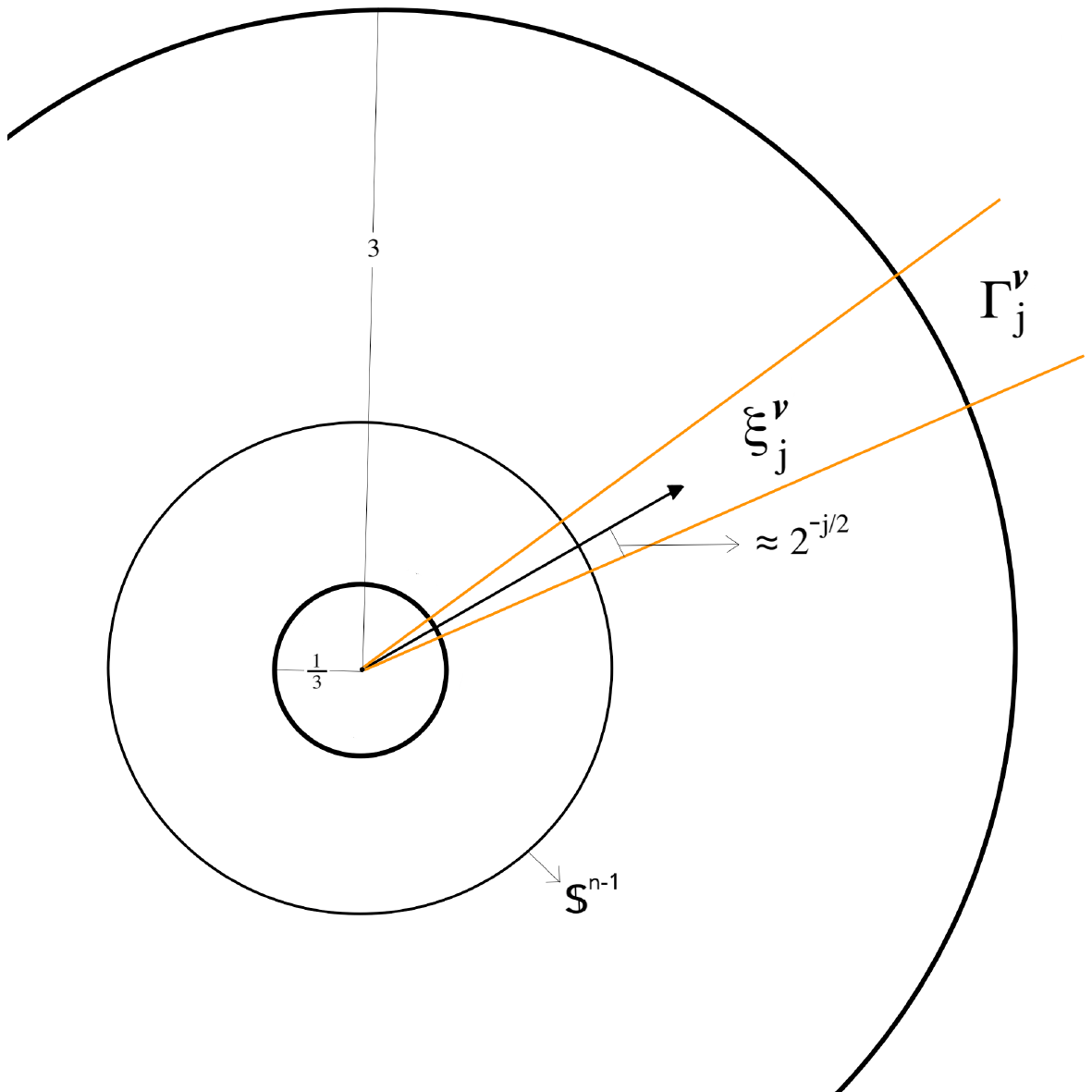}
\end{figure}

Given $\xi_j^\nu$, $\L_\nu$ is an $n\times n$-orthogonal matrix with $\det\L_\nu=1$ and 
$\L_\nu^T \xi^\nu_j=(1,0)^T\in\R\times\R^{n-1}$.
Define
 \bel{Psi mu j alpha}
 \begin{array}{cc}\ds
\Psi^{\alpha}_{j~m}(x)~=~\sum_{\nu\colon\xi^\nu_j\in\mathds{S}^{n-1}}  \Psi^{\alpha~\nu}_{j~m}(x),
 \\\\ \ds 
 \Psi^{\alpha~\nu}_{j~m}(x)~=~\varphi\left[ 2^{-\epsilon j-1}\left|\lambda_m-\left(\L_\nu^T x\right)_1\right|\right] 
\prod_{i=2}^n \varphi\left[2^{-\left[{1\over 2}+\epsilon\right]j}\left|\left(\L_\nu^T x\right)_i\right|\right]
\\ \ds~~~~~~~~~~~~~
~+~\varphi\left[ 2^{-\epsilon j-1}\left|\lambda_m+\left(\L_\nu^T x\right)_1\right|\right] 
\prod_{i=2}^n \varphi\left[2^{-\left[{1\over 2}+\epsilon\right]j}\left|\left(\L_\nu^T x\right)_i\right|\right].
\end{array}
\eeq
Each $\Psi^{\alpha~\nu}_{j~m}$ is supported in the rectangles
\bel{Rectangles}
\begin{array}{lr}\ds
\hbox{R}^\nu_{j~m}~=~
\\\\ \ds
\Bigg\{x\in\R^n\colon\lambda_m-2^{\epsilon j+2}<\left(\L_\nu^T x\right)_1<\lambda_m+2^{\epsilon j+2},~\left|\left(\L_\nu^T x\right)_i\right|<2^{\left[{1\over 2}+\epsilon\right]j+1},~ i=2,\ldots,n
\Bigg\}~\cup
\\ \ds
 \Bigg\{x\in\R^n\colon-\lambda_m-2^{\epsilon j+2}<\left(\L_\nu^T x\right)_1<-\lambda_m+2^{\epsilon j+2},~\left|\left(\L_\nu^T x\right)_i\right|<2^{\left[{1\over 2}+\epsilon\right]j+1},~ i=2,\ldots,n
\Bigg\}.
\end{array}
\eeq
Moreover, $\hbox{R}^\mu_{j~m}\cap \hbox{R}^\nu_{j~m}=\emptyset$ if   $\left|\xi^\mu_j-\xi^\nu_j\right|\ge\c2^{-j/2}2^{\sigma j}$ for  $\c$  large.  
\begin{figure}[h]
\centering
\includegraphics[scale=0.92]{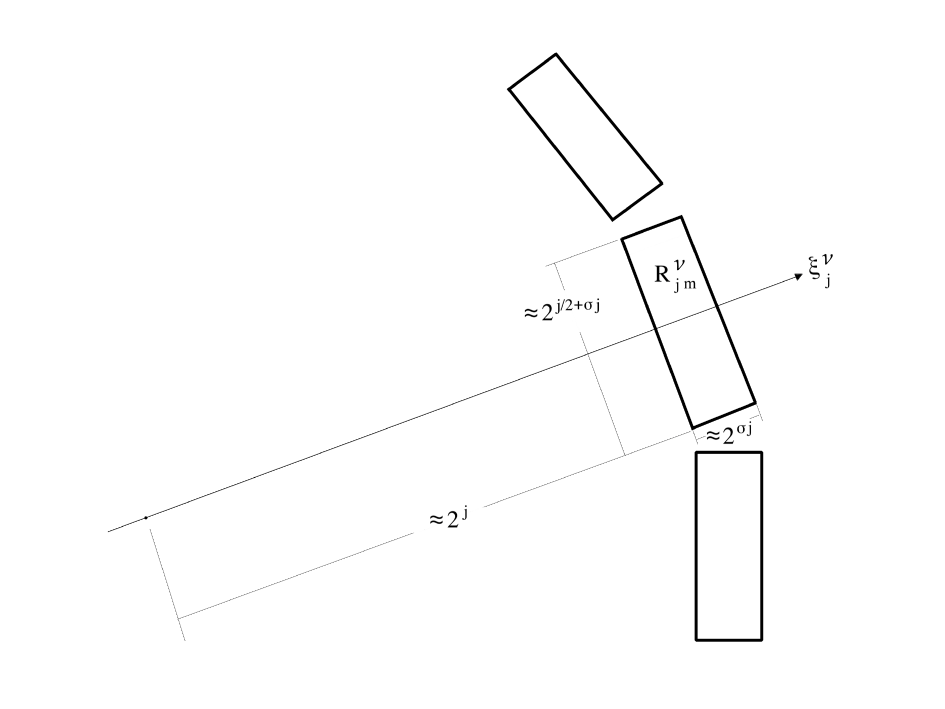}
\end{figure}

Let $\varphi^\nu_j$ defined in (\ref{phi^v_j intro}). We assert
\bel{P v j m}
 \begin{array}{cc}\ds
 \Hat{\P}_{j~m}^\alpha(\xi)~=~
\sum_{\nu\colon \xi^\nu_j\in\mathds{S}^{n-1}}  \Hat{\P}_{j~m}^{\alpha~\nu}(\xi),
\\\\ \ds
 \Hat{\P}_{j~m}^{\alpha~\nu}(\xi)~=~\varphi^\nu_j(\xi) \Hat{\phi}(\xi)\int_{\lambda_{m-1}\leq|r|<\lambda_m} e^{-2\pi\i  r}\Hat{\Omega}^\alpha (r\xi)\omega(r) |r|^{2\alpha-1} dr. 
  \end{array}
\eeq
Recall $\I^\alpha_{j~m} f=f\ast\P_{j~m}^\alpha$ from (\ref{I j m}). We have
\[
\begin{array}{cc}\ds
\I^\alpha_{j~m}f(x)~=~\int_{\R^n} e^{2\pi\i x\cdot\xi} \Hat{f}(\xi)  \Hat{\P}_{j~m}^\alpha(\xi)d\xi,
\qquad
\Hat{\P}_{j~m}^\alpha(\xi)~=~\sum_{\nu\colon\xi^\nu_j\in\mathds{S}^{n-1}} \Hat{\P}_{j~m}^{\alpha~\nu}(\xi).
\end{array}
\]
Let $\Psi^{\alpha}_{j~m}$ and $\Psi^{\alpha~\nu}_{j~m}$ defined in (\ref{Psi mu j alpha}).
Consider
   \bel{II_j U}
 \begin{array}{cc}\ds
 \II^\alpha_{j~m}f(x)~=~\int_{\R^n} f(x-u) \U^\alpha_{j~m}(u)du 
 \\\\ \ds
 \U^\alpha_{j~m}(x) ~=~
 \sum_{\nu\colon\xi^\nu_j\in\mathds{S}^{n-1}}  \U^{\alpha~\nu}_{j~m}(x), 
\qquad
  \U^{\alpha~\nu}_{j~m}(x)~=~\P^{\alpha~\nu}_{j~m}(x) \left[1-\Psi^{\alpha~\nu}_{j~m}(x)\right] 
    \end{array}
  \eeq 
 and
 \bel{III_j V}
 \begin{array}{cc}\ds
 \III^\alpha_{j~m}f(x)~=~\int_{\R^n} f(x-u) \V^\alpha_{j~m}(u)du 
 \\\\ \ds
 \V^\alpha_{j~m}(x) ~=~
 \sum_{\nu\colon\xi^\nu_j\in\mathds{S}^{n-1}}  \V^{\alpha~\nu}_{j~m}(x), 
\qquad
  \V^{\alpha~\nu}_{j~m}(x)~=~\P^{\alpha~\nu}_{j~m}(x) \Psi^{\alpha~\nu}_{j~m}(x). 
    \end{array}
 \eeq 
From (\ref{II_j U})-(\ref{III_j V}), we find
\[
\begin{array}{lr}\ds
\U^{\alpha}_{j~m}(x)+\V^{\alpha}_{j~m}(x)~=~ \sum_{\nu\colon\xi^\nu_j\in\mathds{S}^{n-1}} \U^{\alpha~\nu}_{j~m}(x)+\V^{\alpha~\nu}_{j~m}(x)

\\\\ \ds~~~~~~~~~~~~~~~~~~~~~~~~~~~~~
~=~\sum_{\nu\colon\xi^\nu_j\in\mathds{S}^{n-1}} \P^{\alpha~\nu}_{j~m}(x) \left[1-\Psi^{\alpha~\nu}_{j~m}(x)\right] +
\P^{\alpha~\nu}_{j~m}(x) \Psi^{\alpha~\nu}_{j~m}(x)
\\\\ \ds~~~~~~~~~~~~~~~~~~~~~~~~~~~~~
~=~\sum_{\nu\colon\xi^\nu_j\in\mathds{S}^{n-1}} \P^{\alpha~\nu}_{j~m}(x) 
~=~\P^\alpha_{j~m}(x)
\end{array}
\]
and therefore 
\bel{I=II+III m,j}
\I_{j~m}^\alpha f(x)~=~\II_{j~m}^\alpha f(x)+\III_{j~m}^\alpha f(x).
\eeq
 In section 4, we will explicitly define
\[
 \I^{\alpha~z}_{j~m} f(x)~=~\II^{\alpha~z}_{j~m} f(x)+
 \III^{\alpha~z}_{j~m} f(x),\qquad \hbox{\small{$0\leq\Re z\leq1$}}
\]
of which $\II^{\alpha~z}_{j~m}$ and  $\III^{\alpha~z}_{j~m}$
 are two families of analytic operators.
    In particular, we have
 $\II^{\alpha~{1\over n}}_{j~m}f(x)=\II^{\alpha}_{j~m}f(x)$ and $ \III^{\alpha~{1\over n}}_{j~m}f(x)=\III^{\alpha}_{j~m}f(x)$.

$\bullet$ Given $j>0$, $m=1,\ldots,M$, we have
\bel{II 0 bad}
 \left\| \II^{\alpha~0+\i\Im z}_{j~m} f\right\|_{\L^2(\R^n)} ~\leq~\C_{\Re\alpha}~e^{\c|\Im\alpha|}~2^{-(1-\epsilon)j}2^{j/2}2^{-\ve j}~\left\| f\right\|_{\L^2(\R^n)}
\eeq
and
\bel{III 0 good}
 \left\| \III^{\alpha~0+\i\Im z}_{j~m} f\right\|_{\L^2(\R^n)} ~\leq~\C_{\Re\alpha}~e^{\c|\Im\alpha|}~2^{-(1-\epsilon)j}2^{-\ve j}~\left\| f\right\|_{\L^2(\R^n)}
\eeq
for some $\ve=\ve(\Re\alpha)>0$.

Observe that at $\Re z=0$, $\III_{j~m}^{\alpha~z}$ is a {\it good} partial operator of $\I^{\alpha~z}_{j~m}$ whose $\L^2$-norm is bound by $\C_{\Re\alpha}e^{\c|\Im\alpha|}2^{-(1-\epsilon)j}2^{-\ve j}$. On the other hand,  
$\II^{\alpha~z}_{j~m}$ is the {\it bad} one, satisfying (\ref{II 0 bad}). 

$\bullet$ Given $j>0$, $m=1,\ldots,M$, we find 
\bel{II 1 good}
 \left\| \II^{\alpha~1+\i\Im z}_{j~m} f\right\|_{\L^p(\R^n)} ~\leq~\C_{N~\Re\alpha}~e^{\c|\Im\alpha|}~2^{-(1-\epsilon)j}2^{-N j}~2^{-\ve j}~\left\| f\right\|_{\L^p(\R^n)},\qquad  1\leq p\leq\infty
\eeq
and
\bel{III 1 bad}
 \left\| \III^{\alpha~1+\i\Im z}_{j~m} f\right\|_{\L^p(\R^n)} ~\leq~\C_{\Re\alpha}~e^{\c|\Im\alpha|}~2^{-(1-\epsilon)j}2^{-\ve j}~\left\| f\right\|_{\L^p(\R^n)},\qquad  1\leq p\leq\infty
\eeq
for every $N\ge0$ and some $\ve=\ve(\Re\alpha)>0$.

At $\Re z=1$, $\III^{\alpha~z}_{j~m}$ is a {\it bad} partial operator of $\I^{\alpha~z}_{j~m}$. The best can be concluded is  that it satisfies (\ref{III 1 bad}). On the other hand,  
$\II^{\alpha~z}_{j~m}$ becomes the {\it good} one whose $\L^p$-norm is bounded by $\C_{N~\Re\alpha}e^{\c|\Im\alpha|}2^{-(1-\epsilon)j}2^{-N j}2^{-\ve j}$  for every $N\ge0$.

Recall $\II^{\alpha~{1\over n}}_{j~m}=\II^{\alpha}_{j~m}$ and $\III^{\alpha~{1\over n}}_{j~m}=\III^{\alpha}_{j~m}$. Take into account for ${n-1\over 2n}=0\left({1\over n}\right)+{1\over 2}\left({n-1\over n}\right)$ and ${n+1\over 2n}=1\left({1\over n}\right)+{1\over 2}\left({n-1\over n}\right)$.

From (\ref{II 0 bad}) and (\ref{II 1 good}) with $N$ chosen sufficiently large, we have
\bel{II EST}
 \begin{array}{cc}\ds
  \left\| \II^{\alpha}_{j~m} f\right\|_{\L^p(\R^n)} ~\leq~\C_{\Re\alpha}~e^{\c|\Im\alpha|}~2^{-(1-\epsilon)j}2^{-\ve j}~\left\| f\right\|_{\L^p(\R^n)}, 
  \\\\ \ds
    {n-1\over 2n}~\leq ~{1\over p}~\leq~ {n+1\over 2n}  
 \end{array}
\eeq
by applying Stein interpolation theorem \cite{Stein'}. 
 
From (\ref{III 0 good}) and (\ref{III 1 bad}), Stein interpolation theorem \cite{Stein'} again implies
\bel{III EST}
 \begin{array}{cc}\ds
  \left\| \III^{\alpha}_{j~m} f\right\|_{\L^p(\R^n)} ~\leq~\C_{\Re\alpha}~e^{\c|\Im\alpha|}~2^{-(1-\epsilon)j}2^{-\ve j}~\left\| f\right\|_{\L^p(\R^n)}, 
  \\\\ \ds
    {n-1\over 2n}~\leq~ {1\over p}~\leq~ {n+1\over 2n}. 
 \end{array}
\eeq
Recall (\ref{I=II+III m,j}).
By putting together (\ref{II EST})-(\ref{III EST}), we obtain
(\ref{Result Two j m}) as desired.

The remaining paper is organized as follows. In section 2, we show  {\bf Theorem Two} implying {\bf Theorem One}. In section 3, after a formal derivation of $\I^{\alpha}_{j~m}=\II^{\alpha}_{j~m}+\III^{\alpha}_{j~m}$, we introduce two end-point estimates associated with $\Re z=0$ and $\Re z=1$ in (\ref{II 0 bad})-(\ref{III 1 bad}), stated as {\bf Proposition One} and {\bf Proposition Two}. By using these results, we finish the proof of {\bf Theorem Two} in section 4. We prove {\bf Proposition One} in section 5. 
Section 6 is devoted to certain preliminaries regarding {\bf Proposition Two}. We prove {\bf Proposition Two} in section 7.
Two appendices added in the end for the sake of self-containedness.

\section{Theorem Two implies Theorem One}
\setcounter{equation}{0}
Let $0<\Re\q<{1\over 2}$. 
Consider
  \bel{S_phi}
  \begin{array}{lr}\ds
  \S^\q_\psi f(x)~=~\int_{\R^n} e^{2\pi\i x\cdot\xi} \Hat{f}(\xi) \Hat{\psi}(\xi)\left(1-|\xi|^2\right)^\q_+ d\xi
 \end{array}  
  \eeq  
where $ \Hat{\psi}(\xi)=\left[\Hat{\phi}(\xi)\right]^2$ and $\Hat{\phi}$ is defined in (\ref{hat phi}).

Observe that $\left[1-\Hat{\psi}(\xi)\right] \left(1-|\xi|^2\right)^\q_+$
is a $\L^p$-Fourier multiplier for $1<p<\infty$. In order to prove {\bf Theorem One},
 it is suffice to show
   \bel{Result}
  \begin{array}{cc}\ds
 \left\|  \S^\q_\psi f\right\|_{\L^p(\R^n)}~\leq~\C_{\Re\q}~e^{\c|\Im\q|}~\left\| f\right\|_{\L^p(\R^n)},\qquad
  {n-1\over 2n}~\leq ~{1\over p}~\leq~ {n+1\over 2n}.
 \end{array}
  \eeq 
Let $1-\q=\alpha\in\Cx$. We define
\bel{m+}
\begin{array}{lr}\ds
\m^\alpha_+(\xi)~=~\Hat{\phi}(\xi)\int_0^1  \left({1\over \tau^2-|\xi|^2}\right)^\alpha_+ \tau d\tau
\\\\ \ds~~~~~~~~~~
~=~\Hat{\phi}(\xi)\left\{ \begin{array}{lr}\ds
\int_{|\xi|}^1 \left({1\over \tau^2-|\xi|^2}\right)^\alpha \tau d\tau,\qquad |\xi|<1
\\ \ds~~~~~~~~~~~~~~~~~~
0\qquad~~~~~~~~~~~~~~~ |\xi|\ge1
\end{array}\right.
\\\\ \ds~~~~~~~~~~
~=~{1\over 2} (1-\alpha)^{-1} \Hat{\phi}(\xi) \left(1-|\xi|^2\right)_+^{1-\alpha}
  \end{array}
  \eeq
and
\bel{m-} 
\begin{array}{lr}\ds
\m^\alpha_-(\xi)~=~\Hat{\phi}(\xi)\int_0^1  \left({1\over \tau^2-|\xi|^2}\right)^\alpha_- \tau d\tau 
\\\\ \ds~~~~~~~~~~
~=~(-1)^{-\alpha}\Hat{\phi}(\xi)\left\{ \begin{array}{lr}\ds
\int_0^{|\xi|} \left({1\over \tau^2-|\xi|^2}\right)^\alpha\tau d\tau,\qquad |\xi|\leq1
\\ \ds
\int_0^1 \left({1\over \tau^2-|\xi|^2}\right)^\alpha\tau d\tau,\qquad~ |\xi|>1
\end{array}\right.
\\\\ \ds~~~~~~~~~~
 ~=~{1\over 2}(1-\alpha)^{-1}  \Hat{\phi}(\xi) \left[|\xi|^{2(1-\alpha)}-\left(1-|\xi|^2\right)_-^{1-\alpha}\right].
\end{array}
\eeq
Recall $\Hat{\Lambda}^\alpha(\xi,\tau) $ given at (\ref{Lambda Transform}).
From (\ref{m+})-(\ref{m-}), we find
 \bel{Lambda m}
 \begin{array}{lr}\ds
  \Hat{\phi}(\xi)\int_0^1  \Hat{\Lambda}^\alpha(\xi,\tau) \tau d\tau
  ~=~\pi^{{n-1\over 2}-2\alpha}\Gamma(\alpha) \Bigg\{ \m^\alpha_-(\xi)-\sin\pi\left(\alpha-{1\over 2}\right)\m^\alpha_+(\xi)\Bigg\}.
  \end{array}
 \eeq 
Moreover, define
 \bel{m^alpha}
 \m^\alpha(\xi)~=~\Hat{\phi}(\xi) \Bigg\{  -\left(1-|\xi|^2\right)_-^{1-\alpha}-\sin\pi\left(\alpha-{1\over 2}\right) \left(1-|\xi|^2\right)_+^{1-\alpha}\Bigg\}.
 \eeq 
 From (\ref{m+})-(\ref{m-}), (\ref{Lambda m}) and (\ref{m^alpha}), we have
  \bel{m alpha}
 \begin{array}{lr}\ds
  \Hat{\phi}(\xi)\int_0^1  \Hat{\Lambda}^\alpha(\xi,\tau) \tau d\tau
  ~=~\pi^{{n-2\over 2}-2\alpha}{1\over 2}(1-\alpha)^{-1} \Gamma(\alpha)\Big[  \m^\alpha(\xi)+\Hat{\phi}(\xi)  |\xi|^{2(1-\alpha)} \Big].
  \end{array}
 \eeq 
   \begin{remark} Let 
   \[ \T_zf(x)~=~ \int_{\R^n} e^{2\pi\i x\cdot\xi}\Hat{f}(\xi)\m^z(\xi)d\xi,\qquad z\in\Cx.\]
   We call $\m^z(\xi)$  a restricted $\L^p$-Fourier multiplier if 
\[\left\| \T_z f\right\|_{\L^p(\R^n)}~\leq~\C_{\Re z}~e^{\c|\Im z|}~\left\| f\right\|_{\L^p(\R^n)},\qquad {n-1\over 2n}~\leq~{1\over p}~\leq~{n+1\over 2n}.\]
\end{remark}  
Recall $\I^\alpha$ defined in (\ref{I alpha}). Let ${1\over 2}<\Re\alpha<1$.  {\bf Theorem Two} states that $\I^\alpha f$ satisfies  the $\L^p$-norm inequality in (\ref{Result Two}). Therefore, the left-hand-side of the equation in (\ref{m alpha}) is a restricted $\L^p$-Fourier multiplier.
Because $\Hat{\phi}(\xi)|\xi|^{2(1-\alpha)}$  is  a $\L^p$-Fourier multiplier for $1<p<\infty$, we find $\m^\alpha(\xi)$  as a restricted $\L^p$-Fourier multiplier.

Consider
\[
\Hat{\phi}(\xi)\m^\alpha(\xi)~=~ \Hat{\psi}(\xi) \Bigg\{  -\left(1-|\xi|^2\right)_-^{1-\alpha}-\sin\pi\left(\alpha-{1\over 2}\right) \left(1-|\xi|^2\right)_+^{1-\alpha}\Bigg\}
\]
and
\[
\m^{{1\over 2}+{\alpha\over 2}}(\xi)~=~\Hat{\phi}(\xi) \Bigg\{  -\left(1-|\xi|^2\right)_-^{{1\over 2}-{\alpha\over 2}}-\sin\pi\left({\alpha\over 2}\right) \left(1-|\xi|^2\right)_+^{{1\over 2}-{\alpha\over 2}}\Bigg\}.
\]
  Clearly, 
 both $\Hat{\phi}(\xi)\m^\alpha(\xi)$ and $\m^{{1\over 2}+{\alpha\over 2}}(\xi)$ 
 are  restricted $\L^p$-Fourier multipliers. 
 Furthermore, 
\bel{m_3}
\begin{array}{lr}\ds
\left[\m^{{1\over 2}+{\alpha\over 2}}(\xi)\right]^2
~=~ \Hat{\psi}(\xi)  \Bigg\{  \left(1-|\xi|^2\right)_-^{1-\alpha}+\sin^2\pi\left({\alpha\over 2}\right)\left(1-|\xi|^2\right)_+^{1-\alpha}\Bigg\}  
 \end{array}
\eeq 
is indeed another restricted $\L^p$-Fourier multiplier.

By adding $\left[\m^{{1\over 2}+{\alpha\over 2}}(\xi)\right]^2$ and $\Hat{\phi}(\xi)\m^\alpha(\xi)$ together, we obtain
\bel{subtract}
 \begin{array}{lr}\ds
\Hat{\psi}(\xi) \Bigg[\sin^2\pi\left({\alpha\over 2}\right)-\sin\pi\left(\alpha-{1\over 2}\right)\Bigg] \left(1-|\xi|^2\right)_+^{1-\alpha}.
\end{array}
\eeq 
A direct computation shows
\[
\begin{array}{lr}\ds
\sin^2\pi\left({\alpha\over 2}\right)~=~{1\over 2}-{1\over 4} \left[ e^{-\pi\Im\alpha}+e^{\pi\Im\alpha}\right]\cos \pi\Re\alpha-{\i\over 4} \left[ e^{-\pi\Im\alpha}-e^{\pi\Im\alpha}\right]\sin\pi\Re\alpha,
\\\\ \ds
\sin\pi\left(\alpha-{1\over 2}\right)~=~-{1\over 2} \left[ e^{-\pi\Im\alpha}+e^{\pi\Im\alpha}\right]\cos \pi\Re\alpha-{\i\over 2} \left[ e^{-\pi\Im\alpha}-e^{\pi\Im\alpha}\right]\sin\pi\Re\alpha.
\end{array}
\]
Consequently, we have
\bel{non-zero}
\sin^2\pi\left({\alpha\over 2}\right)-\sin\pi\left(\alpha-{1\over 2}\right)\neq0,\qquad \hbox{\small{$0<\Re\alpha<1$}}.
\eeq
From (\ref{subtract}) and (\ref{non-zero}), we conclude that
 $\Hat{\psi}(\xi)\left(1-|\xi|^2\right)_+^{1-\alpha}$ is a restricted $\L^p$-Fourier multiplier as desired.

\section{End-point estimates regarding $\II^{\alpha}_{j~m}$ and $\III^{\alpha}_{j~m}$}
\setcounter{equation}{0}
Recall $\I^\alpha$ defined in (\ref{I alpha}). 
First, we show $\I^\alpha f$ for ${1\over 2}<\Re\alpha<1$ can be expressed as (\ref{I alpha new express}).

$\Omega^\alpha$ is a distribution defined in $\R^n$ by analytic continuation from
\[
\begin{array}{cc}\ds
\hbox{\small{$\Re\alpha>{n-1\over 2}$}},\qquad \Omega^{\alpha}(x)~=~\hbox{\small{$\pi^{\alpha-{n+1\over 2}}\Gamma^{-1}\left(\alpha-{n-1\over 2}\right)$}} 
 \left({1 \over1-|x|^2}\right)^{{n+1\over 2}-\alpha}_+.
 \end{array}
\]
Equivalently,  it can be  defined by
\bel{Omega^alpha Transform}
\begin{array}{lr}\ds
\Hat{\Omega}^{\alpha}(\xi)~=~\left({1\over|\xi|}\right)^{{n\over 2}-\big[{n+1\over 2}-\alpha\big]} \J_{{n\over 2}-\big[{n+1\over 2}-\alpha\big]}\Big(2\pi|\xi|\Big)
\\\\ \ds~~~~~~~~~~
~=~
\left({1\over|\xi|}\right)^{\alpha-{1\over 2}} \J_{\alpha-{1\over 2}}\Big(2\pi|\xi|\Big),\qquad \alpha\in\Cx.
\end{array}
\eeq
See (\ref{Omega^z})-(\ref{Omega^z Transform}). 
By using the integral formula of Bessel functions in (\ref{Bessel}), we  find
 \bel{r local inte}
 \Hat{\Omega}^\alpha (\xi)~=~ \pi^{\alpha-1}\Gamma\left(\alpha\right)\int_{-1}^1 e^{2\pi\i |\xi|s } (1-s^2)^{\alpha-1} ds.
  \eeq 
On the other hand, $\Lambda^\alpha$ is a distribution defined in $\R^{n+1}$ by analytic continuation from
\bel{Lambda^alpha}
\begin{array}{cc}\ds
\hbox{\small{$\Re\alpha>{n-1\over2}$}},
\qquad
 \Lambda^\alpha(x,t)~=~ \hbox{\small{$\pi^{\alpha-{n+1\over2}}\Gamma^{-1}\left(\alpha-{n-1\over 2}\right)$}}
  \left({1\over t^2-|x|^2}\right)^{{n+1\over 2}-\alpha}_+ .
  \end{array}
\eeq
Let $h$ be a Schwartz function defined in $\R^{n+1}$. 
We have
\bel{Kernel singular}
 \begin{array}{lr}\ds
 h\ast\Lambda^\alpha(x,t)~=~\hbox{\small{$\pi^{\alpha-{n+1\over2}}\Gamma^{-1}\left(\alpha-{n-1\over 2}\right)$}}\iint_{|u|<|r|} h(x-u,t-r) 
  \left({1\over r^2-|u|^2}\right)^{{n+1\over 2}-\alpha}dudr 
 \\\\ \ds~~~~~~~~~~~~~~~~~~
 ~=~\iint_{\R^{n+1}} h(x-u,t-r)\Omega^\alpha\left({u\over r}\right)|r|^{2\alpha-1-n}dudr
  \\\\ \ds~~~~~~~~~~~~~~~~~~
 ~=~\int_\R\left\{\int_{\R^n} e^{2\pi\i x\cdot\xi}\Hat{h}(\xi,t-r) \Hat{\Omega}^\alpha (r\xi)d\xi \right\} |r|^{2\alpha-1}   dr
 \end{array}
\eeq
whenever $\Re\alpha>{n-1\over 2}$. Here, $\Hat{h}(\cdot,t-r)$ is the Fourier transform of $h(x,t-r)$ in $x\in\R^n$.
 \begin{remark}
By the principle of analytic continuation, we must have  (\ref{Kernel singular}) hold for every $\Re\alpha>0$.
\end{remark} 
 Let ${1\over 2}<\Re\alpha<1$. From (\ref{Kernel singular}), we find
\[ h\ast\Lambda^\alpha(x,t)~=~\lim_{N\mt\infty} \iint_{\R^{n+1}} e^{2\pi\i\big[x\cdot\xi+t\tau\big]} \Hat{h}(\xi,\tau) \left\{\int_{-N}^N  e^{-2\pi\i \tau r} \Hat{\Omega}^\alpha (r\xi) |r|^{2\alpha-1} dr\right\} d\xi d\tau.
\]
By using the asymptotic expansion of Bessel functions in (\ref{J asymptotic})-(\ref{J O}), together with (\ref{r local inte}), we write
\bel{Lambda transform Rewrite}
\begin{array}{lr}\ds
\int_\R e^{-2\pi\i \tau r}\Hat{\Omega}^\alpha (r\xi) |r|^{2\alpha-1} dr
~=~\int_\R e^{-2\pi\i \tau r}\left({1\over|r\xi|}\right)^{\alpha-{1\over 2}} \J_{\alpha-{1\over 2}}\Big(2\pi|r\xi|\Big) |r|^{2\alpha-1} dr\qquad\hbox{\small{by (\ref{Omega^alpha Transform})}}
\\\\ \ds~~~~~~~~~~~~~~~~~~~~~~~~~~~~~~~~~~~~~~~~~~
~=~{1\over \pi}\int_\R e^{-2\pi\i \tau r}\left({1\over|r\xi|}\right)^\alpha \cos\left[2\pi|r\xi|-{\pi\over 2}\alpha\right] |r|^{2\alpha-1} dr
\\\\ \ds~~~~~~~~~~~~~~~~~~~~~~~~~~~~~~~~~~~~~~~~~~
~+~\int_\R e^{-2\pi\i \tau r}\mathcal{E}^\alpha\left(2\pi |r\xi|\right) |r|^{2\alpha-1} dr,
\\\\ \ds~~~~~~~~~~~~~~~~~~~~~~~~~~~~~~~~~~~
\left|\mathcal{E}^\alpha(\rho)\right|~\leq~\C_{\Re\alpha}~e^{\c|\Im\alpha|}\left\{\begin{array}{lr}\ds \rho^{-\Re\alpha},\qquad 0<\rho\leq1,
 \\ \ds
 \rho^{-\Re\alpha-1},\qquad \rho>1.
 \end{array}\right.
 \end{array}
\eeq
Consider ${1\over 10}<|\xi|\leq10$. The norm estimate for $\mathcal{E}^\alpha$ implies 
 \[ \left| \int_\R e^{-2\pi\i \tau r}\mathcal{E}^\alpha\left(2\pi |r\xi|\right) |r|^{2\alpha-1} dr\right| 
 ~\leq~\C_{\Re\alpha}~e^{\c|\Im\alpha|}.\]
Because of Euler's formulae, we replace the cosine function in (\ref{Lambda transform Rewrite}) with $e^{-\i\left({\pi\over 2}\right)\alpha} e^{2\pi\i |r\xi|}$ or $e^{\i\left({\pi\over 2}\right)\alpha} e^{-2\pi\i |r\xi|}$.  
By  integration by parts $w.r.t~r$, we find
 \bel{main term norm r}
 \begin{array}{lr}\ds
 \left| \int_{2^{j-1}}^{2^j} e^{2\pi\i \big[|\xi|-\tau\big] r}\left({1\over|r\xi|}\right)^\alpha  |r|^{2\alpha-1} dr\right|
 \\\\ \ds
\leq~ \C_{\Re\alpha}~e^{\c|\Im\alpha|}  \left\{\begin{array}{lr} 2^{\Re\alpha j} ,~~~~~~~~~~~~~~~~~~~~~~\qquad 2^j\leq \big||\tau|-|\xi|\big|^{-1},
 \\ \ds
 \left|{1\over |\tau|-|\xi|}\right| 2^{\big[\Re\alpha-1\big]j},\qquad 2^j> \big||\tau|-|\xi|\big|^{-1}. 
 \end{array}\right.
 \end{array}
 \eeq
For symmetry reason,  (\ref{main term norm r}) further implies
 \bel{Int Norm Est r}
 \begin{array}{lr}\ds
 \left| \int_\R e^{-2\pi\i \tau r}\left({1\over|r\xi|}\right)^\alpha \cos\left[2\pi|r\xi|-{\pi\over 2}\alpha\right] |r|^{2\alpha-1} dr\right|
 \\\\ \ds
 \leq~  \C_{\Re\alpha}~e^{\c|\Im\alpha|} \sum_{2^j\leq \left||\tau|-|\xi|\right|^{-1}} 2^{\Re\alpha j}+\C_{\Re\alpha}~e^{\c|\Im\alpha|} \sum_{2^j> \left||\tau|-|\xi|\right|^{-1}}  \left|{1\over |\tau|-|\xi|}\right| 2^{\big[\Re\alpha-1\big]j}
 \\\\ \ds
  \leq~  \C_{\Re\alpha}~e^{\c|\Im\alpha|}~ \left|{1\over |\tau|-|\xi|}\right|^{\Re\alpha}.
  \end{array}  
 \eeq
 On the other hand, the Fourier transform of $\Lambda^\alpha$ defined by analytic continuation from (\ref{Lambda^alpha}) agrees with the function
\[
\begin{array}{lr}
  \Hat{\Lambda}^\alpha(\xi,\tau)~=~
\pi^{{n-1\over 2}-2\alpha}\Gamma\left(\alpha\right)
 \left\{ ~{\ds \left({1\over \tau^2-|\xi|^2}\right)^\alpha_-} - \sin\pi\left(\alpha-{1\over 2}\right){\ds \left( {1\over \tau^2- |\xi|^2}\right)^\alpha_+}~\right\}
  \end{array}
\]
whenever $|\tau|\neq|\xi|$.  See appendix B.

From (\ref{Kernel singular}) to (\ref{Int Norm Est r}), by applying Lebesgue's dominated convergence theorem, we have
  \bel{Kernel singular region}
 \begin{array}{lr}\ds
 h\ast\Lambda^\alpha(0,0)
 ~=~\lim_{N\mt\infty}\iint_{\R^{n+1}} \Hat{h}(\xi,\tau) 
 \left\{ \int_{-N}^N e^{-2\pi\i \tau r}\Hat{\Omega}^\alpha (r\xi) |r|^{2\alpha-1} dr\right\} d\xi d\tau 
 \\\\ \ds~~~~~~~~~~~~~~~~~~
 ~=~\iint_{\R^{n+1}} \Hat{h}(\xi,\tau) 
 \left\{ \int_{\R} e^{-2\pi\i \tau r}\Hat{\Omega}^\alpha (r\xi) |r|^{2\alpha-1} dr\right\} d\xi d\tau  
 \\\\ \ds~~~~~~~~~~~~~~~~~~
 ~=~\iint_{\R^{n+1}} \Hat{h}(\xi,\tau) 
\Hat{\Lambda}^\alpha(\xi,\tau) d\xi d\tau 
 \end{array}
\eeq
 for every  $h$ defined in $\R^{n+1}$ whose Fourier transform is supported in the cylindrical region: $\left\{(\xi,\tau)\in\R^n\times\R\colon {1\over 10}<|\xi|\leq10\right\}$. 
 
 Let $\varphi\in\mathcal{C}^\infty_o(\R)$   such that $\varphi(t)=1$ if $|t|\leq1$ and $\varphi(t)=0$ if $|t|>2$.  Denote $\ds\c_\varphi^{-1}=\int_\R \varphi(t)dt$.  
Given $(\eta,s)\in\R^n\times\R$ for $|s|\neq|\eta|$ and ${1\over 3}\leq|\eta|\leq3$, we consider a family of {\it good kernels}: 
$\ds\Hat{h}_\ve(\xi,\tau)=\c_\varphi  \ve^{-(n+1)}\varphi \left[\ve^{-1}\sqrt{|\xi-\eta|^2+(\tau-s)^2}\right],~ 0<\ve<{1/ 10}$.

Replace $\Hat{h}(\xi,\tau)$ by $\Hat{h}_\ve(\xi,\tau)$ in (\ref{Kernel singular region}). By taking $\ve\mt0$, we find
\bel{Lambda transform rewrite}
\begin{array}{lr}\ds
\Hat{\Lambda}^\alpha(\xi,\tau)~=~\int_\R e^{-2\pi\i \tau r}\Hat{\Omega}^\alpha (r\xi) |r|^{2\alpha-1} dr,\qquad\hbox{\small{$|\tau|\neq|\xi|$}},\qquad \hbox{\small{${1\over 3}<|\xi|\leq3$}}.
\end{array} 
 \eeq
By applying Lebesgue's dominated convergence theorem again, we have 
 \bel{I alpha rewrite}
 \begin{array}{lr}\ds
~~~~~~~~~~~~~~~~~~~~~~~~~~~~~~~~~~~~~~~~~~~~~ \Hat{\I^\alpha f}(\xi)~=~\Hat{f}(\xi)\Hat{\P}^\alpha(\xi),
 \\\\ \ds 
\Hat{\P}^\alpha(\xi)~=~\Hat{\phi}(\xi)\int_0^1   \Hat{\Lambda}^\alpha(\xi,\tau)\tau d\tau  
   \\\\ \ds~~~~~~~~~
  ~=~ \Hat{\phi}(\xi)\int_{0<\tau<1,~\tau\neq|\xi|}  \left\{\int_\R e^{-2\pi\i \tau r}\Hat{\Omega}^\alpha (r\xi) |r|^{2\alpha-1} dr\right\}\tau d\tau \qquad\hbox{\small{by (\ref{Lambda transform rewrite})}}  
   \\\\ \ds~~~~~~~~~
  ~=~\Hat{\phi}(\xi)\lim_{N\mt\infty}\int_{0<\tau<1,~\tau\neq|\xi|} \left\{ \int_{-N}^N e^{-2\pi\i \tau r}\Hat{\Omega}^\alpha (r\xi) |r|^{2\alpha-1} dr\right\} \tau d\tau
  \\\\ \ds~~~~~~~~~
~=~\Hat{\phi}(\xi) \int_\R \left\{\int_0^1 e^{-2\pi\i \tau r} \tau d\tau\right\} \Hat{\Omega}^\alpha (r\xi) |r|^{2\alpha-1} dr
\\\\ \ds~~~~~~~~~
~=~\Hat{\phi}(\xi)\int_{\R} e^{-2\pi\i  r}\Hat{\Omega}^\alpha (r\xi)\omega(r) |r|^{2\alpha-1} dr,
\\\\ \ds
 \omega(r)~=~e^{2\pi\i r}\int_0^1 e^{-2\pi\i \tau r} \tau d\tau~=~{-1\over 2\pi\i}{1\over r}-{1\over 4\pi^2r^2} \left[1-e^{-2\pi\i r}\right].
\end{array}
 \eeq
Note that $\Hat{\phi}$ defined in (\ref{hat phi}) has 
$\supp\Hat{\phi}\subset\left\{\xi\in\R^n\colon{1\over 3}<|\xi|\leq3\right\}$.

Consider
\bel{P sum j} 
\begin{array}{cc}\ds
\Hat{\P}^\alpha(\xi)
~=~\Hat{\P}^\alpha_<(\xi)+\sum_{j>0}~\Hat{\P}_j^\alpha(\xi),
\\\\ \ds
\Hat{\P}^\alpha_<(\xi)~=~\Hat{\phi}(\xi)\int_{-1}^1 e^{-2\pi\i  r}\Hat{\Omega}^\alpha (r\xi)\omega(r) |r|^{2\alpha-1} dr,
\\\\ \ds
\Hat{\P}_j^\alpha(\xi)~=~\Hat{\phi}(\xi)\int_{2^{j-1}\leq|r|<2^j} e^{-2\pi\i  r}\Hat{\Omega}^\alpha (r\xi)\omega(r) |r|^{2\alpha-1} dr,\qquad j>0.
\end{array}
\eeq 
Furthermore, we assert
\bel{P j m}
\begin{array}{cc}\ds
\Hat{\P}_j^\alpha(\xi)~=~\sum_{m=1}^M \Hat{\P}_{j~m}^\alpha(\xi),
\\\\ \ds
\Hat{\P}_{j~m}^\alpha(\xi)~=~\Hat{\phi}(\xi)\int_{\lambda_{m-1}\leq|r|<\lambda_m} e^{-2\pi\i  r}\Hat{\Omega}^\alpha (r\xi)\omega(r) |r|^{2\alpha-1} dr,
\\\\ \ds
\hbox{$\lambda_m\in[2^{j-1},2^j]$}, 
\qquad
 \hbox{$\lambda_0=2^{j-1}$, $\lambda_M=2^j$ ~{\small and} ~$2^{\epsilon j-1}~\leq~\lambda_m-\lambda_{m-1}<2^{\epsilon j}$}
\end{array}
\eeq
where $0<\sigma=\sigma(\Re\alpha)<{1\over 2}$ can be chosen sufficiently small.

By using (\ref{r local inte}), we find
\bel{kernel <}
 \begin{array}{lr}\ds
 \P^\alpha_<(x)~=~ \int_{\R^n}e^{2\pi\i x\cdot\xi} \Hat{\phi}(\xi)\left\{\int_{-1}^1 e^{-2\pi\i  r}\Hat{\Omega}^\alpha (r\xi)\omega(r) |r|^{2\alpha-1} dr\right\} d\xi
 \\\\ \ds~~~~~~~~~
~=~ \pi^{\alpha-1}\Gamma\left(\alpha\right)
\\ \ds~~~~~~~~~~~~~
\int_{\R^n}e^{2\pi\i x\cdot\xi} \Hat{\phi}(\xi)  \left\{\int_{-1}^1 e^{-2\pi\i  r}\left\{\int_{-1}^1 e^{2\pi\i |r\xi|s } (1-s^2)^{\alpha-1} ds\right\}\omega(r) |r|^{2\alpha-1} dr\right\} d\xi.
\end{array}
 \eeq 
An $N$-fold integration by parts $w.r.t~\xi$ inside (\ref{kernel <}) shows 
$\left|\P^\alpha_<(x)\right|\leq\C_{N~\Re\alpha}e^{\c|\Im\alpha|}\left({1\over 1+|x|}\right)^N$.
Hence that $\P^\alpha_<$ is an $\L^1$-function in $\R^n$.
\v

Given $j>0$,  $\left\{ \xi^\nu_j\right\}_\nu$ is a collection of points  almost equally distributed on  $\mathds{S}^{n-1}$ having a grid length between $2^{-j/2-1}$ and $2^{-j/2}$:   {\bf (1)} $|\xi^\mu_j-\xi^\nu_j|\ge 2^{-j/2-1}$ for every $\xi^\mu_j\neq\xi^\nu_j$. {\bf (2)} For any $u\in \mathds{S}^{n-1}$, there is a $\xi^\nu_j$ in the open set $\{\xi\in\mathds{S}^{n-1}\colon |\xi-u|<2^{-j/2+1}\}$.

Let $\varphi\in\mathcal{C}^\infty_o(\R)$ be a smooth cut-off function  such that $\varphi(t)=1$ if $|t|\leq1$ and $\varphi(t)=0$ if $|t|>2$. Recall (\ref{phi^v_j intro})-(\ref{Gamma_j}). We have
\bel{phi^v_j}
\begin{array}{cc}\ds
\varphi^\nu_j(\xi)~=~{\ds\varphi\Bigg[2^{j/2}\left|{\xi\over |\xi|}-\xi^{\nu}_{j}\right|\Bigg]\over \ds \sum_{\nu\colon\xi^\nu_j\in\mathds{S}^{n-1}} \varphi\Bigg[2^{j/2}\left|{\xi\over |\xi|}-\xi^{\nu}_{j}\right|\Bigg]},
\\\\ \ds
\supp \varphi^\nu_j~\subset~\Gamma^\nu_j~=~\Bigg\{\xi\in\R^n\colon \left|{\xi\over |\xi|}-\xi^{\nu}_{j}\right|<2^{-j/2+1}\Bigg\}.
\end{array}
\eeq

Let
 \bel{P j v}
 \begin{array}{cc}\ds
 \Hat{\P}_{j~m}^{\alpha~\nu}(\xi)~=~\varphi^\nu_j(\xi) \Hat{\phi}(\xi)\int_{\lambda_{m-1}\leq|r|<\lambda_m} e^{-2\pi\i  r}\Hat{\Omega}^\alpha (r\xi)\omega(r) |r|^{2\alpha-1} dr. 
 \end{array}
  \eeq
We write
   \bel{P j m Sum}
    \Hat{\P}_{j~m}^\alpha(\xi)~=~ 
\sum_{\nu\colon \xi^\nu_j\in\mathds{S}^{n-1}}  \Hat{\P}_{j~m}^{\alpha~\nu}(\xi)
\eeq
and   
  \bel{I h j m Sum}
\begin{array}{cc}\ds
\I^{\alpha}_{j~m} f(x)~=~\int_{\R^n} e^{2\pi\i x\cdot\xi} \Hat{f}(\xi)  \Hat{\P}_{j~m}^{\alpha}(\xi)d\xi
\end{array}
\eeq  
for every $j>0$ and $m=1,2,\ldots,M$.

 Given $\xi_j^\nu\in\mathds{S}^{n-1}$, $\L_\nu$ is an $n\times n$-orthogonal matrix with $\det\L_\nu=1$. Moreover, we require
$\L_\nu^T \xi^\nu_j=(1,0)^T\in\R\times\R^{n-1}$. Recall (\ref{Psi mu j alpha})-(\ref{Rectangles}). We have
 \bel{Psi mu j}
 \begin{array}{cc}\ds
 \Psi^{\alpha}_{j~m}(x)~=~\sum_{\nu\colon\xi^\nu_j\in\mathds{S}^{n-1}}  \Psi^{\alpha~\nu}_{j~m}(x),
 \\\\ \ds 
  \Psi^{\alpha~\nu}_{j~m}(x)~=~\varphi\left[ 2^{-\epsilon j-1}\left|\lambda_m-\left(\L_\nu^T x\right)_1\right|\right] 
\prod_{i=2}^n \varphi\left[2^{-\left[{1\over 2}+\epsilon\right]j}\left|\left(\L_\nu^T x\right)_i\right|\right]
\\ \ds~~~~~~~~~~~~~
~+~\varphi\left[ 2^{-\epsilon j-1}\left|\lambda_m+\left(\L_\nu^T x\right)_1\right|\right] 
\prod_{i=2}^n \varphi\left[2^{-\left[{1\over 2}+\epsilon\right]j}\left|\left(\L_\nu^T x\right)_i\right|\right].

 \\\\ \ds
 \supp  \Psi^{\alpha~\nu}_{j~m}~\subset~  \hbox{R}^\nu_{j~m}~=~
 \\\\ \ds
\Bigg\{x\in\R^n\colon \lambda_m-2^{\epsilon j+2}<\left(\L_\nu^T x\right)_1<\lambda_m+2^{\epsilon j+2},~\left|\left(\L_\nu^T x\right)_i\right|<2^{\left[{1\over 2}+\epsilon\right]j+1},~ i=2,\ldots,n
\Bigg\}~\cup
\\ \ds
~~\Bigg\{x\in\R^n\colon -\lambda_m-2^{\epsilon j+2}<\left(\L_\nu^T x\right)_1<-\lambda_m+2^{\epsilon j+2},~\left|\left(\L_\nu^T x\right)_i\right|<2^{\left[{1\over 2}+\epsilon\right]j+1},~ i=2,\ldots,n
\Bigg\}.
 \end{array}
  \eeq
Consider
\bel{I=II+III}
\I^{\alpha}_{j~m}f(x)~=~\II^{\alpha}_{j~m} f(x)+ \III^{\alpha}_{j~m} f(x)
\eeq
of which
\bel{II_j U again}
 \begin{array}{cc}\ds
  \II^{\alpha}_{j~m} f(x)~=~\int_{\R^n} f(x-u) \U^{\alpha}_{j~m}(u)du,
  \\\\ \ds
\U^{\alpha}_{j~m}(x)~=~\sum_{\nu\colon\xi^\nu_j\in\mathds{S}^{n-1}}\U^{\alpha~\nu}_{j~m}(x),\qquad \U^{\alpha~\nu}_{j~m}(x)~=~\P^{\alpha~\nu}_{j~m}(x) \left[1-\Psi^{\alpha~\nu}_{j~m}(x)\right] 
  
    \end{array}
  \eeq 
 and
 \bel{III_j V again}
 \begin{array}{cc}\ds
  \III^{\alpha}_{j~m} f(x)~=~\int_{\R^n} f(x-u) \V^{\alpha}_{j~m}(u)du, 
  \\\\ \ds
  \V^{\alpha}_{j~m} (x)~=~\sum_{\nu\colon\xi^\nu_j\in\mathds{S}^{n-1}} \V^{\alpha~\nu}_{j~m} (x),\qquad        \V^{\alpha~\nu}_{j~m} (x)~=~\Psi^{\alpha~\mu}_{j~m}(x)\P^{\alpha~\nu}_{j~m}(x)
 \end{array}
 \eeq 
for every $j>0$ and $m=1,2,\ldots,M$.

Next, we begin to introduce the two end-point estimates regarding $\II^\alpha_{j~m}$ and $\III^\alpha_{j~m}$.

\subsection{${^\sharp}\U^{\alpha~\beta}_{j~m}$ and ${^\sharp}\V^{\alpha~\beta}_{j~m}$}
${^\sharp}\Omega^\alpha$ is a distribution defined in $\R^n$ by analytic continuation from
\bel{Omega sharp}
\begin{array}{cc}\ds
 {^\sharp}\Omega^{\alpha}(x)~=~\pi^{\alpha-{n+2\over 2}}\hbox{$\Gamma^{-1}\left(\alpha-{n\over 2}\right)$} 
 \left({1 \over1-|x|^2}\right)^{{n+2\over 2}-\alpha}_+,\qquad \hbox{\small{$\Re\alpha>{n\over 2}$}}. 
 \end{array}
\eeq
Equivalently,  it can be  defined by
 \bel{Omega sharp Transform} 
\begin{array}{lr}\ds
{^\sharp}\Hat{\Omega}^{\alpha}(\xi)~=~
\left({1\over|\xi|}\right)^{{n\over 2}-\big[{n+2\over 2}-\alpha\big]} \J_{{n\over 2}-\big[{n+2\over 2}-\alpha\big]}\Big(2\pi|\xi|\Big)
\\\\ \ds~~~~~~~~~~~
~=~\left({1\over|\xi|}\right)^{\alpha-1} \J_{\alpha-1}\Big(2\pi|\xi|\Big),\qquad \alpha\in\Cx.
\end{array}
\eeq
Let $ \Re\alpha\ge\left[{2n\over 2n-1}\right]\Re\beta$ and ${2n-1\over 4n}<\Re\beta<{2n-1\over 2n-2}$. Define
\bel{P sharp j}
\begin{array}{cc}\ds
{^\sharp}\Hat{\P}_j^{\alpha~\beta}(\xi)~=~\Hat{\phi}(\xi)\int_{2^{j-1}\leq|r|<2^j} e^{-2\pi\i  r}{^\sharp}\Hat{\Omega}^\alpha (r\xi)\omega(r) |r|^{2\beta-1} dr,\qquad j>0.
\end{array}
\eeq
\begin{lemma} We have
\bel{Kernel T Est}
\begin{array}{lr}\ds
\sum_{j>0} \left|{^\sharp}\Hat{\P}_j^{\alpha~\beta}(\xi)\right|~\leq~\C_{\Re \alpha~\Re\beta} ~e^{\c|\Im \alpha|}e^{\c|\Im\beta|}~ \left|{1\over 1-|\xi|}\right|^{{1\over 2}-\ve} 
\end{array}
\eeq
for some $\ve=\ve(\Re\alpha,\Re\beta)>0$.
\end{lemma}
Consider
\bel{P j m sharp}
\begin{array}{cc}\ds
{^\sharp}\Hat{\P}_j^{\alpha~\beta}(\xi)~=~\sum_{m=1}^M {^\sharp}\Hat{\P}_{j~m}^{\alpha~\beta}(\xi),
\\\\ \ds
{^\sharp}\Hat{\P}_{j~m}^{\alpha~\beta}(\xi)~=~\Hat{\phi}(\xi)\int_{\lambda_{m-1}\leq|r|<\lambda_m} e^{-2\pi\i  r}{^\sharp}\Hat{\Omega}^\alpha (r\xi)\omega(r) |r|^{2\beta-1} dr,
\\\\ \ds
\hbox{$\lambda_m\in[2^{j-1},2^j]$}, 
\qquad
 \hbox{$\lambda_0=2^{j-1}$, $\lambda_M=2^j$ ~{\small and} ~$2^{\epsilon j-1}~\leq~\lambda_m-\lambda_{m-1}<2^{\epsilon j}$}
\end{array}
\eeq
where $0<\sigma=\sigma(\Re\alpha,\Re\beta)<{1\over 2}$ can be chosen sufficiently small.

Let $\varphi^\nu_j$ defined in (\ref{phi^v_j}).  Assert
\bel{P sharp v}
\begin{array}{cc}\ds
{^\sharp}\Hat{\P}_{j~m}^{\alpha~\beta}(\xi)~=~\sum_{\nu\colon\xi^\nu_j\in\mathds{S}^{n-1}} {^\sharp}\Hat{\P}_{j~m}^{\alpha~\beta~\nu}(\xi),
\\\\ \ds
  {^\sharp}\Hat{\P}_{j~m}^{\alpha~\beta~\nu}(\xi)~=~\varphi^\nu_j(\xi) \Hat{\phi}(\xi)\int_{\lambda_{m-1}\leq|r|<\lambda_m} e^{-2\pi\i  r}{^\sharp}\Hat{\Omega}^\alpha (r\xi)\omega(r) |r|^{2\beta-1} dr.
\end{array}
\eeq
For every $j>0$ and $m=1,\ldots,M$, we define
\bel{U sharp}
{^\sharp}\U^{\alpha~\beta}_{j~m}(x)~=~\sum_{\nu\colon\xi^\nu_j\in\mathds{S}^{n-1}} {^\sharp}\U^{\alpha~\beta~\nu}_{j~m}(x),\qquad
{^\sharp}\U^{\alpha~\beta~\nu}_{j~m}(x)~=~{^\sharp}\P^{\alpha~\beta~\nu}_{j~m}(x) \left[1-\Psi^{\alpha~\nu}_{j~m}(x)\right],
\eeq
\bel{V sharp}
{^\sharp}\V^{\alpha~\beta}_{j~m}(x)~=~\sum_{\nu\colon\xi^\nu_j\in\mathds{S}^{n-1}} {^\sharp}\V^{\alpha~\beta~\nu}_{j~m}(x),\qquad
{^\sharp}\V^{\alpha~\beta~\nu}_{j~m}(x)~=~{^\sharp}\P^{\alpha~\beta~\nu}_{j~m}(x) \Psi^{\alpha~\nu}_{j~m}(x).
 \eeq

{\bf Proposition One}~~{\it Let $ \Re\alpha\ge\left[{2n\over 2n-1}\right]\Re\beta$ and ${2n-1\over 4n}<\Re\beta<{2n-1\over 2n-2}$. We have
\bel{Result One U}
\left| {^\sharp}\Hat{\U}_{j~m}^{\alpha~\beta}(\xi)\right|~\leq~\C_{\Re \alpha~\Re\beta} ~e^{\c|\Im \alpha|}e^{\c|\Im\beta|}~ 2^{-(1-\epsilon)j}2^{j/2} 2^{-\ve j}
\eeq
and
\bel{Result One V}
\left| {^\sharp}\Hat{\V}_{j~m}^{\alpha~\beta}(\xi)\right|~\leq~\C_{\Re \alpha~\Re\beta} ~e^{\c|\Im \alpha|}e^{\c|\Im\beta|} ~2^{-(1-\epsilon)j}2^{-\ve j}
\eeq
for some $\ve=\ve(\Re\alpha,\Re\beta)>0$.}

\subsection{${^\flat}\U^{\alpha~\beta}_{j~m}$ and ${^\flat}\V^{\alpha~\beta}_{j~m}$}
${^\flat}\Omega^\alpha$ is a distribution defined in $\R^n$ by analytic continuation from
\bel{Omega flat}
\begin{array}{cc}\ds
 {^\flat}\Omega^{\alpha}(x)~=~\pi^{\alpha-1}\Gamma^{-1}\left(\alpha\right) 
 \left({1 \over1-|x|^2}\right)^{1-\alpha}_+,\qquad \hbox{\small{$\Re\alpha>0$}}. 
 \end{array}
\eeq
Equivalently, it can be defined by
 \bel{Omega flat Transform} 
\begin{array}{lr}\ds
{^\flat}\Hat{\Omega}^{\alpha}(\xi)~=~\left({1\over|\xi|}\right)^{{n\over 2}-(1-\alpha)} \J_{{n\over 2}-(1-\alpha)}\Big(2\pi|\xi|\Big)
\\\\ \ds~~~~~~~~~~~
~=~\left({1\over|\xi|}\right)^{{n-1\over 2}+\alpha-{1\over2}} \J_{{n-1\over 2}+\alpha-{1\over2}}\Big(2\pi|\xi|\Big),\qquad \alpha\in\Cx.
\end{array}
\eeq
Let $ \Re\alpha>0$ and $0<\Re\beta<{1\over 2}$. Define
\bel{P flat j}
\begin{array}{cc}\ds
{^\flat}\Hat{\P}_j^{\alpha~\beta}(\xi)~=~\Hat{\phi}(\xi)\int_{2^{j-1}\leq|r|<2^j} e^{-2\pi\i  r}{^\flat}\Hat{\Omega}^\alpha (r\xi)\omega(r) |r|^{2\beta-1} dr,\qquad j>0.
\end{array}
\eeq
Consider
\bel{P j m flat}
\begin{array}{cc}\ds
{^\flat}\Hat{\P}_j^{\alpha~\beta}(\xi)~=~\sum_{m=1}^M {^\flat}\Hat{\P}_{j~m}^{\alpha~\beta}(\xi),
\\\\ \ds
{^\flat}\Hat{\P}_{j~m}^{\alpha~\beta}(\xi)~=~\Hat{\phi}(\xi)\int_{\lambda_{m-1}\leq|r|<\lambda_m} e^{-2\pi\i  r}{^\flat}\Hat{\Omega}^\alpha (r\xi)\omega(r) |r|^{2\beta-1} dr,
\\\\ \ds
\hbox{$\lambda_m\in[2^{j-1},2^j]$}, 
\qquad
 \hbox{$\lambda_0=2^{j-1}$, $\lambda_M=2^j$ ~{\small and} ~$2^{\epsilon j-1}~\leq~\lambda_m-\lambda_{m-1}<2^{\epsilon j}$}
\end{array}
\eeq
where $0<\sigma=\sigma(\Re\alpha)<{1\over 2}$ can be chosen sufficiently small.

Let $\varphi^\nu_j$ defined in (\ref{phi^v_j}).  Assert
\bel{P flat v}
\begin{array}{cc}\ds
{^\flat}\Hat{\P}_{j~m}^{\alpha~\beta}(\xi)~=~\sum_{\nu\colon\xi^\nu_j\in\mathds{S}^{n-1}} {^\flat}\Hat{\P}_{j~m}^{\alpha~\beta~\nu}(\xi),
\\\\ \ds
 {^\flat}\Hat{\P}_{j~m}^{\alpha~\beta~\nu}(\xi)~=~\varphi^\nu_j(\xi)\Hat{\phi}(\xi)\int_{\lambda_{m-1}\leq|r|<\lambda_m} e^{-2\pi\i  r}{^\flat}\Hat{\Omega}^\alpha (r\xi)\omega(r) |r|^{2\beta-1} dr.
\end{array}
\eeq
For every $j>0$ and $m=1,\ldots,M$, we define
\bel{U flat}
{^\flat}\U^{\alpha~\beta}_{j~m}(x)~=~\sum_{\nu\colon\xi^\nu_j\in\mathds{S}^{n-1}} {^\flat}\U^{\alpha~\beta~\nu}_{j~m}(x),\qquad
{^\flat}\U^{\alpha~\beta~\nu}_{j~m}(x)~=~{^\flat}\P^{\alpha~\beta~\nu}_{j~m}(x) \left[1-\Psi^{\alpha~\nu}_{j~m}(x)\right]
\eeq
and
\bel{V flat}
{^\flat}\V^{\alpha~\beta}_{j~m}(x)~=~\sum_{\nu\colon\xi^\nu_j\in\mathds{S}^{n-1}} {^\flat}\V^{\alpha~\beta~\nu}_{j~m}(x),\qquad
{^\flat}\V^{\alpha~\beta~\nu}_{j~m}(x)~=~{^\flat}\P^{\alpha~\beta~\nu}_{j~m}(x) \Psi^{\alpha~\nu}_{j~m}(x).
 \eeq

{\bf Proposition Two}~~{\it Let $ \Re\alpha>0$ and $0<\Re\beta<{1\over 2}$.
We have
\bel{Result Two U}
\int_{\R^n}\left| {^\flat}\U_{j~m}^{\alpha~\beta}(x)\right|dx~\leq~\C_{N~\Re \alpha~\Re\beta} ~e^{\c|\Im \alpha|}~ 2^{-(1-\epsilon)j}2^{-Nj} 2^{-\ve j},\qquad \hbox{\small{$N\ge0$}}
\eeq
and
\bel{Result Two V}
\int_{\R^n}\left| {^\flat}\V_{j~m}^{\alpha~\beta}(x)\right|dx~\leq~\C_{\Re \alpha~\Re\beta} ~e^{\c|\Im \alpha|}~  2^{-(1-\epsilon)j}2^{-\ve j}
\eeq
for some $\ve=\ve(\Re\alpha)>0$.}

\section{Proof of Theorem Two}
\setcounter{equation}{0}
Recall $\I^\alpha f$ defined in (\ref{I alpha rewrite}) and $\Hat{\P}^\alpha_<$, $\Hat{\P}^\alpha_j$ and $\Hat{\P}^\alpha_{j~m}$ defined in (\ref{P sum j}) and (\ref{P j m}). 
We have
 \bel{I rewrite again}
 \begin{array}{ccc}\ds
 \I^\alpha f(x)
~=~\int_{\R^n} e^{2\pi\i x\cdot\xi} \Hat{f}(\xi)
\left\{\Hat{\P}^\alpha_<(\xi)+ \sum_{j>0} \Hat{\P}^\alpha_{j}(\xi)\right\} d\xi   
\\\\ \ds
\Hat{\P}^\alpha_<(\xi)
~=~\Hat{\phi}(\xi) \int_{|r|<1} e^{-2\pi\i  r}\Hat{\Omega}^\alpha (r\xi)\omega(r) |r|^{2\alpha-1} dr,
\\\\ \ds
 \Hat{\P}^\alpha_{j}(\xi)
~=~\Hat{\phi}(\xi) \int_{2^{j-1}\leq|r|<2^j} e^{-2\pi\i  r}\Hat{\Omega}^\alpha (r\xi)\omega(r) |r|^{2\alpha-1} dr
\end{array}
 \eeq
 where $\Hat{\phi}$ is defined in (\ref{hat phi}) and  
 \[\supp\Hat{\phi}~\subset~\Bigg\{\xi\in\R^n\colon{1\over 3}<|\xi|\leq3\Bigg\}.\]

 By using the norm estimate of Bessel functions in (\ref{J norm}):
 \[\left|{1\over \rho^{a+\i b}}\J_{a+\i b}\big(2\pi\rho\big)\right|~\leq~\C_a~e^{\c|b|} ~\left({1\over 1+\rho}\right)^{a+{1\over 2}},\qquad a,b\in\R,\qquad \rho>0\]
 we find
 \bel{Omega alpha norm EST}
 \begin{array}{lr}\ds
 \left| \Hat{\Omega}^\alpha (r\xi)\right|~=~\left| \left({1\over |r\xi|}\right)^{\alpha-{1\over2}} \J_{\alpha-{1\over 2}}\big(2\pi|r\xi|\big)\right|\qquad\hbox{\small{by (\ref{Omega^alpha Transform})}}
 \\\\ \ds~~~~~~~~~~~~~~
 ~\leq~\C_{\Re\alpha}~e^{\c|\Im\alpha|}~\left({1\over 1+|r|}\right)^{\Re\alpha},\qquad \hbox{\small{${1\over 3}<|\xi|\leq3$}}.
 \end{array}
 \eeq
 Recall $\omega(r)$ from (\ref{I alpha rewrite}). We have
\[
\omega(r)~=~e^{2\pi\i r}\int_0^1 e^{-2\pi\i \tau r} \tau d\tau~=~{-1\over 2\pi\i}{1\over r}-{1\over 4\pi^2r^2} \left[1-e^{-2\pi\i r}\right]
\]
implying
\bel{omega norm}
|\omega(r)|~\leq~\C\Big[1+|r|\Big]^{-1}.
\eeq
 Let ${1\over2} <\Re\alpha<1$. From (\ref{Omega alpha norm EST})-(\ref{omega norm}), we find
 \bel{Kernel EST}
 \begin{array}{lr}\ds
 \left|\sum_{j>0}  \Hat{\P}^\alpha_{j}(\xi)\right|~=~\left| \Hat{\phi}(\xi) \int_{|r|\ge1} e^{-2\pi\i  r}\Hat{\Omega}^\alpha (r\xi)\omega(r) |r|^{2\alpha-1} dr\right| 
 \\\\ \ds~~~~~~~~~~~~~~~~~
  ~\leq~\C_{\Re\alpha}~e^{\c|\Im\alpha|} \int_{|r|\ge1} \left({1\over 1+|r|}\right)^{\Re\alpha}  \Big[1+|r|\Big]^{-1} |r|^{2\Re\alpha-1}dr
   \\\\ \ds~~~~~~~~~~~~~~~~~
 ~\leq~\C_{\Re\alpha}~e^{\c|\Im\alpha|}. 
 \end{array}
 \eeq
Let $f\in\L^p(\R^n)$  and $p= {2n\over n+1}$. We have
\bel{f Transform}
\Hat{f}~\in~\L^{2n\over n-1}(\R^n).
\eeq
 Therefore, $\Hat{f}(\xi)$ is locally integrable for ${1\over 3}<|\xi|<3$. 
 
 By  applying the dominated convergence theorem, we write
 \bel{I SUM}
 \begin{array}{ccc}\ds
 \I^\alpha f(x)~=~\I^\alpha_< f(x)+\sum_{j>0} \I^\alpha_j f(x),
 \\\\ \ds
 \I^\alpha_< f(x)~=~\int_{\R^n} e^{2\pi\i x\cdot\xi} \Hat{f}(\xi)  \Hat{\P}^\alpha_<(\xi)d\xi,
 \\\\ \ds
   \I^\alpha_j f(x)~=~\int_{\R^n} e^{2\pi\i x\cdot\xi} \Hat{f}(\xi)  \Hat{\P}^\alpha_j(\xi)d\xi,\qquad j>0.
   \end{array}
   \eeq 
 The kernel of $\I^\alpha_<$ as shown in (\ref{kernel <}) is an $\L^1$-function in $\R^n$. 
From now on, we focus on $ \I^\alpha_{j} f, j>0$.

Let ${1\over 2}<\Re\alpha<1$. We can find $\a_i=\a_i(\Re\alpha), \b_i=\b_i(\Re\alpha), i=1,2$ such that
\bel{ab Constraints}
\begin{array}{cc}\ds
\a_1>0,\qquad  0<\b_1<{1\over 2},\qquad \a_2\ge\left[{2n\over 2n-1}\right]\b_2,\qquad {2n-1\over 4n}<\b_2<{2n-1\over 2n-2},
\\\\ \ds
\Re\alpha~=~{\a_1\over n}+\a_2\left({n-1\over n}\right)
~=~{\b_1\over n}+\b_2\left({n-1\over n}\right).
\end{array}
\eeq
 As $\Re\alpha\mt1$, we necessarily have 
$\b_1\mt{1\over 2},~ \b_2\mt {2n-1\over 2n-2}$ and $\a_1\mt0,~ \a_2\mt {n\over n-1}$.

For $0\leq\Re z\leq1$, we define
\bel{Omega^a_z}
\begin{array}{lr}\ds
\Hat{\Omega}^{\alpha~z}(\xi)~=~
 \left({1\over|\xi|}\right)^{\a_1 z+\a_2(1-z)-{1\over 2}+\left({n-1\over 2}\right)z-{1\over 2}(1-z)+\i\Im\alpha}
 \J_{\a_1 z+\a_2(1-z)-{1\over 2}+\left({n-1\over 2}\right)z-{1\over 2}(1-z)+\i\Im\alpha}\left(2\pi|\xi|\right)
 \\\\ \ds
 ~=~\pi^{\a_1 z+\a_2(1-z)-1+\left({n-1\over 2}\right)z-{1\over 2}(1-z)+\i\Im\alpha}\hbox{$\Gamma^{-1}\left(\a_1 z+\a_2(1-z)+\left({n-1\over 2}\right)z-{1\over 2}(1-z)+\i\Im\alpha\right)$}
 \\ \ds~~~~
 \int_{-1}^1 e^{2\pi\i |\xi|s} (1-s^2)^{\a_1 z+\a_2(1-z)-1+\left({n-1\over 2}\right)z-{1\over 2}(1-z)+\i\Im\alpha}ds\qquad\hbox{\small{by (\ref{Bessel})}}
 \end{array}
 \eeq
and
\bel{P^a_z}
\Hat{\P}^{\alpha~z}_j(\xi)~=~e^{\big[z-{1\over n}\big]^2}\Hat{\phi}(\xi)\int_{2^{j-1}\leq|r|<2^j} e^{-2\pi\i r} \Hat{\Omega}^{\alpha~z}(r\xi) \omega(r)|r|^{2\big[\b_1z+\b_2(1-z)\big]-1+2\i\Im\alpha} dr,\qquad\hbox{\small{$j>0$}}.
\eeq
Let $\lambda_m\in[2^{j-1},2^j]$: $\lambda_0=2^{j-1}$, $\lambda_M=2^j$ and
$2^{\epsilon j-1}\leq\lambda_m-\lambda_{m-1}<2^{\epsilon j}$ for  which $0<\sigma=\sigma(\Re\alpha)<{1\over 2}$ can be chosen sufficiently small. We have
\[\begin{array}{cc}\ds
\I^\alpha_j f(x)~=~\sum_{m=1}^M \I^\alpha_{j~m} f(x),
\\\\ \ds

\I^\alpha_{j~m} f(x)~=~\int_{\R^n} e^{2\pi\i x\cdot\xi} \Hat{f}(\xi)  \Hat{\P}^\alpha_{j~m}(\xi)d\xi,
\\\\ \ds
\Hat{\P}^\alpha_{j~m}(\xi)~=~\Hat{\phi}(\xi) \int_{\lambda_{m-1}\leq|r|<\lambda_m} e^{-2\pi\i  r}\Hat{\Omega}^\alpha (r\xi)\omega(r) |r|^{2\alpha-1} dr.
\end{array}\]

Consider
\bel{P^a_z m}
\begin{array}{cc}\ds
\Hat{\P}^{\alpha~z}_j(\xi)~=~\sum_{m=1}^M \Hat{\P}^{\alpha~z}_{j~m}(\xi),
\\\\ \ds
\Hat{\P}^{\alpha~z}_{j~m}(\xi)~=~e^{\big[z-{1\over n}\big]^2}\Hat{\phi}(\xi)\int_{\lambda_{m-1}\leq|r|<\lambda_m} e^{-2\pi\i r} \Hat{\Omega}^{\alpha~\a}_z (r\xi) \omega(r)|r|^{2\big[\b_1z+\b_2(1-z)\big]-1+2\i\Im\alpha} dr.
\end{array}
\eeq
Observe that $\Hat{\P}^{\alpha~{1\over n}}_{j~m}(\xi)=\Hat{\P}^{\alpha}_{j~m}(\xi)$.
 
Define
 \bel{I^a_z,j}
 \I^{\alpha~z}_{j~m} f(x)~=~\int_{\R^n} e^{2\pi\i x\cdot\xi} \Hat{f}(\xi) \Hat{\P}^{\alpha~z}_{j~m}(\xi)d\xi
\eeq
for every $ j>0$, $m=1,2,\ldots,M$.

Given $j>0$, $\left\{ \xi^\nu_j\right\}_\nu$ is a collection of points almost equally distributed on  $\mathds{S}^{n-1}$ having  a grid length between $2^{-j/2-1}$ and $2^{-j/2}$. See  {\bf (1)}-{\bf(2)} above (\ref{phi^v_j intro})-(\ref{Gamma_j}).

Let $\varphi^\nu_j$ defined in (\ref{phi^v_j}).  Assert
\bel{P z Transform j}
\begin{array}{cc}\ds
\Hat{\P}^{\alpha~z}_{j~m}(\xi)~=~\sum_{\nu\colon\xi^\nu_j\in\mathds{S}^{n-1}} ~\Hat{\P}^{\alpha~\nu~z}_{j~m}(\xi),
\\\\ \ds
\Hat{\P}^{\alpha~\nu~z}_{j~m}(\xi)~=~e^{\big[z-{1\over n}\big]^2}\Hat{\phi}(\xi)\varphi^\nu_j(\xi)\int_{\lambda_{m-1}\leq|r|<\lambda_m} e^{-2\pi\i r} \Hat{\Omega}^{\alpha~z} (r\xi) \omega(r)|r|^{2\big[\b_1z+\b_2(1-z)\big]-1+2\Im\alpha} dr.
\end{array}
\eeq
Let $\Psi^{\alpha}_{j~m}, \Psi^{\alpha~\mu}_{j~m}$ defined in (\ref{Psi mu j}).
For every $j>0$ and $m=1,\ldots,M$, we define 
\bel{II_j z}
 \begin{array}{cc}\ds
  \II^{\alpha~z}_{j~m} f(x)~=~\int_{\R^n} f(x-u) \U^{\alpha~z}_{j~m}(u)du,
  \\\\ \ds
   \U^{\alpha~z}_{j~m}(x)~=~\sum_{\nu\colon\xi^\nu_j\in\mathds{S}^{n-1}}  \U^{\alpha~\nu~z}_{j~m}(x),\qquad  
  \U^{\alpha~\nu~z}_{j~m}(x)~=~\P^{\alpha~\nu~z}_{j~m}(x) \left[1-\Psi^{\alpha~\nu}_{j~m}(x)\right] 
    \end{array}
  \eeq 
 and
 \bel{III_j z}
 \begin{array}{cc}\ds
  \III^{\alpha~z}_{j~m} f(x)~=~\int_{\R^n} f(x-u) \V^{\alpha~z}_{j~m}(u)du,
  \\\\ \ds
   \V^{\alpha~z}_{j~m}(x)~=~\sum_{\nu\colon\xi^\nu_j\in\mathds{S}^{n-1}}  \V^{\alpha~\nu~z}_{j~m}(x),\qquad  
  \V^{\alpha~\nu~z}_{j~m}(x)~=~\P^{\alpha~\nu~z}_{j~m}(x) \Psi^{\alpha~\nu}_{j~m}(x). 
    \end{array}
 \eeq 
In particular, at $z={1\over n}$, we find  
\[\Hat{\P}^{\alpha~\nu~{1\over n}}_{j~m}(\xi)~=~\Hat{\P}^{\alpha~\nu}_{j~m}(\xi),\qquad \II^{\alpha~{1\over n}}_{j~m} f(x)~=~\II^{\alpha}_j f(x),\qquad \III^{\alpha~{1\over n}}_{j~m} f(x)~=~\III^{\alpha}_j f(x)\]
 as defined in (\ref{II_j U again})-(\ref{III_j V again}).

Let $\Re z=0$. From (\ref{II_j z})-(\ref{III_j z}), we have
$\II^{\alpha~0+\i\Im z}_{j~m} f(x)=e^{\big[\i\Im z-{1\over n}\big]^2}f\ast{^\sharp}\U^{\A_2~\B_2}_{j~m}(x)$ and $\III^{\alpha~0+\i\Im z}_{j~m} f(x)=e^{\big[\i\Im z-{1\over n}\big]^2} f\ast {^\sharp}\V^{\A_2~\B_2}_j f(x)$
of which
 ${^\sharp}\U^{\A_2~\B_2}_{j~m}$, ${^\sharp}\V^{\A_2~\B_2}_{j~m}$  are defined in 
(\ref{U sharp})-(\ref{V sharp})  with $\A_2=\a_2+\i\left[ \a_1-\a_2+\left({n-1\over 2}\right)+{1\over 2}\right]\Im z+\i\Im\alpha$ and $\B_2=\b_2+\i\left[ \b_1-\b_2\right]\Im z+\i\Im \alpha$.  
 
Recall $\a_2\ge\left[{2n\over 2n-1}\right]\b_2$, ${2n-1\over 4n}<\b_2<{2n-1\over 2n-2}$ in (\ref{ab Constraints}).
By applying  (\ref{Result One U})  in {\bf Proposition One} together with  Plancherel theorem, we find
\bel{0 EST II}
\begin{array}{lr}\ds
\left\| \II^{\alpha~0+\i\Im z}_{j~m}f\right\|_{\L^2(\R^n)} ~=~\left\| \Hat{\II^{\alpha~0+\i\Im z}_{j~m}f}\right\|_{\L^2(\R^n)}
\\\\ \ds~~~~~~~~~~~~~~~~~~~~~~~~~~~~~~~
~\leq~\C_{\a_2~\b_2} ~e^{\c|\Im\alpha|}~2^{-(1-\epsilon)j}2^{j/2} 2^{-\ve j}~\left\| \Hat{f}\right\|_{\L^2(\R^n)} 
\\\\ \ds~~~~~~~~~~~~~~~~~~~~~~~~~~~~~~~
~\leq~\C_{\Re\alpha} ~e^{\c|\Im\alpha|}~2^{-(1-\epsilon)j}2^{j/2} 2^{-\ve j} ~\left\| f\right\|_{\L^2(\R^n)} 
\end{array}
\eeq
for some $\ve=\ve(\a_2,\b_2)=\ve(\Re\alpha)>0$.

On the other hand, by using (\ref{Result One V}) in {\bf Proposition One} together with Plancherel theorem, we have
\bel{0 EST III}
\left\| \III^{\alpha~0+\i\Im z}_{j~m}f\right\|_{\L^2(\R^n)} ~\leq~\C_{\Re\alpha} ~e^{\c|\Im\alpha|}~2^{-(1-\epsilon)j}2^{-\ve j} ~\left\| f\right\|_{\L^2(\R^n)}. 
\eeq
Let $\Re z=1$. From (\ref{II_j z})-(\ref{III_j z}), we find
$\II^{\alpha~1+\i\Im z}_{j~m} f(x)=e^{\big[\i\Im z+{n-1\over n}\big]^2}f\ast{^\flat}\U^{\A_1~\B_1}_{j~m}(x)$ and $\III^{\alpha~1+\i\Im z}_{j~m} f(x)=e^{\big[\i\Im z+{n-1\over n}\big]^2} f\ast {^\flat}\V^{\A_1~\B_1}_{j~m} f(x)$
where
 ${^\flat}\U^{\A_1~\B_1}_{j~m}$, ${^\flat}\V^{\A_1~\B_1}_{j~m}$ are defined in 
(\ref{U flat})-(\ref{V flat})  with $\A_1=\a_1+\i\left[ \a_1-\a_2+\left({n-1\over 2}\right)+{1\over 2}\right]\Im z+\i\Im\alpha$ and $\B_1=\b_1+\i\left[ \b_1-\b_2\right]\Im z+\i\Im \alpha$. 

Note that  $\a_1>0, 0<\b_1<{1\over 2}$ in (\ref{ab Constraints}). By applying  (\ref{Result Two U})  in {\bf Proposition Two}  and using Minkowski integral  inequality, we have
\bel{1 EST II}
\begin{array}{lr}\ds
\left\|\II^{\alpha~1+\i\Im z}_{j~m}f\right\|_{\L^1(\R^n)}
~\leq~\C_{N~\a_1~\b_1}~ e^{\c|\Im\alpha|}~2^{-(1-\epsilon)j}2^{-jN} 2^{-j\ve}~\left\| f\right\|_{\L^1(\R^n)} 
\\\\ \ds~~~~~~~~~~~~~~~~~~~~~~~~~~~~~~~~~~~~~~~
~=~\C_{N~\Re\alpha}~2^{-(1-\epsilon)j}2^{-Nj} 2^{-\ve j} \left\| f\right\|_{\L^1(\R^n)},\qquad\hbox{\small{$N\ge0$}}
\end{array}
\eeq
for some $\ve=\ve(\a_1,\b_1)=\ve(\Re\alpha)>0$. 

On the other hand, by applying (\ref{Result Two V}) in {\bf Proposition Two} and using Minkowski integral  inequality, we find
\bel{1 EST III}
\left\|\III^{\alpha~1+\i\Im z}_{j~m}f\right\|_{\L^1(\R^n)}
~\leq~\C_{\Re\alpha}~2^{-(1-\epsilon)j} 2^{-\ve j} \left\| f\right\|_{\L^1(\R^n)}.
\eeq
Recall {\bf Remark 1.2}. There are a constant multiple of $2^{(1-\epsilon)j}$ many $1\leq m\leq M$. Take into account ${1\over p}={n+1\over 2n}={1\over n}+{1\over 2}{n-1\over n}$.

From (\ref{0 EST II}) and (\ref{1 EST II}) with  $N\ge{n-1\over 2}$, by applying Stein interpolation theorem \cite{Stein'}, we obtain
\bel{EST II}
\begin{array}{lr}\ds
\left\| \II^\alpha_j f\right\|_{\L^p(\R^n)}~\leq~\sum_{m=1}^M \left\| \II^{\alpha}_{j~m} f\right\|_{\L^p(\R^n)}
~\leq~\C_{\Re\alpha} ~e^{\c|\Im\alpha|}~ 2^{-\ve j}~\left\| f\right\|_{\L^p(\R^n)}.
\end{array}
\eeq
From (\ref{0 EST III}) and (\ref{1 EST III}), by applying Stein interpolation theorem \cite{Stein'} again,  we have
\bel{EST III}
\left\| \III^\alpha_j f\right\|_{\L^p(\R^n)}~\leq~\sum_{m=1}^M \left\| \II^{\alpha}_{j~m} f\right\|_{\L^p(\R^n)}
~\leq~\C_{\Re\alpha} ~e^{\c|\Im\alpha|}~ 2^{-\ve j}~\left\| f\right\|_{\L^p(\R^n)}.
\eeq
By summing over every $j>0$ and using (\ref{EST II})-(\ref{EST III}), we conclude 
\bel{I End-point p}
\left\| \I^\alpha f\right\|_{\L^p(\R^n)}~\leq~\C_{\Re\alpha} ~e^{\c|\Im\alpha|}~\left\| f\right\|_{\L^p(\R^n)}.
\eeq
Lastly, the adjoint operator of $\I^\alpha$ defined in (\ref{I alpha}) can be simply given by $\I^{\bar{\alpha}}$ where $\bar{\alpha}$ is the complex conjugate of $\alpha$. Clearly, it satisfies all above estimates. By duality, we obtain 
\bel{I End-point p'}
\left\| \I^\alpha f\right\|_{\L^{p\over p-1}(\R^n)}~\leq~\C_{\Re\alpha} ~e^{\c|\Im\alpha|}~\left\| f\right\|_{\L^{p\over p-1}(\R^n)}.
\eeq
The $\L^p$-norm inequality in (\ref{Result Two}) follows by applying Riesz-Thorin interpolation theorem.

\section{Proof of Proposition One}
\setcounter{equation}{0}
Recall  (\ref{P j m sharp}). Let $\Re\alpha\ge\left[{2n\over 2n-1}\right]\Re\beta$ and  ${2n-1\over 4n}<\Re\beta<{2n-1\over 2n-2}$. We have
\[
\begin{array}{cc}\ds
{^\sharp}\Hat{\P}_j^{\alpha~\beta}(\xi)~=~\sum_{m=1}^M \Hat{\P}_{j~m}^{\alpha~\beta}(\xi),
\\\\ \ds
  {^\sharp}\Hat{\P}_{j~m}^{\alpha~\beta}(\xi)~=~ \Hat{\phi}(\xi)\int_{\lambda_{m-1}\leq|r|<\lambda_m} e^{-2\pi\i  r}{^\sharp}\Hat{\Omega}^\alpha (r\xi)\omega(r) |r|^{2\beta-1} dr
\end{array}
\]
where ${^\sharp}\Hat{\Omega}^\alpha $ is defined in (\ref{Omega sharp Transform}).

Moreover, $\Hat{\phi}$ is defined in (\ref{hat phi}).  $\omega$ is given in (\ref{I alpha rewrite}). 
$\lambda_m\in[2^{j-1},2^j]$: $\lambda_0=2^{j-1}$, $\lambda_M=2^j$ and $2^{\epsilon j-1}\leq\lambda_m-\lambda_{m-1}<2^{\epsilon j}$  for which
$0<\sigma=\sigma(\Re\alpha,\Re\beta)<{1\over 2}$ can be chosen sufficiently small. 
\v
{\bf Lemma One}~~{\it Given $j>0$ and $m=1,2,\ldots,M$, we have
\bel{Result 1 > j}
\left| {^\sharp}\Hat{\P}_{j~m}^{\alpha~\beta}(\xi)\right|~\leq~\C_{\Re \alpha~\Re\beta} ~e^{\c|\Im \alpha|}~ 2^{-(1-\epsilon)j}2^{j/2} 2^{-\ve j},\qquad \left|1-|\xi|\right|\leq 2^{-j}
\eeq
and
\bel{Result 1 > xi}
\left| {^\sharp}\Hat{\P}_{j~m}^{\alpha~\beta}(\xi)\right|~\leq~\C_{\Re \alpha~\Re\beta} ~e^{\c|\Im \alpha|}e^{\c|\Im\beta|}~ 2^{-(1-\epsilon)j}\left|{1\over 1-|\xi|}\right|^{1\over 2} 2^{-\ve j},\qquad \left|1-|\xi|\right|> 2^{-j}
\eeq
for some $\ve=\ve(\Re\alpha,\Re\beta)>0$.}
\begin{remark} As required for later estimates, we choose $\sigma=\sigma(\Re\alpha,\Re\beta)$ to satisfy $n\sigma\leq{1\over 2}\ve$.
\end{remark} 
Recall {\bf Remark 1.2}. There are 
 at most a constant multiple of $ 2^{(1-\epsilon)j}$ many $\lambda_m\in[2^{j-1},2^j]$.
Consequently, by applying  (\ref{Result 1 > j}) and (\ref{Result 1 > xi}), we find
\bel{P_j sharp EST<}
\begin{array}{lr}\ds
\left| {^\sharp}\Hat{\P}_j^{\alpha~\beta}(\xi)\right|~\leq~\sum_{m=1}^M \left| {^\sharp}\Hat{\P}_{j~m}^{\alpha~\beta}(\xi)\right|
\\\\ \ds~~~~~~~~~~~~~~~~
~\leq~\C_{\Re \alpha~\Re\beta} ~e^{\c|\Im \alpha|}~ 2^{j/2} 2^{-\ve j},\qquad \left|1-|\xi|\right|\leq 2^{-j}
\end{array}
\eeq
and
\bel{P_j sharp EST>}
\begin{array}{lr}\ds
\left| {^\sharp}\Hat{\P}_j^{\alpha~\beta}(\xi)\right|~\leq~\sum_{m=1}^M \left| {^\sharp}\Hat{\P}_{j~m}^{\alpha~\beta}(\xi)\right|
\\\\ \ds~~~~~~~~~~~~~~~
~\leq~\C_{\Re \alpha~\Re\beta} ~e^{\c|\Im \alpha|}e^{\c|\Im\beta|}~ \left|{1\over 1-|\xi|}\right|^{1\over 2} 2^{-\ve j},\qquad \left|1-|\xi|\right|> 2^{-j}.
\end{array}
\eeq
\vsk
Recall {\bf Lemma 3.1}. Consider
\[
\begin{array}{lr}\ds
\sum_{j>0} \left|{^\sharp}\Hat{\P}_j^{\alpha~\beta}(\xi)\right|~\leq~\sum_{j>0\colon |1-|\xi||\leq2^{-j}} \left|{^\sharp}\Hat{\P}_j^{\alpha~\beta}(\xi)\right|+\sum_{j>0\colon |1-|\xi||>2^{-j}} \left|{^\sharp}\Hat{\P}_j^{\alpha~\beta}(\xi)\right|.
\end{array}
\]
We have
\bel{T est.1}
\begin{array}{lr}\ds
\sum_{j>0\colon |1-|\xi||\leq2^{-j}} \left|{^\sharp}\Hat{\P}_j^{\alpha~\beta}(\xi)\right|
~\leq~\C_{\Re\alpha~\Re\beta}~e^{\c\Im\alpha}  \sum_{j>0\colon |1-|\xi||\leq2^{-j}}  2^{j/2} 2^{-\ve j}\qquad\hbox{\small{by (\ref{P_j sharp EST<})}}
\\\\ \ds~~~~~~~~~~~~~~~~~~~~~~~~~~~~~~~~~~~
~\leq~\C_{\Re\alpha~\Re\beta}~e^{\c\Im\alpha} \left|{1\over 1-|\xi|}\right|^{{1\over 2}-{\ve\over 2}} ~\sum_{j>0} 2^{-\big({\ve\over2}\big)j}
\\\\ \ds~~~~~~~~~~~~~~~~~~~~~~~~~~~~~~~~~~~
~\leq~\C_{\Re\alpha~\Re\beta}~e^{\c\Im\alpha} \left|{1\over 1-|\xi|}\right|^{{1\over 2}-{\ve\over 2}}
\end{array}
\eeq
and
\bel{T est.2}
\begin{array}{lr}\ds
\sum_{j>0\colon |1-|\xi||>2^{-j}} \left|{^\sharp}\Hat{\P}_j^{\alpha~\beta}(\xi)\right|
~\leq~\C_{\Re\alpha~\Re\beta}~e^{\c\Im\alpha}e^{\c\Im\beta} 
 \sum_{j>0\colon |1-|\xi||>2^{-j}}\left|{1\over 1-|\xi|}\right|^{1\over 2} 2^{-\ve j}
\qquad
 \hbox{\small{by (\ref{P_j sharp EST>})}}
\\\\ \ds~~~~~~~~~~~~~~~~~~~~~~~~~~~~~~~~~~~
~\leq~\C_{\Re\alpha~\Re\beta}~e^{\c\Im\alpha}e^{\c\Im\beta} \left|{1\over 1-|\xi|}\right|^{{1\over 2}-{\ve\over 2}} ~\sum_{j>0} 2^{-\big({\ve\over2}\big)j}
\\\\ \ds~~~~~~~~~~~~~~~~~~~~~~~~~~~~~~~~~~~
~\leq~\C_{\Re\alpha~\Re\beta}~e^{\c\Im\alpha}e^{\c\Im\beta} \left|{1\over 1-|\xi|}\right|^{{1\over 2}-{\ve\over 2}}.
\end{array}
\eeq
From (\ref{T est.1})-(\ref{T est.2}), we conclude (\ref{Kernel T Est}).

\subsection{Proof of Lemma One}
Recall ${^\sharp}\Omega^\alpha$ defined in (\ref{Omega sharp})-(\ref{Omega sharp Transform}).
We have
\[
\begin{array}{lr}\ds
{^\sharp}\Hat{\P}^{\alpha~\beta}_{j~m}(\xi)~=~\Hat{\phi}(\xi)\int_{\lambda_{m-1}\leq|r|<\lambda_m} e^{-2\pi\i  r}{^\sharp}\Omega^\alpha(r\xi) \omega(r) |r|^{2\beta-1} dr
\\\\ \ds~~~~~~~~~~~~~
~=~\Hat{\phi}(\xi)\int_{\lambda_{m-1}\leq|r|<\lambda_m} e^{-2\pi\i  r} \left({1\over|r\xi|}\right)^{\alpha-1} \J_{\alpha-1}\left(2\pi|r\xi|\right)
 \omega(r) |r|^{2\beta-1}  dr.
 \end{array}
\]
Note that 
$\supp\Hat{\phi}\subset\Big\{\xi\in\R^n\colon{1\over 3}<|\xi|\leq3\Big\}$ 
and
 \[\omega(r)~=~{-1\over 2\pi\i}{1\over r}-{1\over 4\pi^2r^2} \left[1-e^{-2\pi\i r}\right]\]
as shown in (\ref{I alpha rewrite}).

 For symmetry reason, we consider $r>0$ only. 
 Define
    \bel{Q_j int sharp}
{^\sharp}\Hat{\Q}_{j~m}^{\alpha~\beta}(\xi)
 ~=~\Hat{\phi}(\xi)\int_{\lambda_{m-1}}^{\lambda_m} e^{-2\pi\i  r} \left({1\over r|\xi|}\right)^{\alpha-1} \J_{\alpha-1}\left(2\pi r|\xi|\right)
 r^{2\beta-2}  dr
\eeq
and
    \bel{R_j int sharp}
    \begin{array}{lr}\ds
{^\sharp}\Hat{\hbox{\bf R}}_{j~m}^{\alpha~\beta}(\xi)
 ~=~\Hat{\phi}(\xi)\int_{\lambda_{m-1}}^{\lambda_m} e^{-2\pi\i  r} \left[1-e^{-2\pi\i r}\right] \left({1\over r|\xi|}\right)^{\alpha-1} \J_{\alpha-1}\left(2\pi r|\xi|\right)
 r^{2\beta-3}  dr.
\end{array}
\eeq
Because $ \Re\alpha\ge\left[{2n\over 2n-1}\right]\Re\beta$ and ${2n-1\over 4n}<\Re\beta<{2n-1\over 2n-2}$, we find 
\bel{Constraints}
\Re\alpha~\ge~\left[{2n\over 2n-1}\right]\Re\beta~>~{1\over 2},\qquad  2\Re\beta-\Re\alpha~\leq~\left[{2n-2\over 2n-1}\right]\Re\beta~<~1.
\eeq 
By using  the norm estimate of Bessel functions in (\ref{J norm}), we have
\[
 \begin{array}{lr}\ds
\left| {^\sharp}\Omega^\alpha(r\xi)\right|    ~=~\left({1\over r|\xi|}\right)^{\Re\alpha-1} \left|\J_{\alpha-1}\left(2\pi r|\xi|\right)\right|
\\\\ \ds~~~~~~~~~~~~~~~
~\leq~\C_{\Re \alpha}~e^{\c|\Im \alpha|} ~ \Bigg\{{1\over 1+|r\xi|}\Bigg\}^{\Re \alpha-{1\over 2}}.
 \end{array}
\]
This further implies
\bel{R_j EST sharp}
\begin{array}{lr}\ds
\left| {^\sharp}\Hat{\hbox{\bf R}}_j^{\alpha~\beta}(\xi)\right|~\leq~
\C_{\Re \alpha} ~e^{\c|\Im \alpha|}~
\int_{\lambda_{m-1}}^{\lambda_m} 
\left\{{1\over1+r}\right\}^{\Re\alpha-{1\over 2}}  r^{2\Re\beta-3}  dr\qquad\hbox{\small{$\left(~{1\over 3}<|\xi|\leq3~\right)$}}
\\\\ \ds~~~~~~~~~~~~~~~~
~\leq~\C_{\Re \alpha} ~e^{\c|\Im \alpha|}~
\int_{\lambda_{m-1}}^{\lambda_m}   r^{2\Re\beta-\Re\alpha-2-{1\over 2}}  dr 
\\\\ \ds~~~~~~~~~~~~~~~~
~\leq~\C_{\Re \alpha~\Re\beta} ~e^{\c|\Im \alpha|}~2^{\left[2\Re\beta-\Re\alpha-2-{1\over 2}\right]j} \left|\lambda_m-\lambda_{m-1}\right|
\\\\ \ds~~~~~~~~~~~~~~~~
~\leq~\C_{\Re \alpha~\Re\beta} ~e^{\c|\Im \alpha|}~2^{\left[2\Re\beta-\Re\alpha-2-{1\over 2}\right]j} 2^{\epsilon j}
\\\\ \ds~~~~~~~~~~~~~~~~
~=~\C_{\Re \alpha~\Re\beta} ~e^{\c|\Im \alpha|}~2^{-(1-\epsilon)j}2^{\left[2\Re\beta-\Re\alpha-1\right]j}2^{-j/2},\qquad\hbox{\small{$2\Re\beta-\Re\alpha-1<0$}}.
\end{array}
\eeq
Next, we aim to show
\bel{Q_j EST sharp j <}
\left| {^\sharp}\Hat{\Q}_{j~m}^{\alpha~\beta}(\xi)\right|~\leq~\C_{\Re \alpha~\Re\beta} ~e^{\hbox{\small{{\bf c}}}|\Im \alpha|}~2^{-(1-\epsilon)j}2^{j/2} 2^{-j\ve},\qquad |1-|\xi||\leq 2^{-j}
\eeq
and
\bel{Q_j EST sharp j >}
\left| {^\sharp}\Hat{\Q}_{j~m}^{\alpha~\beta}(\xi)\right|~\leq~\C_{\Re \alpha~\Re\beta} ~e^{\hbox{\small{{\bf c}}}|\Im \alpha|}e^{\hbox{\small{{\bf c}}}|\Im \beta|}~2^{-(1-\epsilon)j}\left|{1\over 1-|\xi|}\right|^{1\over 2} 2^{-j\ve},\qquad |1-|\xi||> 2^{-j}
\eeq
for some $\ve=\ve(\Re\alpha,\Re\beta)>0$.

Later, it should be clear that the same estimates regarding ${^\sharp}\Hat{\Q}_{j~m}^{\alpha~\beta}$ and ${^\sharp}\Hat{\hbox{\bf R}}_{j~m}^{\alpha~\beta}$ remain to be true for $r\in[-\lambda_m, -\lambda_{m-1}]$.
Therefore, we conclude (\ref{Result 1 > j}) and (\ref{Result 1 > xi})
from (\ref{R_j EST sharp}) and (\ref{Q_j EST sharp j <})-(\ref{Q_j EST sharp j >}).

By using the asymptotic expansion of Bessel functions in (\ref{J asymptotic}), 
we write 
\bel{Q=S+E}
\begin{array}{lr}\ds
{^\sharp}\Hat{\Q}_{j~m}^{\alpha~\beta}(\xi)~=~\Hat{\phi}(\xi)\int_{\lambda_{m-1}}^{\lambda_m} 
e^{-2\pi\i r} 
\Big[{^\sharp}\mathfrak{S}^\alpha\left(r|\xi|\right)+{^\sharp}\mathcal{E}^\alpha\left(r|\xi|\right)\Big]
r^{2\beta-2}  dr
\end{array}
\eeq
where
\bel{S principal}
{^\sharp}\mathfrak{S}^\alpha\left(\rho\right)~=~{1\over \pi}  \left({1\over \rho}\right)^{\alpha-{1\over 2}} \cos\left[2\pi \rho-{\pi\over 2}\alpha+{\pi\over 4}\right],\qquad\rho>0
\eeq
and
\bel{E}
\left| {^\sharp}\mathcal{E}^\alpha\left(\rho\right)\right|~\leq~\C_{\Re\alpha}~e^{\c|\Im\alpha|} \left\{ \begin{array}{lr}\ds \rho^{-\Re\alpha+{1\over 2}},\qquad 0<\rho\leq1,
\\ \ds
\rho^{-\Re\alpha-{1\over 2}},\qquad~~~ \rho>1.
\end{array}
\right.
\eeq
By using  (\ref{E}), we find   
 \bel{E est}
\begin{array}{lr}\ds
\left| \Hat{\phi}(\xi)\int_{\lambda_{m-1}}^{\lambda_m} 
e^{-2\pi\i r} 
{^\sharp}\mathcal{E}^\alpha\left(r|\xi|\right)
r^{2\beta-2}  dr  \right|
 ~\leq~\C_{\Re\alpha}~e^{\c|\Im\alpha|} 
  \int_{\lambda_{m-1}}^{\lambda_m}  \left({1\over r}\right)^{\Re\alpha+{1\over 2}} r^{2\Re\beta-2} dr 
\\\\ \ds~~~~~~~~~~~~~~~~~~~~~~~~~~~~~~~~~~~~~~~~~~~~~~~~~~~~~~~~
 ~\leq~\C_{\Re\alpha~\Re\beta}~e^{\c|\Im\alpha|}~2^{\left[2\Re\beta-\Re\alpha-2-{1\over2}\right]j}2^{\epsilon j}
 \\\\ \ds~~~~~~~~~~~~~~~~~~~~~~~~~~~~~~~~~~~~~~~~~~~~~~~~~~~~~~~~
 ~=~\C_{\Re\alpha~\Re\beta}~e^{\c|\Im\alpha|}~2^{-(1-\epsilon)j}2^{\left[2\Re\beta-\Re\alpha-1\right]j}2^{-j/2}. 
\end{array}
\eeq
Because of Euler's formulae, we replace the cosine function in (\ref{S principal}) with $e^{\pm2\pi\i \rho}$ multiplied by  $e^{\pm\pi\i \big[{1\over 4}-{\alpha\over 2}\big]}$ whose amplitude is controlled by $\C_{\Re\alpha}e^{\c|\Im\alpha|}$. 
Assert
\bel{A_j sharp}
\begin{array}{lr}\ds
{^\sharp}\Hat{\A}_{j~m}^{\alpha~\beta}(\xi)~=~
\Hat{\phi}(\xi)\int_{\lambda_{m-1}}^{\lambda_m}  e^{2\pi\i \big[|\xi|-1\big]r} \left({1\over r|\xi|}\right)^{\alpha-{1\over 2}}r^{2\beta-2}  dr
\\\\ \ds~~~~~~~~~~~~~
~=~|\xi|^{{1\over 2}-\alpha} \Hat{\phi}(\xi)\int_{\lambda_{m-1}}^{\lambda_m}  e^{2\pi\i \big[|\xi|-1\big]r} r^{2\beta-\alpha+{1\over 2}-2}  dr.
\end{array}
\eeq
Suppose $|1-|\xi||\leq2^{-j}$. We have
\bel{A_j sharp Est.1}
\begin{array}{lr}\ds
\left| {^\sharp}\Hat{\A}_{j~m}^{\alpha~\beta}(\xi)\right|~\leq~\C_{\Re\alpha} \int_{\lambda_{m-1}}^{\lambda_m}   r^{2\Re\beta-\Re\alpha+{1\over 2}-2}  dr
\\\\ \ds~~~~~~~~~~~~~~~~~
~\leq~\C_{\Re\alpha} ~2^{j/2} 2^{\left[2\Re\beta-\Re\alpha-2\right]j}2^{\epsilon j}
\\\\ \ds~~~~~~~~~~~~~~~~~
~=~\C_{\Re\alpha} ~2^{-(1-\epsilon)j}2^{j/2} 2^{\left[2\Re\beta-\Re\alpha-1\right]j}.
\end{array}
\eeq
On the other hand, suppose $|1-|\xi||>2^{-j}$. By integration by parts $w.r.t~r$ in (\ref{A_j sharp}), we find
\bel{A_j sharp by parts}
\begin{array}{lr}\ds
{^\sharp}\Hat{\A}_{j~m}^{\alpha~\beta}(\xi)~=~{1\over 2\pi\i} |\xi|^{{1\over 2}-\alpha}\Hat{\phi}(\xi) {1\over |\xi|-1} e^{2\pi\i \big[|\xi|-1\big]r} r^{2\beta-\alpha+{1\over 2}-2} \Bigg|^{\lambda_m}_{\lambda_{m-1}}
\\\\ \ds~~~~~~~~~~~~~~~
- {1\over 2\pi\i}\Big[\hbox{$2\beta-\alpha-{3\over 2}$}\Big] |\xi|^{{1\over 2}-\alpha}\Hat{\phi}(\xi) {1\over |\xi|-1}\int_{\lambda_{m-1}}^{\lambda_m}  e^{2\pi\i \big[|\xi|-1\big]r} r^{2\beta-\alpha-{1\over 2}-2}  dr.
\end{array}
\eeq
From (\ref{A_j sharp by parts}), by using the mean-value theorem on the boundary term, we have
\bel{A_j sharp Est.2}
\begin{array}{lr}\ds
\left| {^\sharp}\Hat{\A}_{j~m}^{\alpha~\beta}(\xi)\right|~\leq~\C_{\Re\alpha~\Re\beta}~e^{\c\Im\alpha}e^{\c\Im\beta}~ \left|{1\over 1-|\xi|}\right| 2^{j\big[2\Re\beta-\Re\alpha-{1\over 2}-2\big]}\left|\lambda_m-\lambda_{m-1}\right|
\\\\ \ds~~~~~~~~~~~~~~~~~
~+~\C_{\Re\alpha~\Re\beta}~e^{\c\Im\alpha}e^{\c\Im\beta} \left|{1\over 1-|\xi|}\right| \int_{\lambda_{m-1}}^{\lambda_m}   r^{2\Re\beta-\Re\alpha-{1\over 2}-2}  dr
\\\\ \ds~~~~~~~~~~~~~~~~~
~\leq~\C_{\Re\alpha~\Re\beta}~e^{\c\Im\alpha}e^{\c\Im\beta}~~ \left|{1\over 1-|\xi|}\right|
~2^{\left[2\Re\beta-\Re\alpha-{1\over 2}-2\right]j} 2^{\epsilon j}
\\\\ \ds~~~~~~~~~~~~~~~~~
~\leq~\C_{\Re\alpha~\Re\beta}~e^{\c\Im\alpha}e^{\c\Im\beta} ~2^{-(1-\epsilon)j}\left|{1\over 1-|\xi|}\right| ~2^{-j/2} 2^{\left[2\Re\beta-\Re\alpha-1\right]j}
\\\\ \ds~~~~~~~~~~~~~~~~~
~\leq~\C_{\Re\alpha~\Re\beta}~e^{\c\Im\alpha}e^{\c\Im\beta} ~2^{-(1-\epsilon)j}\left|{1\over 1-|\xi|}\right|^{1\over 2} ~ 2^{\left[2\Re\beta-\Re\alpha-1\right]j}.
\end{array}
\eeq
Our estimates in (\ref{A_j sharp Est.1})-(\ref{A_j sharp Est.2}) remain to be true if $e^{2\pi\i r|\xi|}$ in (\ref{A_j sharp}) is replaced by $e^{-2\pi\i r|\xi|}$. 
Together with (\ref{E est}), we conclude (\ref{Q_j EST sharp j <})-(\ref{Q_j EST sharp j >}). This finishes the proof of {\bf Lemma One}.

\subsection{Size of ${^\sharp}\Hat{\U}^{\alpha~\beta}_{j~m}(\xi)$ and ${^\sharp}\Hat{\V}^{\alpha~\beta}_{j~m}(\xi)$}
We begin to prove (\ref{Result One U}) and (\ref{Result One V}) in {\bf Proposition One}.
Recall (\ref{P sharp v}). We have 
\[
\begin{array}{lr}\ds
{^\sharp}\Hat{\P}_{j~m}^{\alpha~\beta~\nu}(\xi)~=~\varphi^\nu_j(\xi) \Hat{\phi}(\xi)\int_{\lambda_{m-1}\leq|r|<\lambda_m} e^{-2\pi\i  r}{^\sharp}\Hat{\Omega}^\alpha (r\xi)\omega(r) |r|^{2\beta-1} dr
  \\\\ \ds~~~~~~~~~~~~~~~
  ~=~\varphi^\nu_j(\xi)    {^\sharp}\Hat{\P}_{j~m}^{\alpha~\beta}(\xi)
\end{array}\]
where $\varphi^\nu_j$ is defined in (\ref{phi^v_j}). Moreover, 
\[{^\sharp}\Hat{\P}_{j~m}^{\alpha~\beta}(\xi)~=~\sum_{\nu\colon\xi^\nu_j\in\mathds{S}^{n-1}} {^\sharp}\Hat{\P}_{j~m}^{\alpha~\beta~\nu}(\xi)~=~  {^\sharp}\Hat{\P}_{j~m}^{\alpha~\beta}(\xi) \sum_{\nu\colon\xi^\nu_j\in\mathds{S}^{n-1}}   \varphi^\nu_j(\xi).   
\]
\begin{remark} For $\varphi^\nu_j$ defined in (\ref{phi^v_j}), it is easy to see
\[0\leq\varphi^\nu_j(\xi)\leq1,\qquad \sum_{\nu\colon\xi^\nu_j\in\mathds{S}^{n-1}} \varphi^\nu_j(\xi)=1.\]
\end{remark}
By applying (\ref{Result 1 > j})-(\ref{Result 1 > xi}) in {\bf Lemma One}, we find
\bel{Kernel l Sum Est}
\begin{array}{lr}\ds
\left|{^\sharp}\Hat{\P}_{j~m}^{\alpha~\beta}(\xi)\right|
~\leq~\C_{\Re \alpha~\Re\beta} ~e^{\c|\Im \alpha|}e^{\c|\Im\beta|} ~ 2^{-(1-\epsilon)j}2^{j/2}2^{-\ve j}.
\end{array}
\eeq

Recall
 $\Psi^{\alpha}_{j~m}, \Psi^{\alpha~\nu}_{j~m}$  defined in (\ref{Psi mu j}) and ${^\sharp}\U^{\alpha~\beta}_{j~m}$, ${^\sharp}\V^{\alpha~\beta}_{j~m}$ defined in (\ref{U sharp})-(\ref{V sharp}). 
 We find
\[
\begin{array}{lr}\ds
{^\sharp}\U^{\alpha~\beta}_{j~m}(x)+{^\sharp}\V^{\alpha~\beta}_{j~m}(x)~=~ \sum_{\nu\colon\xi^\nu_j\in\mathds{S}^{n-1}}{^\sharp}\U^{\alpha~\beta}_{j~m}(x)+{^\sharp}\V^{\alpha~\beta}_{j~m}(x)
\\\\ \ds~~~~~~~~~~~~~~~~~~~~~~~~~~~~~~~
~=~\sum_{\nu\colon\xi^\nu_j\in\mathds{S}^{n-1}}
{^\sharp}\P^{\alpha~\beta~\nu}_{j~m}(x) \left[1-\Psi^{\alpha~\nu}_{j~m}(x)\right] +
{^\sharp}\P^{\alpha~\beta~\nu}_{j~m}(x)  \Psi^{\alpha~\nu}_{j~m}(x)
\\\\ \ds~~~~~~~~~~~~~~~~~~~~~~~~~~~~~~~
~=~\sum_{\nu\colon\xi^\nu_j\in\mathds{S}^{n-1}}
{^\sharp}\P^{\alpha~\beta~\nu}_{j~m}(x)
~=~{^\sharp}\P^{\alpha~\beta}_{j~m}(x).
\end{array}
\]
Hence, (\ref{Result One U}) can be obtained by putting together (\ref{Result One V}) and  
 (\ref{Kernel l Sum Est}). 
 
Our task is left to show
\bel{V Est ve}
\left| {^\sharp}\Hat{\V}_{j~m}^{\alpha~\beta}(\xi)\right|~\leq~\C_{\Re \alpha~\Re\beta} ~e^{\c|\Im \alpha|}e^{\c|\Im\beta|} ~2^{-(1-\epsilon)j}2^{-j\ve},\qquad \ve=\ve(\Re\alpha,\Re\beta)>0.
\eeq
Let ${^\sharp}\V^{\alpha~\beta}_{j~m}$, ${^\sharp}\V^{\alpha~\beta~\nu}_{j~m}$ defined in (\ref{V sharp}). We have
\bel{V sharp Transform}
\begin{array}{lr}\ds
{^\sharp}\Hat{\V}^{\alpha~\beta}_{j~m}(\xi)~=~\int_{\R^n}e^{-2\pi\i x\cdot\xi}\left\{\sum_{\nu\colon\xi^\nu_j\in\mathds{S}^{n-1}} \Psi^{\alpha~\nu}_{j~m}(x){^\sharp}\P^{\alpha~\beta~\nu}_{j~m}(x)\right\}dx
\\\\ \ds~~~~~~~~~~~~~
~=~\sum_{\nu\colon\xi^\nu_j\in\mathds{S}^{n-1}}\int_{\R^n} \Hat{\Psi}^{\alpha~\nu}_{j~m}(\xi-\zeta)
{^\sharp}\Hat{\P}^{\alpha~\beta~\nu}_{j~m}(\zeta) d\zeta.
\end{array}
\eeq
On the other hand, for $\Psi^{\alpha~\nu}_{j~m}$ defined in (\ref{Psi mu j}), we find
\bel{Psi mu Transform}
\begin{array}{lr}\ds
\Hat{\Psi}^{\alpha~\nu}_{j~m}(\xi-\zeta)~=~
\int_{\R^n} e^{2\pi\i x\cdot(\xi-\zeta)} \varphi\left[ 2^{-\epsilon j-1}\left|\lambda_m-\left(\L_\nu^T x\right)_1\right|\right]
\prod_{i=2}^n \varphi\left[2^{-\left[{1\over 2}+\epsilon\right]j}\left|\left(\L_\nu^T x\right)_i\right|\right] dx
\\ \ds~~~~~~~~~~~~~~~~~~~
~+~\int_{\R^n} e^{2\pi\i x\cdot(\xi-\zeta)}\varphi\left[ 2^{-\epsilon j-1}\left|\lambda_m+\left(\L_\nu^T x\right)_1\right|\right]
\prod_{i=2}^n \varphi\left[2^{-\left[{1\over 2}+\epsilon\right]j}\left|\left(\L_\nu^T x\right)_i\right|\right] dx.
\end{array}
\eeq
$\L_\nu$ is an $n\times n$-orthogonal matrix  with $\det\L_\nu=1$.  Moreover,  $\L^T_\nu\xi^\nu_j=(1,0)^T\in\R\times\R^{n-1}$. 
Write $\zeta=\L_\nu \eta$ and $x=\L_\nu u$ inside (\ref{V sharp Transform})-(\ref{Psi mu Transform}). 
We have
\bel{Psi mu j Transform}
\begin{array}{lr}\ds
\Hat{\Psi}^{\alpha~\nu}_{j~m}(\xi-\L_\mu\eta)
~=~
\int_{\R^n}e^{-2\pi\i \big[\L_\nu^T\xi-\eta\big]\cdot  u} \varphi\left[ 2^{-\epsilon j-1}\left|\lambda_m-u_1\right|\right]
\prod_{i=2}^n \varphi\left[2^{-\left[{1\over 2}+\epsilon\right]j}\left|u_i\right|\right] du
\\ \ds~~~~~~~~~~~~~~~~~~~~~~~~
~+~\int_{\R^n}e^{-2\pi\i \big[\L_\nu^T\xi-\eta\big]\cdot  u} \varphi\left[ 2^{-\epsilon j-1}\left|\lambda_m+u_1\right|\right]
\prod_{i=2}^n \varphi\left[2^{-\left[{1\over 2}+\epsilon\right]j}\left|u_i\right|\right] du
\end{array}
\eeq
and
\[\begin{array}{lr}\ds
{^\sharp}\Hat{\P}_{j~m}^{\alpha~\beta~\nu}(\L_\nu \eta)~=~\varphi^\nu_j(\L_\nu \eta) \Hat{\phi}(\L_\nu \eta)\int_{\lambda_{m-1}\leq|r|<\lambda_m} e^{-2\pi\i  r}~{^\sharp}\Hat{\Omega}^\alpha (r\L_\nu \eta)\omega(r) |r|^{2\beta-1} dr.
\end{array}\]
Recall $\varphi\in\mathcal{C}^\infty_o(\R)$ such that $\varphi(t)=1$ if $|t|\leq1$ and $\varphi(t)=0$ for $|t|>2$.  We find 
\bel{Supp varphi+}
\begin{array}{lr}\ds
\vol~\supp \left\{~\Bigg[\varphi\left[ 2^{-\epsilon j-1}\left|\lambda_m-u_1\right|\right]+\varphi\left[ 2^{-\epsilon j-1}\left|\lambda_m+u_1\right|\right]\Bigg]\prod_{i=2}^n \varphi\left[2^{-\left[{1\over 2}+\epsilon\right]j}\left|u_i\right|\right]\right\}
\\\\ \ds
\lesssim~2^{(n-1)\big[{1\over 2}+\epsilon\big]j} 2^{\epsilon j} ~=~2^{\left({n-1\over 2}\right)j}2^{ n\epsilon j}.
\end{array}
\eeq
From (\ref{Psi mu j Transform}) and (\ref{Supp varphi+}), we have
\bel{a priori}
\left| \Hat{\Psi}^{\alpha~\nu}_{j~m}(\xi-\L_\mu\eta)\right|~\lesssim~2^{\left({n-1\over 2}\right)j}2^{ n\epsilon j}.
\eeq
A direct computation shows
\bel{varphi j EST}
\begin{array}{lr}\ds
\left| \partial_{u_1}^N \varphi\left[ 2^{-\epsilon j-2}\left|\lambda_m- u_1\right|\right]\right| 
~\leq~\C_N~2^{-N\epsilon j},
\\\\ \ds
\left| \partial_{u_i}^N \varphi\left[2^{-\left[{1\over 2}+\epsilon\right]j}\left|u_i\right|\right]\right|
~\leq~\C_N~2^{-N\big[{1\over 2}+\epsilon\big]j},\qquad i=2,\ldots,n
\end{array}
\eeq
for every $N\ge0$.

{\bf Case 1}~~Suppose $|\xi|\leq {1\over 10}$ or $|\xi|>10$. For ${1\over 3}<|\eta|\leq3$, we either have  $\left| \big[\L_\mu^T\xi-\eta\big]_1\right|>{1\over 5\sqrt{n}}$ or $\left| \big[\L_\mu^T\xi-\eta\big]_i\right|>{1\over 5\sqrt{n}}$ for some $i=2,\ldots,n$.

By integration by parts $w.r.t~u_1$ or $u_i$ inside (\ref{Psi mu j Transform}) and using (\ref{Supp varphi+})-(\ref{varphi j EST}), we find
\[
\begin{array}{lr}\ds
\left|\Hat{\Psi}^{\alpha~\nu}_{j~m}(\xi-\L_\mu\eta)\right|~\leq~\C_N ~2^{\left({n-1\over 2}\right)j}2^{ n\epsilon j}~  2^{- N\sigma j},\qquad\hbox{\small{$N\ge0$}}.
\end{array}
\]
By choosing $N$ sufficiently large depending on $\epsilon=\epsilon(\Re\alpha,\Re\beta)>0$, we obtain
\bel{Psi mu j Transform norm by parts'}
\begin{array}{lr}\ds
\left|\Hat{\Psi}^{\alpha~\nu}_{j~m}(\xi-\L_\mu\eta)\right|~\leq~\C_{\Re\alpha~\Re\beta}.
\end{array}
\eeq
 
{\bf Case 2}~~Let ${1\over 10}<|\xi|\leq 10$. We must have 
\[\xi~\in~ \Gamma^{\mu}_j~=~\Bigg\{\xi\in\R^n\colon \left|{\xi\over |\xi|}-\xi^{\mu}_{j}\right|<2^{-j/2+1}\Bigg\}\qquad\hbox{ \small{ for some}}~~ \xi^\mu_j\in\mathds{S}^{n-1}.\] 
\v
\begin{lemma}   Let ${1\over 10}<|\xi|\leq 10$ and ${1\over 3}<|\eta|\leq3$. Suppose $\left|\xi^\mu_j-\xi^\nu_j\right|\ge \c 2^{-j/2}$ for some $\c>0$ large. Either we have $\left|\Big[\L_\nu^T\xi-\eta\Big]_1\right|\approx2$ or $\left|\Big[\L_\nu^T\xi-\eta\Big]_\imath\right|\gtrsim2^{-j/2}$ for some  $\imath\in\{2,\ldots,n\}$ whenever $\L_\nu \eta\in \Gamma^\nu_j$.
\end{lemma}

{\bf Proof}~~Note that $\L^T_\nu\xi^\nu_j=(1,0)^T\in\R\times\R^{n-1}$. Consider $\L_\nu \eta\in \Gamma^\nu_j$.
We have 
$\left| \L_\nu^T \xi^\nu_j-(1,0)^T\right|=0$.
 By using the triangle inequality, we find
\bel{eq size est1}
\begin{array}{lr}\ds
{|\eta_i|\over|\eta|}~\leq~\left|{\eta\over|\eta|}-(1,0)^T\right|~\leq~\left|{\eta\over|\eta|}-\L_\nu^T \xi^\nu_j\right|+\left|\L_\nu^T \xi^\nu_j-(1,0)^T\right|
\\\\ \ds~~~~~~~~~~~~~~~~~~~~~~~~~~~~~~~~~~~
~\leq ~2^{-j/2+1},\qquad i=2,\ldots,n.
\end{array}
\eeq
On the other hand, recall $\left|\xi^\mu_j-\xi^\nu_j\right|\ge\c2^{-j/2}$ for some $\c$ large. 
Let $\xi\in\Gamma^{\mu}_j$. By using the triangle inequality, we find
\bel{eq size est2}
\begin{array}{lr}\ds
\left| {\L_\nu^T\xi\over |\xi|}- (1,0)^T\right|~=~\left|{\xi\over |\xi|}- \xi^\nu_j\right|
\\\\ \ds~~~~~~~~~~~~~~~~~~~~~~~~
~\ge~\left|\xi^{\mu}_j-\xi^\nu_j\right|-\left| {\xi\over |\xi|}- \xi^{\mu}_j\right|
\\\\ \ds~~~~~~~~~~~~~~~~~~~~~~~~
~\ge~\c 2^{-j/2}-2^{-j/2+1}     
\\\\ \ds~~~~~~~~~~~~~~~~~~~~~~~~
~=~\big[\c-2\big]~ 2^{-j/2}.
\end{array}
\eeq
Suppose $\left| \left( {\L_\nu^T\xi\over |\xi|}\right)_1-1\right|\leq\left[\sum_{i=2}^n \left( {\L_\nu^T\xi\over |\xi|}\right)_i^2\right]^{1\over2}$. The estimate in (\ref{eq size est2})
 further implies
\bel{eq size est3}
\begin{array}{lr}\ds
\left| {\left(\L_\nu^T\xi\right)_\imath\over |\xi|}\right|~\ge~{1\over \sqrt{2}\sqrt{n-1}} \left| {\L_\nu^T\xi\over |\xi|}- (1,0)^T\right|
\\\\ \ds~~~~~~~~~~~~~~
~\ge~{1\over \sqrt{2}\sqrt{n-1}} \big[\c-2\big]~ 2^{-j/2}\qquad \hbox{\small{for some}}\qquad \imath\in\{2,\ldots,n\}.
\end{array}
\eeq
Suppose $\left| \left( {\L_\nu^T\xi\over |\xi|}\right)_1-1\right|>\left[\sum_{i=2}^n \left( {\L_\nu^T\xi\over |\xi|}\right)_i^2\right]^{1\over2}$. As a geometric fact, ${\L_\nu^T\xi\over |\xi|}$ belongs to the semi-sphere opposite to $(1,0)^T$. We either have $\left( {\L_\nu^T\xi\over |\xi|}\right)_1\approx-1$ or $\left| \left( {\L_\nu^T\xi\over |\xi|}\right)_\imath\right|>{1\over \sqrt{2}\sqrt{n-1}} \big[\c-2\big]~ 2^{-j/2}$ for some $\imath\in\{2,\ldots,n\}$ as (\ref{eq size est3}).

By taking into account ${1\over 3}<|\eta|\leq3$ and ${1\over 10}<|\xi|\leq 10$ in (\ref{eq size est1}) and (\ref{eq size est2})-(\ref{eq size est3}), we find
\[
|\eta_\imath|~\leq~3~ 2^{-j/2+1},\qquad \left| \left(\L_\nu^T\xi\right)_\imath\right|~\ge~{1\over 10\sqrt{2}\sqrt{n-1}} \big[\c-2\big]~ 2^{-j/2}.
\]
Consequently, $\ds\left|\big[\L_\nu^T\xi-\eta\big]_\imath\right|\gtrsim2^{-j/2}$  provided that $\c$ is  large.
\endproof
\v
Let $\xi\in\Gamma^{\mu}_j$ and $\left|\xi^\mu_j-\xi^\nu_j\right|\ge\c2^{-j/2}$ for some $\c$ large. By applying {\bf Lemma 5.1} and  using the estimates in (\ref{Supp varphi+})-(\ref{varphi j EST}), an $N$-fold integration by parts either $w.r.t~u_1$ or $u_\imath$ inside (\ref{Psi mu j Transform}) shows
\bel{Psi mu j Transform norm by parts}
\begin{array}{lr}\ds
\left|\Hat{\Psi}^{\alpha~\mu}_{j~m}(\xi-\L_\mu\eta)\right|~\leq~\C_N ~2^{\left({n-1\over 2}\right)j}2^{ n\epsilon j}~\left[2^{-N\sigma j}+ 2^{N\big({1\over 2}\big)j} 2^{-N\big[{1\over2}+\epsilon\big] j}\right]
\\\\ \ds~~~~~~~~~~~~~~~~~~~~~~~~~~
~\leq~\C_N ~2^{\left({n-1\over 2}\right)j}2^{ n\epsilon j}~  2^{-N\epsilon  j}
\\\\ \ds~~~~~~~~~~~~~~~~~~~~~~~~~~
~\leq~\C_{\Re\alpha~\Re\beta}
\end{array}
\eeq
for $N$ sufficiently large depending on $\epsilon=\epsilon(\Re\alpha,\Re\beta)>0$.

Now, we assert
\bel{Sum S_k}
\begin{array}{cc}\ds
\int_{\R^n} \Hat{\Psi}^{\alpha~\nu}_{j~m}(\xi-\L_\mu\eta)
{^\sharp}\Hat{\P}^{\alpha~\beta~\nu}_{j~m}(\L_\mu\eta) d\eta
~=~\sum_{k=0}^\infty \int_{\mathcal{S}_k} \Hat{\Psi}^{\alpha~\nu}_{j~m}(\xi-\L_\mu\eta)
{^\sharp}\Hat{\P}^{\alpha~\beta~\nu}_{j~m}(\L_\mu\eta) d\eta,
\\\\ \ds
\mathcal{S}_k~=~\Bigg\{\eta\in\R^n\colon 2^{-k}\leq|1-|\eta||<2^{-k+1}\Bigg\}.
\end{array}
\eeq
Observe that
\bel{Supp}
\vol\Bigg\{\mathcal{S}_k\cap\supp {^\sharp}\Hat{\P}^{\alpha~\beta~\nu}_{j~m} \Bigg\}~\lesssim~2^{-\left({n-1\over 2}\right)j} 2^{-k}.
\eeq
Recall ${^\sharp}\Hat{\P}^{\alpha~\beta~\nu}_{j~m}(\L_\mu \eta)=\varphi^\nu_j(\L_\mu \eta) {^\sharp}\Hat{\P}^{\alpha~\beta}_{j~m}(\L_\mu\eta)$ where $0\leq\varphi^\nu_j(\L_\mu \eta)\leq1$. See {\bf Remark 5.2}. 

For $j\leq k$, we have $\left|1-|\eta|\right|\leq 2^{-k+1}\leq 2^{-j+1}$. 
By applying (\ref{Result 1 > j})  in 
{\bf Lemma One}, we find
\bel{P sharp mu j >}
\begin{array}{lr}\ds
\left|{^\sharp}\Hat{\P}^{\alpha~\beta~\nu}_{j~m}(\L_\mu \eta)\right|~\leq~\C_{\Re\alpha~\Re\beta} ~e^{\c|\Im\alpha|}~2^{-(1-\epsilon)j}2^{j/2}2^{-\ve j}
\\\\ \ds~~~~~~~~~~~~~~~~~~~~~~
~\leq~\C_{\Re\alpha~\Re\beta} ~e^{\c|\Im\alpha|}~2^{-(1-\epsilon)j}2^{k/2}2^{-\ve j}.
\end{array}
\eeq
For $j>k$, we have $\left|1-|\eta|\right|>2^{-k}>2^{-j}$. By applying
(\ref{Result 1 > xi}) in {\bf Lemma One}, we find
\bel{P sharp mu xi <}
\begin{array}{lr}\ds
\left|{^\sharp}\Hat{\P}^{\alpha~\beta~\nu}_j(\L_\mu \eta)\right|~\leq~\C_{\Re\alpha~\Re\beta} ~e^{\c|\Im\alpha|}e^{\c|\Im\beta|}~2^{-(1-\epsilon)j}\left|{1\over 1-|\eta|}\right|^{1\over 2} 2^{-\ve j}
\\\\ \ds~~~~~~~~~~~~~~~~~~~~~~
~\leq~\C_{\Re\alpha~\Re\beta} ~e^{\c|\Im\alpha|}e^{\c|\Im\beta|}~2^{-(1-\epsilon)j}2^{k/2}2^{-\ve j}.
\end{array}
\eeq
From (\ref{Supp}) and (\ref{P sharp mu j >})-(\ref{P sharp mu xi <}), we obtain
\bel{P sharp Int EST}
\begin{array}{lr}\ds
\sum_{k=0}^\infty \int_{\mathcal{S}_k} \left|
{^\sharp}\Hat{\P}^{\alpha~\beta~\nu}_{j~m}(\L_\mu\eta)\right| d\eta
\\\\ \ds
~\leq~\sum_{k=0}^\infty~\left|{^\sharp}\Hat{\P}^{\alpha~\beta~\nu}_j(\L_\mu \eta)\right| \vol\Bigg\{\mathcal{S}_k\cap\supp {^\sharp}\Hat{\P}^{\alpha~\beta~\nu}_{j~m} \Bigg\}
\\\\ \ds
~\leq~\sum_{k=0}^\infty~\C_{\Re\alpha~\Re\beta} ~e^{\c|\Im\alpha|}e^{\c|\Im\beta|}  2^{-(1-\epsilon)j}2^{k/2}2^{-\ve j} ~2^{-\left({n-1\over 2}\right)j} 2^{-k}
\\\\ \ds
~=~\C_{\Re\alpha~\Re\beta} ~e^{\c|\Im\alpha|}e^{\c|\Im\beta|} ~2^{-(1-\epsilon)j}2^{-\left({n-1\over 2}\right)j}2^{-\ve j} \sum_{k=0}^\infty 2^{-k/2}
\\\\ \ds
~\leq~\C_{\Re\alpha~\Re\beta} ~e^{\c|\Im\alpha|}e^{\c|\Im\beta|}~2^{-(1-\epsilon)j}2^{-\left({n-1\over 2}\right)j}2^{-\ve j}.
\end{array}
\eeq
Recall {\bf Remark 1.3}.  There are at most a constant multiple of $2^{\big({n-1\over 2}\big) j}$ many $\xi^\nu_j\in\mathds{S}^{n-1}$.

Consider $|\xi|\leq{1\over 10}$ or $|\xi|>10$. From (\ref{V sharp Transform}) and (\ref{Sum S_k}), we have
\bel{V sharp Transform bound Case1}
\begin{array}{lr}\ds
\left|{^\sharp}\Hat{\V}^{\alpha~\beta}_{j~m}(\xi)\right|~\leq~\sum_{\nu\colon\xi^\nu_j\in\mathds{Z}^{n-1}} \sum_{k=0}^\infty \left|\int_{\mathcal{S}_k} \Hat{\Psi}^{\alpha~\nu}_{j~m}(\xi-\L_\nu\eta)
{^\sharp}\Hat{\P}^{\alpha~\beta~\nu}_{j~m}(\L_\nu\eta) d\eta\right|
\\\\ \ds~~~~~~~~~~~~~~~~~
~\leq~\C_{\Re\alpha~\Re\beta}~\sum_{\nu\colon\xi^\nu_j\in\mathds{S}^{n-1}} \sum_{k=0}^\infty \int_{\mathcal{S}_k} \left|
{^\sharp}\Hat{\P}^{\alpha~\beta~\nu}_{j~m}(\L_\mu\eta)\right| d\eta\qquad \hbox{\small{by (\ref{Psi mu j Transform norm by parts'})}}
\\\\ \ds~~~~~~~~~~~~~~~~
~\leq~\C_{\Re\alpha~\Re\beta} ~e^{\c|\Im\alpha|}e^{\c|\Im\beta|}~\sum_{\nu\colon\xi^\nu_j\in\mathds{S}^{n-1}}~2^{-(1-\epsilon)j}2^{-\left({n-1\over 2}\right)j}2^{-\ve j}\qquad\hbox{\small{by (\ref{P sharp Int EST})}}
\\\\ \ds~~~~~~~~~~~~~~~~
~\leq~\C_{\Re\alpha~\Re\beta} ~e^{\c|\Im\alpha|}e^{\c|\Im\beta|}~  2^{\left({n-1\over 2}\right)j} ~2^{-(1-\epsilon)j}2^{-\left({n-1\over 2}\right)j}2^{-\ve j}\qquad\hbox{\small{by {\bf Remark 1.3}}}
\\\\ \ds~~~~~~~~~~~~~~~~
~=~\C_{\Re\alpha~\Re\beta} ~e^{\c|\Im\alpha|}e^{\c|\Im\beta|} ~2^{-(1-\epsilon)j}2^{-\ve j}.
\end{array}
\eeq
Let ${1\over 10}<|\xi|\leq10$ and $\xi\in\Gamma^\mu_j$ for some $\xi^\mu_j\in\mathds{S}^{n-1}$. From (\ref{V sharp Transform}) and (\ref{Sum S_k}), we have
\bel{V sharp Transform bound Case2}
\begin{array}{lr}\ds
\left|{^\sharp}\Hat{\V}^{\alpha~\beta}_{j~m}(\xi)\right|~\leq~\sum_{\nu\colon\xi^\nu_j\in\mathds{S}^{n-1}}\sum_{k=0}^\infty \left|\int_{\mathcal{S}_k} \Hat{\Psi}^{\alpha~\nu}_{j~m}(\xi-\L_\nu\eta)
{^\sharp}\Hat{\P}^{\alpha~\beta~\nu}_{j~m}(\L_\nu\eta) d\eta\right|
\\\\ \ds
=~\Bigg[\sum_{\nu\colon\xi^\nu_j\in\mathds{S}^{n-1}, \left|\xi^\mu_j-\xi^\nu_j\right|\ge\c2^{-j/2}}+\sum_{\nu\colon\xi^\nu_j\in\mathds{S}^{n-1}, \left|\xi^\mu_j-\xi^\nu_j\right|\leq\c2^{-j/2}}\Bigg]\sum_{k=0}^\infty \left|\int_{\mathcal{S}_k} \Hat{\Psi}^{\alpha~\nu}_{j~m}(\xi-\L_\nu\eta)
{^\sharp}\Hat{\P}^{\alpha~\beta~\nu}_{j~m}(\L_\nu\eta) d\eta\right|~~~\hbox{\small{for some $\c$ large}}
\\\\ \ds
\leq~\C_{\Re\alpha~\Re\beta}\sum_{\nu\colon\xi^\nu_j\in\mathds{S}^{n-1}, \left|\xi^\mu_j-\xi^\nu_j\right|\ge\c2^{-j/2}}\sum_{k=0}^\infty \int_{\mathcal{S}_k} 
\left|{^\sharp}\Hat{\P}^{\alpha~\beta~\nu}_{j~m}(\L_\nu\eta) \right|d\eta\qquad\hbox{\small{by (\ref{Psi mu j Transform norm by parts})}}
\\ \ds
+~\C \sum_{\nu\colon\xi^\nu_j\in\mathds{S}^{n-1}, \left|\xi^\mu_j-\xi^\nu_j\right|\leq\c2^{-j/2}} 2^{\left({n-1\over 2}\right)j}2^{n\epsilon j}  \sum_{k=0}^\infty \int_{\mathcal{S}_k} 
\left|{^\sharp}\Hat{\P}^{\alpha~\beta~\nu}_{j~m}(\L_\nu\eta) \right|d\eta\qquad\hbox{\small{by (\ref{a priori})}}
\\\\ \ds
\leq~\C_{\Re\alpha~\Re\beta} ~e^{\c|\Im\alpha|}e^{\c|\Im\beta|}\sum_{\nu\colon\xi^\nu_j\in\mathds{S}^{n-1}, \left|\xi^\mu_j-\xi^\nu_j\right|\ge\c2^{-j/2}} 2^{-(1-\epsilon)j}2^{-\left({n-1\over 2}\right)j}2^{-\ve j}
\\ \ds
+~\C_{\Re\alpha~\Re\beta} ~e^{\c|\Im\alpha|}e^{\c|\Im\beta|}\sum_{\nu\colon\xi^\nu_j\in\mathds{S}^{n-1}, \left|\xi^\mu_j-\xi^\nu_j\right|\leq\c2^{-j/2}} 2^{\left({n-1\over 2}\right)j}2^{n\epsilon j}~ 2^{-(1-\epsilon)j}2^{-\left({n-1\over 2}\right)j}2^{-\ve j} \qquad\hbox{\small{by (\ref{P sharp Int EST})}}
\\\\ \ds
\leq~\C_{\Re\alpha~\Re\beta} ~e^{\c|\Im\alpha|}e^{\c|\Im\beta|}  ~2^{\big({n-1\over 2}\big)j}  ~2^{-(1-\epsilon)j} 2^{-\left({n-1\over 2}\right)j} 2^{-\ve j}
\\ \ds
+~\C_{\Re\alpha~\Re\beta} e^{\c|\Im\alpha|}e^{\c|\Im\beta|}  
~2^{\left({n-1\over 2}\right)j}2^{n\epsilon j}  ~2^{-(1-\epsilon)j}2^{-\left({n-1\over 2}\right)j}2^{-\ve j} \qquad \hbox{\small{by {\bf Remark 1.3}}}
\\\\ \ds
\leq~\C_{\Re\alpha~\Re\beta} e^{\c|\Im\alpha|}e^{\c|\Im\beta|} ~2^{-(1-\epsilon)j}	2^{n\epsilon j} 2^{-\ve j}. 
\end{array}
\eeq
We obtain (\ref{V Est ve}) provided that $n\sigma\leq{1\over2}\ve$ as discussed in {\bf Remark 5.1}.

\section{On the $\L^1$-norm of ${^\flat}\U^{\alpha~\beta}_{j~m}$ and ${^\flat}\V^{\alpha~\beta}_{j~m}$}
\setcounter{equation}{0}
Recall ${^\flat}\Omega^\alpha$ defined in (\ref{Omega flat})-(\ref{Omega flat Transform}) and 
${^\flat}\Hat{\P}^{\alpha~\beta}_{j~m}$ defined in (\ref{P flat v}) for $\Re\alpha>0$ and $0<\Re\beta<{1\over 2}$.
Given $j>0$ and $m=1,2,\ldots,M$, we have
\[
\begin{array}{cc}\ds
{^\flat}\P^{\alpha~\beta}_{j~m}(x)~=~\sum_{\nu\colon\xi^\nu_j\in\mathds{S}^{n-1}} {^\flat}\P^{\alpha~\beta~\nu}_{j~m}(x),
\\\\ \ds
{^\flat}\P^{\alpha~\beta~\nu}_{j~m}(x)~=~\int_{\R^n}e^{2\pi\i x\cdot\xi} \left\{\int_{\lambda_{m-1}\leq|r|<\lambda_m} e^{-2\pi\i  r} \varphi^\nu_j(\xi)\Hat{\phi}(\xi){^\flat}\Hat{\Omega}^\alpha (r\xi) \omega(r) |r|^{2\beta-1} dr\right\} d\xi,
\\\\ \ds
{^\flat}\Hat{\Omega}^\alpha (\xi)~=~\left({1\over|\xi|}\right)^{{n-1\over 2}+\Re\alpha-{1\over2}} \J_{{n-1\over 2}+\alpha-{1\over2}}\Big(2\pi|\xi|\Big).
\end{array}
\]
 $\lambda_m\in[2^{j-1},2^j]$:  $\lambda_0=2^{j-1}$, $\lambda_M=2^j$ and
$2^{\epsilon j-1}\leq\lambda_m-\lambda_{m-1}<2^{\epsilon j}$ for some $0<\sigma=\sigma(\Re\alpha)<{1\over 2}$ which can be chosen sufficiently small.  

$\varphi^\nu_j$ is defined in (\ref{phi^v_j}) whose support is contained in
$\Gamma^\nu_j=\Bigg\{\xi\in\R^n\colon \left|{\xi\over |\xi|}-\xi^{\nu}_{j}\right|<2^{-j/2+1}\Bigg\}$.

 Recall {\bf Remark 1.3}.
There are at most  a constant multiple of $2^{\big({n-1\over 2}\big)j}$ many $\xi^\nu_j\in\mathds{S}^{n-1}$ .

$\Hat{\phi}$ is defined in  (\ref{hat phi})  and $\supp\Hat{\phi}\subset\left\{\xi\in\R^n\colon{1\over 3}<|\xi|\leq3\right\}$. Moreover, 
 $\omega(r)$  given   in (\ref{I alpha rewrite}) satisfies 
$|\omega(r)|\lesssim \Big[1+|r|\Big]^{-1}$ as shown in (\ref{omega norm}).

By using the norm estimate of Bessel functions in (\ref{J norm}), we find 
\bel{Omega flat Transform norm}
\begin{array}{lr}\ds
\left| {^\flat}\Hat{\Omega}^\alpha (r\xi)\right|~=~\left({1\over|r\xi|}\right)^{{n-1\over 2}+\Re\alpha-{1\over2}} \left|\J_{{n-1\over 2}+\alpha-{1\over2}}\Big(2\pi |r\xi|\Big)\right|
\\\\ \ds~~~~~~~~~~~~~~~
~\leq~\C_{\Re\alpha}~e^{\c|\Im\alpha|}~\left({1\over 1+|r\xi|}\right)^{{n-1\over 2}+\Re\alpha}.
\end{array}
\eeq
Fubini's theorem allows us to write
\bel{P_j flat rewrite}
\begin{array}{cc}\ds
{^\flat}\P^{\alpha~\beta}_{j~m}(x)~=~ \int_{\lambda_{m-1}\leq|r|<\lambda_m} e^{-2\pi\i  r}~{^\flat}\P^{\alpha}_{r~j}(x) \omega(r) |r|^{2\beta-1} dr
\end{array}
\eeq
for which 
\bel{P_j v flat rewrite}
\begin{array}{cc}\ds
{^\flat}\P^{\alpha}_{r~j}(x)~=~\sum_{\nu\colon\xi^\nu_j\in\mathds{S}^{n-1}} {^\flat}\P^{\alpha~\nu}_{r~j}(x)
\\\\ \ds
{^\flat}\P^{\alpha~\nu}_{r~j}(x)~=~\int_{\R^n}e^{2\pi\i x\cdot\xi}  \varphi^\nu_j(\xi)\Hat{\phi}(\xi){^\flat}\Hat{\Omega}^\alpha (r\xi) d\xi.
\end{array}
\eeq

Recall $\Psi^{\alpha}_{j~m}$, $\Psi^{\alpha~\nu}_{j~m}$  defined in (\ref{Psi mu j}) and ${^\flat}\U^{\alpha~\beta}_{j~m},~{^\flat}\V^{\alpha~\beta}_{j~m}$ defined in (\ref{U flat})-(\ref{V flat}). We have
\[\begin{array}{cc}\ds
{^\flat}\U^{\alpha~\beta}_{j~m}(x)~=~\sum_{\nu\colon\xi^\nu_j\in\mathds{S}^{n-1}} {^\flat}\U^{\alpha~\beta~\nu}_{j~m}(x),\qquad 
{^\flat}\U^{\alpha~\beta~\nu}_{j~m}(x)~=~
{^\flat}\P^{\alpha~\beta~\nu}_{j~m}(x) \left[1-\Psi^{\alpha~\nu}_{j~m}(x)\right], 
\\\\ \ds
{^\flat}\V^{\alpha~\beta}_{j~m}(x)~=~\sum_{\nu\colon\xi^\nu_j\in\mathds{S}^{n-1}} {^\flat}\V^{\alpha~\beta~\nu}_{j~m}(x),\qquad 
{^\flat}\V^{\alpha~\beta~\nu}_{j~m}(x)~=~
{^\flat}\P^{\alpha~\beta~\nu}_{j~m}(x) \Psi^{\alpha~\nu}_{j~m}(x).
 \end{array}
\]
In order to prove (\ref{Result Two U})-(\ref{Result Two V}) in {\bf Proposition Two}, it is suffice to show
\bel{Result Three U r}
\begin{array}{lr}\ds
\sum_{\nu\colon\xi^\nu_j\in\mathds{S}^{n-1}} \int_{\R^n}\left|{^\flat}\P^{\alpha~\nu}_{r~j}(x)\right| \left[1-\Psi^{\alpha~\nu}_{j~m}(x)\right]dx 
~\leq~\C_{N~\Re \alpha} ~e^{\c|\Im \alpha|}~ 2^{-N j} 2^{-\ve j},\qquad \hbox{\small{$N\ge0$}}
\end{array}
\eeq
and
\bel{Result Three V r}
\sum_{\nu\colon\xi^\nu_j\in\mathds{S}^{n-1}}  \int_{\R^n}\Psi^{\alpha~\nu}_{j~m}(x)\left|{^\flat}\P^{\alpha~\nu}_{r~j}(x)\right|dx~\leq~\C_{\Re \alpha} ~e^{\c|\Im \alpha|}~  2^{-\ve j}
\eeq
for some $\ve=\ve(\Re\alpha)>0$ whenever $\lambda_{m-1}\leq|r|<\lambda_m$.

This is an immediate consequence of using Minkowski integral inequality together with $|\omega(r)|\leq\C \Big[1+|r|\Big]^{-1}$, $0<\Re\beta<{1\over 2}$ and $\left|\lambda_m-\lambda_{m-1}\right|\leq2^{\epsilon j}$.

Let ${^\flat}\P^{\alpha~\nu}_{r~j}$  defined in (\ref{P_j v flat rewrite}). We write
\bel{P_r Sum}
\begin{array}{lr}\ds
{^\flat}\P^{\alpha~\nu}_{r~j}(x)
~=~\int_{\R^n}e^{2\pi\i x\cdot\xi}  \varphi^\nu_j(\xi)\Hat{\phi}(\xi)\left({1\over|r\xi|}\right)^{{n-1\over 2}+\alpha-{1\over2}} \J_{{n-1\over 2}+\alpha-{1\over2}}\Big(2\pi|r\xi|\Big) d\xi
\\\\ \ds~~~~~~~~~~~~
~=~{^\flat}\Q^{\alpha~\nu}_{r~j}(x)+\sum_{k=1}^L {^\flat}\Q^{\alpha~\nu}_{r~j~k}(x)+{^\flat}\hbox{\bf R}^{\alpha~\nu}_{r~j~L}(x)
\end{array}
\eeq
by using
the asymptotic expansion of Bessel functions  in (\ref{J asymptotic})-(\ref{J O}).

First,  
\bel{Q_r v}
\begin{array}{lr}\ds
{^\flat}\Q^{\alpha~\nu}_{r~j}(x)~=~{1\over \pi}\int_{\R^n}e^{2\pi\i x\cdot\xi} \varphi^\nu_j(\xi)\Hat{\phi}(\xi) \left({1\over |r\xi|}\right)^{{n-1\over 2}+\alpha} \cos\left[2\pi |r\xi|-\left({n-1\over 2}+\alpha-{1\over2}\right) {\pi\over 2}-{\pi\over 4}\right] d\xi.
\end{array}
\eeq
Next, there are  finitely many 
\bel{Q_r v k}
\begin{array}{lr}\ds
{^\flat}\Q^{\alpha~\nu}_{r~j~k}(x)~=~
\\ \ds
~~~~(2\pi)^{-2k}{\a_k\over \pi}\int_{\R^n}e^{2\pi\i x\cdot\xi} \varphi^\nu_j(\xi)\Hat{\phi}(\xi) \left({1\over |r\xi|}\right)^{{n-1\over 2}+\alpha+2k} \cos\left[2\pi |r\xi|-\left({n-1\over 2}+\alpha-{1\over2}\right) {\pi\over 2}-{\pi\over 4}\right] d\xi
\\ \ds
+~(2\pi)^{-2k+1}{\b_k\over \pi}\int_{\R^n}e^{2\pi\i x\cdot\xi} \varphi^\nu_j(\xi)\Hat{\phi}(\xi) \left({1\over |r\xi|}\right)^{{n-1\over 2}+\alpha+2k-1} \sin\left[2\pi |r\xi|-\left({n-1\over 2}+\alpha-{1\over2}\right) {\pi\over 2}-{\pi\over 4}\right] d\xi
\end{array}
\eeq
where $\a_k, \b_k$ are constants given in (\ref{J asymptotic}) for $k=1,2,\ldots,L$. 

Finally, the remainder term 
\bel{remainder term}
\begin{array}{cc}\ds
{^\flat}\hbox{\bf R}^{\alpha~\nu}_{r~j~L}(x)~=~\int_{\R^n}e^{2\pi\i x\cdot\xi} \varphi^\nu_j(\xi)\Hat{\phi}(\xi) \H^\alpha\Big(2\pi |r\xi|\Big) d\xi,
\\\\ \ds
 \partial_\xi^\gamma \H^\alpha(\rho)~=~\O \left(\rho^{-2L-{1\over 2}}\right),\qquad  \rho\mt\infty
\end{array}
\eeq
where the implied constant  is bounded by $\C_{\gamma~\Re\alpha} e^{\c|\Im\alpha|}$.

Recall $\varphi^\nu_j$  defined in  (\ref{phi^v_j}). We find $\partial_\xi \varphi^\nu_j( \xi)=0$ whenever $\left|{\xi\over |\xi|}-\xi^\nu_j\right|\leq2^{-j/2}$. A direct computation shows
$\left|\p_\xi^\gamma \varphi_j^{\nu}(\xi)\right|\leq\C_{\gamma} 2^{|\gamma| j/2}$ 
for every multi-index $\gamma$ whenever ${1\over 3}<|\xi|\leq3$.

We have $|r|\approx 2^j$ if $\lambda_{m-1}\leq|r|<\lambda_m$. An $(n+1)$-fold integration by parts $w.r.t~ \xi$ in (\ref{remainder term}) gives
\bel{R EST}
\begin{array}{cc}\ds
\left| {^\flat}\hbox{\bf R}^{\alpha~\nu}_{r~j~L}(x)\right|~\leq~\C_{\Re\alpha}~ e^{\c|\Im\alpha|}~\left({1\over 1+|x|}\right)^{n+1}~2^{-\big[2L+{1\over 2}-{n+1\over 2}\big]j}.
\end{array}
\eeq
\begin{remark} There are a constant multiple of $2^{\big({n-1\over 2}\big)j}$ many $\xi^\nu_j\in\mathds{S}^{n-1}$. Moreover, $0\leq \Psi^{\alpha~\nu}_{j~m}(x) \leq1$. By using (\ref{R EST}), we obtain
 \[
\begin{array}{lr}\ds
\sum_{\nu\colon\xi^\nu_j\in\mathds{S}^{n-1}} \int_{\R^n}\left| {^\flat}\hbox{\bf R}^{\alpha~\nu}_{r~j~L}(x)\right| \left[1-\Psi^{\alpha~\nu}_{j~m}(x)\right]dx 
~\leq~\C_{\Re\alpha}~e^{\c|\Im\alpha|}~2^{- j/2} 2^{-\big[2L-n\big]j},
\end{array}
\]
\[
\sum_{\nu\colon\xi^\nu_j\in\mathds{S}^{n-1}}  \int_{\R^n}\Psi^{\alpha~\nu}_{j~m}(x)\left| {^\flat}\hbox{\bf R}^{\alpha~\nu}_{r~j~L}(x)\right|dx~\leq~\C_{\Re\alpha}~e^{\c|\Im\alpha|}~2^{- j/2} 2^{-\big[2L-n\big]j}.\]
\end{remark}
Let $\Re\alpha>0$ and $0<\Re\beta<{1\over2}$. We aim to show
\bel{Result U flat Q}
\begin{array}{lr}\ds
\sum_{\nu\colon\xi^\nu_j\in\mathds{S}^{n-1}} \int_{\R^n}\left|{^\flat}\Q^{\alpha~\nu}_{r~j}(x)\right| \left[1-\Psi^{\alpha~\nu}_{j~m}(x)\right]dx 
~\leq~\C_{N~\Re \alpha} ~e^{\c|\Im \alpha|}~ 2^{-N j} 2^{-\ve j},\qquad \hbox{\small{$N\ge0$}}
\end{array}
\eeq
and
\bel{Result V flat Q}
\sum_{\nu\colon\xi^\nu_j\in\mathds{S}^{n-1}}  \int_{\R^n}\Psi^{\alpha~\nu}_{j~m}(x)\left|{^\flat}\Q^{\alpha~\nu}_{r~j}(x)\right|dx~\leq~\C_{\Re \alpha} ~e^{\c|\Im \alpha|}~  2^{-\ve j}
\eeq
for some $\ve=\ve(\Re\alpha)>0$ whenever $\lambda_{m-1}\leq|r|<\lambda_m$.
\begin{remark}  Later, it should be clear that (\ref{Result U flat Q})-(\ref{Result V flat Q}) hold if  ${^\flat}\Q^{\alpha~\nu}_{r~j}$  is replaced by ${^\flat}\Q^{\alpha~\nu}_{r~j~k}$. Furthermore, an extra  term $2^{-j\big[2k-1\big]}$ appears on the right hand side of these inequalities:
\[
\begin{array}{lr}\ds
\sum_{\nu\colon\xi^\nu_j\in\mathds{S}^{n-1}} \int_{\R^n}\left|{^\flat}\Q^{\alpha~\nu}_{r~j~k}(x)\right| \left[1-\Psi^{\alpha~\nu}_{j~m}(x)\right]dx 
~\leq~\C_{N~\Re \alpha} ~e^{\c|\Im \alpha|}~ 2^{-N j} 2^{-\ve j}2^{-j\big[2k-1\big]},\qquad \hbox{\small{$N\ge0$}},
\end{array}
\]
\[
\sum_{\nu\colon\xi^\nu_j\in\mathds{S}^{n-1}}  \int_{\R^n}\Psi^{\alpha~\nu}_{j~m}(x)\left|{^\flat}\Q^{\alpha~\nu}_{r~j~k}(x)\right|dx~\leq~\C_{\Re \alpha} ~e^{\c|\Im \alpha|}~  2^{-\ve j} 2^{-j\big[2k-1\big]}.\]
\end{remark}
From (\ref{Result U flat Q})-(\ref{Result V flat Q}) and {\bf Remark 6.1}-{\bf Remark 6.2}, we conclude (\ref{Result Three U r})-(\ref{Result Three V r}).

\section{Proof of Proposition Two}
\setcounter{equation}{0}
Let $\Re\alpha>0$ and $0<\Re\beta<{1\over2}$.
Recall ${^\flat}\Q^{\alpha~\nu}_{r~j}(x)$ defined in (\ref{Q_r v}):
\[
\begin{array}{lr}\ds
{^\flat}\Q^{\alpha~\nu}_{r~j}(x)~=~{1\over \pi}\int_{\R^n}e^{2\pi\i x\cdot\xi} \varphi^\nu_j(\xi)\Hat{\phi}(\xi) \left({1\over |r\xi|}\right)^{{n-1\over 2}+\alpha} \cos\left[2\pi |r\xi|-\left({n-1\over 2}+\alpha-{1\over2}\right) {\pi\over 2}-{\pi\over 4}\right] d\xi.
\end{array}
\]
Because of Euler's formulae, we replace the cosine function  with $e^{\pm2\pi\i |r\xi|}$  multiplied by  $e^{\pm\pi\i \big[{1\over 4}+{1\over 2}\big({n-1\over 2}+\alpha-{1\over2}\big)\big]}$ whose amplitude is controlled by $\C_{\Re\alpha}e^{\c|\Im\alpha|}$.    
This assertion leads us to study
\bel{A flat j nu}
\begin{array}{ccc}\ds
{^\flat}\A^{\alpha}_{r~j}(x)~=~\sum_{\nu\colon\xi^\nu_j\in\mathds{S}^{n-1}} {^\flat}\A^{\alpha~\nu}_{r~j}(x),
\\\\ \ds
{^\flat}\A^{\alpha~\nu}_{r~j}(x)~=~\int_{\R^n}e^{2\pi\i \big[x\cdot\xi- r|\xi|\big]} \varphi^\nu_j(\xi) \Hat{\phi}(\xi)\left({1\over r|\xi|}\right)^{{n-1\over 2}+\alpha}  d\xi,
\\\\ \ds
{^\flat}\B^{\alpha}_{r~j}(x)~=~\sum_{\nu\colon\xi^\nu_j\in\mathds{S}^{n-1}} {^\flat}\B^{\alpha~\nu}_{r~j}(x),
\\\\ \ds
{^\flat}\B^{\alpha~\nu}_{r~j}(x)~=~\int_{\R^n}e^{2\pi\i \big[x\cdot\xi+ r|\xi|\big]} \varphi^\nu_j(\xi) \Hat{\phi}(\xi)\left({1\over r|\xi|}\right)^{{n-1\over 2}+\alpha}  d\xi
\end{array}
\eeq
where $\lambda_{m-1}\leq r<\lambda_m$ for $\lambda_m\in[2^{j-1},2^j]$ as before.

Note that $\supp\Hat{\phi}\subset\Big\{\xi\in\R^n\colon {1\over 3}<|\xi|\leq3\Big\}$  and
$\ds\supp\varphi^\nu_j~\subset~\Gamma^\nu_j~=~\Bigg\{\xi\in\R^n\colon \left|{\xi\over |\xi|}-\xi^{\nu}_{j}\right|<2^{-j/2+1}\Bigg\}$.
We have
\bel{A nu flat norm}
\begin{array}{lr}\ds
\left|{^\flat}\A^{\alpha~\nu}_{r~j}(x)\right|~\leq~~r^{-\big[{n-1\over 2}+\Re\alpha\big]} \int_{\supp\Hat{\phi}\varphi^\nu_j} \left({1\over |\xi|}\right)^{{n-1\over 2}+\Re\alpha}  d\xi
\\\\ \ds~~~~~~~~~~~~~~~~
~\leq~\C_{\Re\alpha} ~2^{-\big[{n-1\over 2}+\Re\alpha\big]j} 2^{-\left({n-1\over 2}\right)j},
\\\\ \ds
\left|{^\flat}\B^{\alpha~\nu}_{r~j}(x)\right|~\leq~~r^{-\big[{n-1\over 2}+\Re\alpha\big]} \int_{\supp\Hat{\phi}\varphi^\nu_j} \left({1\over |\xi|}\right)^{{n-1\over 2}+\Re\alpha}  d\xi
\\\\ \ds~~~~~~~~~~~~~~~~
~\leq~\C_{\Re\alpha} ~2^{-\big[{n-1\over 2}+\Re\alpha\big]j} 2^{-\left({n-1\over 2}\right)j}.
\end{array}
\eeq
Recall $\Psi^{\alpha}_{j~m}$ and $\Psi^{\alpha~\nu}_{j~m}$ defined in (\ref{Psi mu j}). We find 
\bel{supp Psi a u j m}
\vol~\supp  \Psi^{\alpha~\nu}_{j~m}~\lesssim~ 2^{\left({n-1\over 2}\right) j}2^{ n\epsilon j}.
\eeq
There are a constant multiple of $2^{\big({n-1\over 2}\big)j}$ many $\xi^\nu_j\in\mathds{S}^{n-1}$. 

We have
\bel{A flat j nu norm Sum}
\begin{array}{lr}\ds
\sum_{\nu\colon\xi^\nu_j\in\mathds{S}^{n-1}} \int_{\R^n} \Psi^{\alpha~\nu}_{j~m}(x)\left|{^\flat}\A^{\alpha~\nu}_{r~j}(x)\right| dx
\\\\ \ds
~\leq~\sum_{\nu\colon\xi^\nu_j\in\mathds{S}^{n-1}} \int_{\supp \Psi^{\alpha~\nu}_{j~m}} \left|{^\flat}\A^{\alpha~\nu}_{r~j}(x)\right|dx
\\\\ \ds
~\leq~\C_{\Re\alpha} \sum_{\nu\colon\xi^\nu_j\in\mathds{S}^{n-1}} 2^{-\big[{n-1\over 2}+\Re\alpha\big]j} 2^{-\left({n-1\over 2}\right)j}~2^{\left({n-1\over 2}\right)j}2^{ n\epsilon j}
\qquad
\hbox{\small{by (\ref{A nu flat norm})-(\ref{supp Psi a u j m})}}
\\\\ \ds
~\leq~\C_{\Re\alpha}~2^{\left({n-1\over 2}\right)j} 2^{-\big[{n-1\over 2}+\Re\alpha\big]j} 2^{-\left({n-1\over 2}\right)j}~2^{\left({n-1\over 2}\right)j}2^{n\epsilon j}
\\\\ \ds
~=~\C_{\Re\alpha}~2^{-\big[\Re\alpha-n\sigma\big]j}.
\end{array}
\eeq
Clearly, the same estimates in (\ref{A nu flat norm})-(\ref{A flat j nu norm Sum}) apply to ${^\flat}\B^{\alpha}_{r~j}(x)$. We obtain (\ref{Result V flat Q}) by choosing $n\sigma<{1\over2}\Re\alpha$.

\v
{\bf Lemma Two}~~{\it Given $j>0$ and $m=1,\ldots,M$, we have
\bel{Result Two A flat}
\begin{array}{cc}\ds
\sum_{\nu\colon\xi^\nu_j\in\mathds{S}^{n-1}}\int_{\R^n} \left|{^\flat}\A^{\alpha~\nu}_{r~j}(x)\right|\left[1-\Psi^{\alpha~\nu}_{j~m}(x)\right] dx
~\leq~\C_{N~\Re\alpha}~e^{\c|\Im\alpha|}~ 2^{-Nj} 2^{-\ve j},\qquad N\ge0,
\\\\ \ds
\sum_{\nu\colon\xi^\nu_j\in\mathds{S}^{n-1}}\int_{\R^n} \left|{^\flat}\B^{\alpha~\nu}_{r~j}(x)\right|\left[1-\Psi^{\alpha~\nu}_{j~m}(x)\right] dx
~\leq~\C_{N~\Re\alpha}~e^{\c|\Im\alpha|}~ 2^{-Nj} 2^{-\ve j},\qquad N\ge0
\end{array}
\eeq
for  some $\ve=\ve(\Re\alpha)>0$ whenever $\lambda_{m-1}\leq r<\lambda_m$.}

\subsection{Proof of Lemma Two}
For symmetry reason, we only focus on ${^\flat}\A^{\alpha}_{r~j}$.

Let ${^\flat}\A^{\alpha}_{r~j}$, ${^\flat}\A^{\alpha~\nu}_{r~j}(x)$ defined in (\ref{A flat j nu}). By taking  $\xi\mt r^{-1}\xi$, we find
\bel{A flat v}
\begin{array}{cc}\ds
{^\flat}\A^{\alpha}_{r~j}(x)~=~\sum_{\nu\colon \xi^\nu_j\in\mathds{S}^{n-1}} {^\flat}\A^{\alpha~\nu}_{r~j}(x),
\\\\ \ds
{^\flat}\A^{\alpha~\nu}_{r~j}(x)~=~r^{-n}\int_{\R^n}e^{2\pi\i \big[\left(r^{-1}x\right)\cdot\xi- |\xi|\big]} \varphi^\nu_j(\xi)\Hat{\phi}\left(r^{-1}\xi\right) \left({1\over |\xi|}\right)^{{n-1\over 2}+\alpha}  d\xi.
\end{array}
\eeq
Given $\xi^\nu_j\in\mathds{S}^{n-1}$, $\L_\nu$ is an $n\times n$-orthogonal matrix with $\det\L_\nu=1$ and   
\[\L_\nu^T \xi^\nu_j=(1,0)^T\in\R\times\R^{n-1}.\]
Let $\xi=\L_\nu \eta$ and  $x=\L_\nu u$ in (\ref{A flat v}). We have
 \bel{A flat v rewrite}
 \begin{array}{cc}\ds
 {^\flat}\A^{\alpha~\nu}_{r~j}(\L_\nu u)~=~r^{-n}\int_{\R^n}e^{2\pi\i \big[\left(r^{-1}u\right)\cdot\eta- |\eta|\big]} \varphi^\nu_j(\L_\nu \eta)\Hat{\phi}\left(r^{-1}\eta\right) \left({1\over |\eta|}\right)^{{n-1\over 2}+\alpha}  d\eta.
 \end{array}
 \eeq 
Note that  $\varphi^\nu_j$ is defined in (\ref{phi^v_j}).  $\Hat{\phi}$ defined in (\ref{hat phi}) is radially symmetric. 
\begin{remark} Write $\eta=(\eta_1,\eta')\in\R\times\R^{n-1}$. 
 We have 
 ${1\over 3} r<\eta_1<3 r$ and $ |\eta'|\lesssim 2^{j/2}$
  whenever $\varphi^\nu_j(\L_\nu \eta)\Hat{\phi}\left(r^{-1}\eta\right) \neq0$. 
\end{remark}
Consider
 \bel{Theta v j}
\begin{array}{lr}\ds
{^\flat}\A^{\alpha~\nu}_{r~j}(\L_\nu u)~=~r^{-n}\int_{\R^n}e^{2\pi\i  \big[\left(r^{-1}u\right)\cdot\eta- \eta_1\big]}\Theta^{\alpha~\nu}_{r~j}(\eta) d\eta,
\\\\ \ds
\Theta^{\alpha~\nu}_{r~j}(\eta)~=~e^{-2\pi\i \big[|\eta|-\eta_1\big]} \varphi^\nu_j(\L_\nu\eta)\Hat{\phi}\left(r^{-1}\eta\right) \left({1\over |\eta|}\right)^{{n-1\over 2}+\alpha}.
 \end{array}
\eeq
{\bf Remark 7.1} implies
\bel{vol Theta}
\vol ~\Bigg\{\supp \Theta^{\alpha~\nu}_{r~j}\Bigg\}~\lesssim~2^j 2^{\left({n-1\over 2}\right)j}.
\eeq
Recall $\lambda_{m-1}\leq r<\lambda_m$ where $\lambda_m\in[2^{j-1},2^j]:$  
$2^{\epsilon j-1}\leq\lambda_m-\lambda_{m-1}<2^{\epsilon j}$ where $0<\sigma<{1\over 2}$ can be chosen sufficiently small.
For brevity, write 
\[\lambda~=~\lambda_m.\]
We aim to show
\bel{Diff Ineq Theta N}
\left| \left(\lambda\partial_{\eta_1}\right)^N \left(\lambda^{1\over 2} \partial_{\eta'}\right)^\gamma  \Theta^{\alpha~\nu}_{r~j}(\eta)\right|~\leq~\C_{N~\gamma~\Re\alpha}~ e^{\c|\Im\alpha|}~2^{-\big[{n-1\over 2}+\Re\alpha\big]j}
\eeq
for every $N\ge0$ and every multi-index $\gamma$.
\v

{\bf 1.}~~Denote 
\[\chi(\eta)~=~ |\eta|-\eta_1,\qquad \eta_1>0.\] 
We claim
\bel{d Est Psi}
\left|\p_\eta^\gamma\chi(\eta)\right|~\leq~\C_\gamma~2^{-\gamma_1 j}2^{-\big[|\gamma|-\gamma_1\big]j/2}
\eeq
for every multi-index $\gamma$ with $|\gamma|\ge1$ whenever 
\[\L_\nu\eta~\in~\Gamma_j^\nu\cap \left\{{1\over 3} r<|\eta|\leq 3r\right\}.\]
Note that $\chi(\eta)$ is homogeneous of degree $1$ in $\eta$: $\chi(\rho\eta)=\rho\chi(\eta),~\rho>0$.  We have
\bel{Psi est1}
\left|\partial_\eta^\gamma \chi(\eta)\right|~\leq~\C_\gamma~|\eta|^{1-|\gamma|}
\eeq
for every multi-index $\gamma$.

Indeed, $\partial_\eta^\gamma \chi(\rho\eta)=\rho^{|\gamma|}\partial_\xi^\gamma\chi(\xi)$ for $\xi=\rho\eta$. On the other hand, $\partial_\eta^\gamma \chi(\rho\eta)=\partial_\eta^\gamma\left[\rho\chi(\eta)\right]=\rho\partial_\eta^\gamma\chi(\eta)$. Choose $\rho=|\eta|^{-1}$. We find $\partial_\eta^\gamma\chi(\eta)=|\eta|^{1-|\gamma|} \partial_\xi^\gamma\chi(\xi)$ where $|\xi|=1$.

Suppose $|\gamma|-\gamma_1\ge2$. From (\ref{Psi est1}), we find
\bel{Psi est2}
\begin{array}{lr}\ds
\left|\partial_\eta^\gamma \chi(\eta)\right|~\leq~\C_\gamma~2^{-\big(|\gamma|-1\big)j}
\\\\ \ds~~~~~~~~~~~~~~
~=~\C_\gamma~2^{-\gamma_1 j}2^{-\big[|\gamma|-\gamma_{1}-1\big]j} 
\\\\ \ds~~~~~~~~~~~~~~
~\leq~\C_\gamma~2^{-\gamma_1 j}2^{-\big[|\gamma|-\gamma_{1}\big]j/2}. 
\end{array}
\eeq
Suppose $|\gamma|-\gamma_{1}=0$ or $1$. Denote  $\eta^1=(\eta_1,0)^T\in\R\times\R^{n-1}$.
A direct computation shows
\[
 \begin{array}{cc}\ds
\left( \nabla_\eta \chi\right)(\eta^1)~=~0,
\\\\ \ds
\Big(\partial_{\eta_{1}}^N \nabla_{\eta'}\chi\Big)(\eta^1)
 ~=~ \Big(\nabla_{\eta'}\partial_{\eta_{1}}^N \chi\Big)(\eta^1)~=~0,\qquad \hbox{\small{$N\ge0$}}.
 \end{array}
\]
By writing out  the Taylor expansion of $\partial_{\eta_{1}}^N \chi(\eta)$ and $\nabla_{\eta'}\partial_{\eta_{1}}^N \chi(\eta)$  in the $\eta'$-subspace and using $\Big(\nabla_{\eta'}\partial_{\eta_{1}}^N \chi\Big)(\eta^1)=0$, we have
 \bel{p_eta_1 N norm}
 \begin{array}{lr}\ds
 \partial_{\eta_{1}}^N \chi(\eta)~=~ \O\left(|\eta'|^2|\eta|^{-N-1}\right),
\qquad
 \nabla_{\eta'} \partial_{\eta_{1}}^N \chi(x,\eta)~=~ \O\left(|\eta'||\eta|^{-N-1}\right),\qquad \hbox{\small{$N\ge1$}}.
 \end{array}
  \eeq
From (\ref{p_eta_1 N norm}) and {\bf Remark 7.1}, we conclude
\[
  \left|\partial_{\eta_{1}}^N \chi(\eta)\right|~\leq~\C_N  2^{- Nj},\qquad \left|\nabla_{\eta'}\partial_{\eta_{1}}^N \chi(\eta)\right|~\leq~\C_N  2^{-N j}2^{-j/2},\qquad \hbox{\small{$N\ge1$}}.
\]
\v
 {\bf 2.}~~ For $\Hat{\phi}$ defined in (\ref{hat phi}), we easily find
\bel{Diff Ineq phi}
\left| \partial_\eta^\gamma \left\{\Hat{\phi}\left(r^{-1}\eta\right) \left({1\over |\eta|}\right)^{{n-1\over 2}+\alpha}\right\}\right|~\leq~\C_{\gamma~\Re\alpha} ~e^{\c|\Im\alpha|}~2^{-\big[{n-1\over 2}+\Re\alpha\big]j} 2^{-|\gamma|j}
\eeq
for every multi-index $\gamma$. 
  \v
  
{\bf 3.}~~Recall $\varphi^\nu_j$  defined in  (\ref{phi^v_j}). We have
\[
\varphi^\nu_j(\L_\nu \eta)~=~{\ds\varphi\Bigg[2^{j/2}\left|{\eta\over |\eta|}-\eta^{\nu}_{j}\right|\Bigg]\over \ds \sum_{\nu\colon\xi^\nu_j\in\mathds{S}^{n-1}} \varphi\Bigg[2^{j/2}\left|{\eta\over |\eta|}-\eta^{\nu}_{j}\right|\Bigg]}.
\]
Observe that $\partial_\eta \varphi^\nu_j(\L_\nu \eta)=0$ whenever $\left|{\eta\over |\eta|}-\eta^\nu_j\right|\leq2^{-j/2}$.

Suppose $\eta_1\approx2^j$ and $|\eta'|\lesssim2^{j/2}$. 
 From direct computation, we find
\bel{deri computation 1}
\begin{array}{lr}\ds
{\partial\over \partial \eta_1}{\eta_1\over |\eta|}~=~{1\over |\eta|}-{\eta_1^2\over |\eta|^3}~=~
{|\eta'|^2\over |\eta|^3},
\\\\ \ds
{\partial\over \partial \eta_1}{\eta_i\over |\eta|}~=~-{\eta_1\eta_i\over |\eta|^3},\qquad i~=~2,\ldots,n
\end{array}
\eeq
and
\bel{deri eta norm 1}
\begin{array}{lr}\ds
\left| {\partial\over \partial \eta_1}{\eta\over |\eta|}\right|~\approx~2^{j/2}|\eta|^{-2}~\approx~2^{-j/2}|\eta|^{-1}~\approx~2^{-j/2}2^{-j},
\\\\ \ds
\left| \left({\partial\over \partial \eta_1}\right)^N{\eta\over |\eta|}\right|~\leq~\C_N~2^{-j/2}2^{-j} 2^{-(N-1)j}~=~\C_N~2^{-j/2} 2^{-Nj},\qquad N\ge1.
\end{array}
\eeq
By carrying out the differentiation $w.r.t$ $(\eta_1,\eta')^T$ and taking into account for $\left|{\eta\over |\eta|}-\eta^\nu_j\right|>2^{-j/2}$, we have
\bel{nu Diff est 1}
\begin{array}{cc}\ds
\left|\partial_{\eta_1}^N \partial_{\eta'}^\gamma \varphi_j^\nu \left(\L_\nu\eta\right)\right|~\leq~\C_{N~\gamma} ~2^{-Nj} 2^{|\gamma| j/2}2^{-|\gamma|j}
\end{array}
\eeq
for every $N\ge0$ and every multi-index $\gamma$.

By putting together the above estimates in {\bf 1}-{\bf 3}, we obtain (\ref{Diff Ineq Theta N}). 
\v

Define a differential operator
\bel{D}
\D_\lambda f~=~Id+\lambda^2 \p_{\eta_1}^2 f+\lambda \Delta_{\eta'} f.
\eeq
Let ${^\flat}\A^{\alpha~\nu}_{r~j}(\L_\nu u)$ be given in (\ref{Theta v j}). An $N$-fold integration by parts $w.r.t~\D_\lambda$ shows
\bel{INT by parts}
\begin{array}{lr}\ds
{^\flat}\A^{\alpha~\nu}_{r~j}(\L_\nu u)~=~
(2\pi\i)^{-2N} \left[1+\left[ \left(\lambda r^{-1}\right)u_1- \lambda\right]^2+\lambda \sum_{i=2}^n \left(r^{-1}u_i\right)^2\right]^{-N} 
\\\\ \ds~~~~~~~~~~~~~~~~~~~~~~~~~
r^{-n}\int_{\R^n}e^{2\pi\i  \big[\left(r^{-1}u\right)\cdot\eta- \eta_1\big]}\D^N_\lambda \Theta^{\alpha~\nu}_{r~j}(\eta) d\eta.
\end{array}
\eeq
By using (\ref{vol Theta})-(\ref{Diff Ineq Theta N}), we have
\bel{Int norm}
\begin{array}{lr}\ds
\left|\int_{\R^n}e^{2\pi\i  \big[\left(r^{-1}u\right)\cdot\eta- \eta_1\big]}\D^N_\lambda \Theta^{\alpha~\nu}_{r~j}(\eta) d\eta\right|~\leq~\C_{N~\Re\alpha}~ e^{\c|\Im\alpha|}~2^{-\big[{n-1\over 2}+\Re\alpha\big]j} 2^j 2^{\left({n-1\over 2}\right)j},\qquad N\ge0.
\end{array}
\eeq
Recall  $\Psi^{\alpha~\nu}_{j~m}$ defined in (\ref{Psi mu j}):
\[ \begin{array}{cc}\ds
\Psi^{\alpha}_{j~m}(x)~=~\sum_{\nu\colon\xi^\nu_j\in\mathds{S}^{n-1}}  \Psi^{\alpha~\nu}_{j~m}(x),
 \\\\ \ds 
 \Psi^{\alpha~\nu}_{j~m}(x)~=~\varphi\left[ 2^{-\epsilon j-1}\left|\lambda_m-\left(\L_\nu^T x\right)_1\right|\right] 
\prod_{i=2}^n \varphi\left[2^{-\left[{1\over 2}+\epsilon\right]j}\left|\left(\L_\nu^T x\right)_i\right|\right]
\\ \ds~~~~~~~~~~~~~
~+~\varphi\left[ 2^{-\epsilon j-1}\left|\lambda_m+\left(\L_\nu^T x\right)_1\right|\right] 
\prod_{i=2}^n \varphi\left[2^{-\left[{1\over 2}+\epsilon\right]j}\left|\left(\L_\nu^T x\right)_i\right|\right].
\end{array}\]
We find
$0\leq\Psi^{\alpha~\nu}_{j~m}(x)\leq1$.
Define
\bel{Tilde Psi}
\Tilde{\Psi}^{\alpha~\nu}_{j~m}(x)~=~\left\{\begin{array}{lr}\ds 1,\qquad \Psi^{\alpha~\nu}_{j~m}(x)=1,
\\ \ds
0,\qquad \Psi^{\alpha~\nu}_{j~m}(x)\neq1.
\end{array}\right.
\eeq
Observe that $\Tilde{\Psi}^{\alpha~\nu}_{j~m}\left(\L_\nu u\right)$ is supported in the rectangles
\bel{Rec}
\Tilde{\hbox{R}}_{j~m}~=~\Tilde{\hbox{R}}_{j~m}^+\cup \Tilde{\hbox{R}}_{j~m}^-,
\eeq
\bel{Rec+-}
\begin{array}{lr}\ds
\Tilde{\hbox{R}}_{j~m}^+~=~
\Bigg\{u\in\R^n\colon \lambda-2^{\epsilon j+1}\leq u_1\leq \lambda+2^{\epsilon j+1}, ~~|u_i|\leq 2^{\big[{1\over 2}+\epsilon\big] j},~~i=2,\ldots,n\Bigg\},
\\\\ \ds
\Tilde{\hbox{R}}_{j~m}^-~=~
\Bigg\{u\in\R^n\colon -\lambda-2^{\epsilon j+1}\leq u_1\leq -\lambda+2^{\epsilon j+1}, ~~|u_i|\leq 2^{\big[{1\over 2}+\epsilon\big] j},~~i=2,\ldots,n\Bigg\}
\end{array}
\eeq
From (\ref{Tilde Psi})-(\ref{Rec+-}), we have
\bel{INT SUM EST}
\begin{array}{lr}\ds
\sum_{\nu\colon\xi^\nu_j\in\mathds{S}^{n-1}} \int_{\R^n} \left| {^\flat}\A^{\alpha~\nu}_{r~j}(x)\right| \left[ 1-\Psi^{\alpha~\nu}_{j~m}(x)\right]dx
\\\\ \ds~~~~~~~
~\leq~\sum_{\nu\colon\xi^\nu_j\in\mathds{S}^{n-1}}\int_{\R^n} \left| {^\flat}\A^{\alpha~\nu}_{r~j}(x)\right| \left[ 1-\Tilde{\Psi}^{\alpha~\nu}_{j~m}(x)\right]dx
\\\\ \ds~~~~~~~
~=~\sum_{\nu\colon\xi^\nu_j\in\mathds{S}^{n-1}}\int_{\R^n\setminus  \supp \Tilde{\Psi}^{\alpha~\nu}_{j~m}} \left| {^\flat}\A^{\alpha~\nu}_{r~j}(x)\right| dx
\\\\ \ds~~~~~~~
~=~\sum_{\nu\colon\xi^\nu_j\in\mathds{S}^{n-1}}\int_{\R^n\setminus  \Tilde{\hbox{\small{R}}}_{j~m}} \left| {^\flat}\A^{\alpha~\nu}_{r~j}(\L_\nu u)\right| du
\\\\ \ds~~~~~~~
~\leq~\sum_{\nu\colon\xi^\nu_j\in\mathds{S}^{n-1}}\int_{\R^n\setminus  \Tilde{\hbox{\small{R}}}_{j~m}^+} \left| {^\flat}\A^{\alpha~\nu}_{r~j}(\L_\nu u)\right| du.
\end{array}
\eeq
Let $u\notin\Tilde{\hbox{R}}_{j~m}^+$ defined in (\ref{Rec+-}).  Suppose $u_1\notin[\lambda-2^{\epsilon j+1}, \lambda+2^{\epsilon j+1}]$.
We find
\bel{lambda r ratio}
\begin{array}{lr}\ds
\Big(\lambda r^{-1}\Big)u_1-\lambda~>~\Big[\lambda+2^{\epsilon j+1}\Big]-\lambda
\\\\ \ds~~~~~~~~~~~~~~~~~~~~~
~=~2^{\epsilon j+1},\qquad \hbox{\small{$u_1>\lambda+2^{\epsilon j+1}$}};
\\\\ \ds
\Big(\lambda r^{-1}\Big)u_1-\lambda~<~{\lambda\over \lambda-2^{\epsilon j}} \Big[\lambda-2^{\epsilon j+1}\Big]-\lambda
\\\\ \ds~~~~~~~~~~~~~~~~~~~~~
~=~\lambda \left\{ {\lambda-2^{\epsilon j+1}\over \lambda-2^{\epsilon j}}-1\right\}
\\\\ \ds~~~~~~~~~~~~~~~~~~~~~
~=~\lambda \left[ { -2^{\epsilon j+1}+2^{\epsilon j}\over \lambda-2^{\epsilon j}}\right]
\\\\ \ds~~~~~~~~~~~~~~~~~~~~~
~<~-2^{\epsilon j},\qquad\hbox{\small{$u_1< \lambda-2^{\epsilon j+1}$}}.
\end{array}
\eeq
There are a constant multiple of $ 2^{\left({n-1\over 2}\right)j}$ many $\xi^\nu_j\in\mathds{S}^{n-1}$. 
By using (\ref{INT by parts})-(\ref{Int norm}), we have
\bel{INT EST}
\begin{array}{lr}\ds
\sum_{\nu\colon\xi^\nu_j\in\mathds{S}^{n-1}}\int_{\R^n\setminus  \Tilde{\hbox{\small{R}}}_{j~m}^+} \left| {^\flat}\A^{\alpha~\nu}_{r~j}(\L_\nu u)\right| du
\\\\ \ds

~\leq~\C_{N~\Re\alpha}~ e^{\c|\Im\alpha|} \sum_{\nu\colon\xi^\nu_j\in\mathds{S}^{n-1}}r^{-n}2^{-\big[{n-1\over 2}+\Re\alpha\big]j}~ 2^j 2^{\left({n-1\over 2}\right)j}
\\ \ds~~~~~~~~~~~~~~~~~~~~~~~~~~~~~
\int_{\R^n\setminus  \Tilde{\hbox{\small{R}}}_{j~m}^+}\left[1+\left[ \left(\lambda r^{-1}\right)u_1- \lambda\right]^2+\lambda \sum_{i=2}^n \left(r^{-1}u_i\right)^2\right]^{-N}  du
\\\\ \ds
~\leq~\C_{N~\Re\alpha}~ e^{\c|\Im\alpha|}~2^{\left({n-1\over 2}\right)j}~2^{-nj}2^{-\big[{n-1\over 2}+\Re\alpha\big]j} ~ 2^j2^{\left({n-1\over 2}\right)j}
\\ \ds~~~~~~~~~~~~~~~~~~~~~~~~~~~~~
\int_{\R^n\setminus  \Tilde{\hbox{\small{R}}}_{j~m}^+}\left[1+\left[ \left(\lambda r^{-1}\right)u_1- \lambda\right]^2+\lambda \sum_{i=2}^n \left(r^{-1}u_i\right)^2\right]^{-N} du
\\\\ \ds
~=~\C_{N~\Re\alpha}~ e^{\c|\Im\alpha|}~2^{-\big[{n-1\over 2}+\Re\alpha\big]j} \int_{\R^n\setminus  \Tilde{\hbox{\small{R}}}_{j~m}^+}\left[1+\left[ \left(\lambda r^{-1}\right)u_1- \lambda\right]^2+\lambda \sum_{i=2}^n \left(r^{-1}u_i\right)^2\right]^{-N}  du
\\\\ \ds
~=~\C_{N~\Re\alpha}~ e^{\c|\Im\alpha|}~2^{-\big[{n-1\over 2}+\Re\alpha\big]j} r^{{n-1\over 2}}
\\ \ds~~~~~
\int_{\R^n\setminus \big\{\lambda-2^{\epsilon j+1}\leq u_1\leq \lambda+2^{\epsilon j+1},~|u_i|\leq (2^jr^{-1})^{1\over 2}2^{\epsilon j}\big\}}\left[1+\left[ \left(\lambda r^{-1}\right)u_1- \lambda\right]^2+\left(\lambda r^{-1}\right) \sum_{i=2}^n u_i^2\right]^{-N}  du
\\ \ds~~~~~~~~~~~~~~~~~~~~~~~~~~~~~~~~~~~~~~~~~~~~~~~~~~~~~~~~~~~~~~~~~~~~~~~~~~~~~~~~~~~~~~~~~~~~~\hbox{\small{$u_i\mt r^{1\over 2} u_i$, $i=2,\ldots,n$}}
\\\\ \ds
~\approx~\C_{N~\Re\alpha}~ e^{\c|\Im\alpha|}~2^{-\big[{n-1\over 2}+\Re\alpha\big]j} 2^{\big({n-1\over 2}\big)j}
\\ \ds~~~~~
\int_{\R^n\setminus \big\{\lambda-2^{\epsilon j+1}\leq u_1\leq \lambda+2^{\epsilon j+1},~|u_i|\leq (2^jr^{-1})^{1\over 2}2^{\epsilon j}\big\}}\left[1+\left[ \left(\lambda r^{-1}\right)u_1- \lambda\right]^2+\left(\lambda r^{-1}\right) \sum_{i=2}^n u_i^2\right]^{-N}  du
\\\\ \ds
~\leq~\C_{N~\Re\alpha}~ e^{\c|\Im\alpha|}~2^{-\big[{n-1\over 2}+\Re\alpha\big]j} 2^{\big({n-1\over 2}\big)j} 2^{-2\epsilon \big[N-n-1\big]j}
\\ \ds~~~~~
\int_{\R^n\setminus \big\{\lambda-2^{\epsilon j+1}\leq|u_1|\leq \lambda+2^{\epsilon j+1},~|u_i|\leq (2^jr^{-1})^{1\over 2}2^{\epsilon j}\big\}}\left[1+\left[ \left(\lambda r^{-1}\right)u_1- \lambda\right]^2+\left(\lambda r^{-1}\right) \sum_{i=2}^n u_i^2\right]^{-n-1}  du
\\ \ds~~~~~~~~~~~~~~~~~~~~~~~~~~~~~~~~~~~~~~~~~~~~~~~~~~~~~~~~~~~~~~~~~~~~~~~~~~~~~~~~~~~~~~~~~~
 \hbox{\small{by (\ref{lambda r ratio}) and $\lambda\approx r\approx 2^j$}}
\\\\ \ds
~=~\C_{N~\Re\alpha}~ e^{\c|\Im\alpha|}~2^{-\Re\alpha j} 2^{-2\epsilon \big[N-n-1\big]j}\int_{\R^n} \left[1+|u_1- r|^2+\sum_{i=2}^n u_i^2\right]^{-n-1} du.
\end{array}
\eeq
From (\ref{INT SUM EST})-(\ref{INT EST}), by taking $N$ sufficiently large depending $\epsilon=\sigma(\Re\alpha)>0$, we conclude
\bel{Final}
\sum_{\nu\colon\xi^\nu_j\in\mathds{S}^{n-1}} \int_{\R^n} \left| {^\flat}\A^{\alpha~\nu}_{r~j}(x)\right| \left[ 1-\Psi^{\alpha~\nu}_{j~m}(x)\right]dx ~\leq~\C_{N~\Re\alpha}~ e^{\c|\Im\alpha|}~2^{-\Re\alpha j} 2^{-N j},\qquad\hbox{\small{$N\ge 0$}}
\eeq
as desired.

\appendix 
  \section{Some estimates regarding Bessel functions}
 \setcounter{equation}{0}
$\bullet$ For  $a>-{1\over 2}, b\in\R$ and $\rho>0$,  a Bessel function $\J_{a+\i b}$ has an integral formula 
\bel{Bessel}
\J_{a+\i b}(\rho)~=~{(\rho/2)^{a+\i b}\over \pi^{1\over 2}\Gamma\left(a+{1\over 2}+\i b\right)} \int_{-1}^1 e^{\i \rho s } (1-s^2)^{a-{1\over 2}+\i b} ds.
\eeq
$\bullet$ For every $a>-{1\over 2}, b\in\R$ and $\rho>0$, we have 
\bel{J asymptotic}
\begin{array}{lr}\ds
\J_{a+\i b}(\rho) ~\sim~\left({\pi \rho\over 2}\right)^{-{1\over 2}} \cos\left[\rho-(a+\i b) {\pi\over 2}-{\pi\over 4}\right]
\\\\ \ds~~~~~~~~~~~~
~+~\left({\pi \rho\over 2}\right)^{-{1\over 2}} \sum_{k=1}^\infty \cos\left[\rho-(a+\i b) {\pi\over 2}-{\pi\over 4}\right] \a_k \rho^{-2k}+\sin\left[\rho-(a+\i b) {\pi\over 2}-{\pi\over 4}\right]~ \b_k \rho^{-2k+1}
\end{array}
\eeq
where
\bel{a_kb_k}
\begin{array}{cc}\ds
\a_k~=~(-1)^k[a+\i b, 2k]2^{-2k},\qquad \b_k~=~(-1)^{k+1}[a+\i b, 2k-1]2^{-2k+1},
\\\\ \ds
[a+\i b,m]~=~{\Gamma\left({1\over 2}+a+\i b+m\right)\over m!\Gamma\left({1\over 2}+a+\i b-m\right)},\qquad  m~=~0,1,2,\ldots
\end{array}
\eeq
in the sense of that
\bel{J O}
\begin{array}{lr}\ds
\left({d\over d\rho}\right)^\ell\left\{ \J_{a+\i b}(\rho) ~-~\left({\pi \rho\over 2}\right)^{-{1\over 2}} \cos\left[\rho-(a+\i b) {\pi\over 2}-{\pi\over 4}\right]\right.
\\ \ds~~~
\left.~-~\left({\pi \rho\over 2}\right)^{-{1\over 2}}\sum_{k=1}^N  \cos\left[\rho-(a+\i b) {\pi\over 2}-{\pi\over 4}\right] \a_k \rho^{-2k}+\sin\left[\rho-(a+\i b) {\pi\over 2}-{\pi\over 4}\right]~ \b_k \rho^{-2k+1}\right\}
\\\\ \ds
~=~\O\left(\rho^{-2N-{1\over 2}}\right)\qquad N\ge0,\qquad \ell\ge0,\qquad \rho\mt\infty.
\end{array}
\eeq
The implied constant is bounded by $\C_a e^{\c|b|}$.

By putting together (\ref{Bessel}) and (\ref{J asymptotic})-(\ref{J O}), we find a norm estimate in below.

$\bullet$ For  every $a>-{1\over 2}$, $b\in\R$ and $\rho>0$, 
\bel{J norm}
\left| {1\over \rho^{a+\i b}}~\J_{a+\i b}(\rho)\right|~\leq~\C_a~\left({1\over 1+\rho}\right)^{{1\over 2}+a}e^{\c|b|}.
\eeq  
\begin{remark}  In fact, (\ref{J norm}) is  true for every $a, b\in\R$ and $\rho>0$:
\[
\left| {1\over \rho^{a+\i b}}~\J_{a+\i b}(\rho)\right|~\leq~\C_a~\left({1\over 1+\rho}\right)^{{1\over 2}+a}e^{\c|b|},\qquad a,b\in\R,\qquad \rho>0.
\]
\end{remark}
This is a consequence of using the following identity.

$\bullet$  For every $a,b\in\R$ and $\rho>0$, we have
\bel{J identity}
\J_{a-1+\i b} (\rho)~=~2\left[{a+\i b\over \rho}\right] \J_{a+\i b}(\rho)-\J_{a+1+\i b}(\rho).
\eeq
More discussion of Bessel functions can be found in  the book of Watson \cite{Watson}.
\v
Let $z\in\Cx$.  
$\Omega^z$  is a distribution defined by analytic continuation from
\bel{Omega^z}
\begin{array}{ccc}\ds
\Re z<1,\qquad\Omega^z(x)~=~\pi^{-z}\Gamma^{-1}\left(1-z\right)  \left({1\over 1-|x|^2}\right)^{z}_+.
\end{array}
\eeq
$\bullet$  $\Omega^z$ can be equivalently defined by 
 \bel{Omega^z Transform} 
\begin{array}{lr}\ds
\Hat{\Omega}^z(\xi)~=~\left({1\over|\xi|}\right)^{{n\over 2}-z} \J_{{n\over 2}-z}\Big(2\pi|\xi|\Big),\qquad z\in\Cx.
\end{array}
\eeq
\begin{remark} From (\ref{Bessel}), (\ref{J asymptotic})-(\ref{J O}) and (\ref{J identity}), we conclude that $\Hat{\Omega}^z(\xi)$  for given $\xi\in\R^n$   is an analytic function of $z\in\Cx$.
\end{remark} 
Denote $\omega_{n-2}=2 \pi^{n-1\over 2} \Gamma^{-1}\left({n-1\over 2}\right)$ which is the area of $\mathds{S}^{n-2}$.
From (\ref{Omega^z}), we have
\bel{Omega^lambda Fourier transform}
\begin{array}{lr}\ds
\Hat{\Omega}^z(\xi)~=~\pi^{-z}  \Gamma^{-1}\left(1-z\right)\int_{|x|<1} e^{-2\pi\i x\cdot\xi} \left({1\over 1-|x|^2}\right)^z dx
\\\\ \ds~~~~~~~~~
~=~\pi^{-z}  \Gamma^{-1}\left(1-z\right)\int_0^\pi\left\{\int_0^1 e^{-2\pi\i |\xi|r\cos\vartheta} (1-r^2)^{-z} r^{n-1}dr\right\} \omega_{n-2} \sin^{n-2}\vartheta d\vartheta
\\\\ \ds~~~~~~~~~
~=~\omega_{n-2}\pi^{-z}  \Gamma^{-1}\left(1-z\right)\int_{-1}^1\left\{\int_0^1 e^{2\pi\i|\xi|rs} (1-r^2)^{-z} r^{n-1}dr\right\} (1-s^2)^{n-3\over 2} ds~~~~ \hbox{\small{($-s=\cos\vartheta$)}}
\\\\ \ds~~~~~~~~~
~=~\omega_{n-2}\pi^{-z}  \Gamma^{-1}\left(1-z\right)\int_0^1\left\{\int_{-1}^1 e^{2\pi\i|\xi|rs} (1-s^2)^{n-3\over 2} ds\right\}(1-r^2)^{-z} r^{n-1}dr
\\\\ \ds~~~~~~~~~
~=~\omega_{n-2}\pi^{-z}  \Gamma^{-1}\left(1-z\right)\int_0^1\left\{\int_{-1}^1 \cos(2\pi|\xi|rs) (1-s^2)^{n-3\over 2} ds\right\}(1-r^2)^{-z} r^{n-1}dr.
\end{array}
\eeq
Recall the Beta function identity:
\bel{Beta}
{\Gamma(z)\Gamma(w)\over \Gamma(z+w)}~=~\int_0^1 r^{z-1}(1-r)^{w-1}dr
\eeq
for every $\Re z>0$ and $\Re w>0$.

By writing out the Taylor expansion of the cosine function inside (\ref{Omega^lambda Fourier transform}), we find
\bel{Omega^lambda Fourier transform Sum}
\begin{array}{lr}\ds
\omega_{n-2}\pi^{-z}  \Gamma^{-1}\left(1-z\right)\int_0^1\left\{\int_{-1}^1 \cos(2\pi|\xi|rs) (1-s^2)^{n-3\over 2} ds\right\}(1-r^2)^{-z} r^{n-1}dr
\\\\ \ds
=~\omega_{n-2}\pi^{-z}  \Gamma^{-1}\left(1-z\right)
\\ \ds~~~~~
\sum_{k=0}^\infty (-1)^k {(2\pi|\xi|)^{2k}\over (2k)!}\left\{\int_{-1}^1 s^{2k} (1-s^2)^{n-3\over 2} ds\right\}\left\{\int_0^1 r^{2k+n-1}(1-r^2)^{-z} dr\right\}
\\\\ \ds
=~{1\over 2}\omega_{n-2}\pi^{-z}  \Gamma^{-1}\left(1-z\right)
\\ \ds~~~~~
\sum_{k=0}^\infty (-1)^k {(2\pi|\xi|)^{2k}\over (2k)!}\left\{\int_0^1 t^{k+{1\over 2}-1} (1-t)^{{n-1\over 2}-1} dt\right\}\left\{\int_0^1\rho^{k+{n\over 2}-1}(1-\rho)^{1-z-1} d\rho\right\}
\\ \ds~~~~~~~~~~~~~~~~~~~~~~~~~~~~~~~~~~~~~~~~~~~~~~~~~~~~~~~~~~~~~~~~~~~~~~~~~~~~~~~~~~~~~~~~~~~~~~~~~~~~~~
\hbox{\small{($t=s^2$,~$\rho=r^2$)}}
\\\\ \ds
=~{1\over 2}\omega_{n-2}\pi^{-z}  \Gamma^{-1}\left(1-z\right)\sum_{k=0}^\infty (-1)^k {(2\pi|\xi|)^{2k}\over (2k)!}~{\Gamma\left(k+{1\over 2}\right)\Gamma\left({n-1\over 2}\right)\over\Gamma\left(k+{n\over 2}\right)}~{ \Gamma\left(k+{n\over 2}\right) \Gamma\left(1-z\right)         \over    \Gamma\left(k+{n\over 2}+1-z\right)} 
\\ \ds~~~~~~~~~~~~~~~~~~~~~~~~~~~~~~~~~~~~~~~~~~~~~~~~~~~~~~~~~~~~~~~~~~~~~~~~~~~~~~~~~~~~~~~~~~~~~~~~~~~~~~~~~~~~~~
 \hbox{\small{by (\ref{Beta})}}
\\\\ \ds
=~ \pi^{{n-1\over 2}-z}\sum_{k=0}^\infty (-1)^k {(2\pi|\xi|)^{2k}\over (2k)!}~{\Gamma\left(k+{1\over 2}\right)        \over    \Gamma\left(k+{n\over 2}+1-z\right)}.
\end{array}
\eeq
On the other hand, we have
\bel{Omega cos}
\begin{array}{lr}\ds
\left({1\over|\xi|}\right)^{{n\over 2}-z} \J_{{n\over 2}-z}\Big(2\pi|\xi|\Big)
~=~\hbox{\small{$\pi^{{n-1\over 2}-z}  \Gamma^{-1}\left({n+1\over 2}-z\right)$}}\int_{-1}^1 e^{2\pi\i  |\xi| s} (1-s^2)^{{n-1\over 2}-z} ds\qquad\hbox{\small{by (\ref{Bessel})}}
\\\\ \ds~~~~~~~~~~~~~
~=~\hbox{\small{$\pi^{{n-1\over 2}-z}  \Gamma^{-1}\left({n+1\over 2}-z\right)$}}\int_{-1}^1 \cos\left(2\pi  |\xi| s\right) (1-s^2)^{{n-1\over 2}-z} ds
\\\\ \ds~~~~~~~~~~~~~
~=~\hbox{\small{$\pi^{{n-1\over 2}-z}  \Gamma^{-1}\left({n+1\over 2}-z\right)$}} \sum_{k=0}^\infty (-1)^k {(2\pi|\xi|)^{2k}\over (2k)!}\int_{-1}^1 s^{2k} (1-s^2)^{{n-1\over 2}-z} ds
\\\\ \ds~~~~~~~~~~~~~
~=~\hbox{\small{$\pi^{{n-1\over 2}-z}  \Gamma^{-1}\left({n+1\over 2}-z\right)$}} \sum_{k=0}^\infty (-1)^k {(2\pi|\xi|)^{2k}\over (2k)!}\int_0^1 \rho^{k+{1\over 2}-1} (1-\rho)^{{n+1\over 2}-z-1} d\rho\qquad \hbox{\small{($\rho=s^2$)}}
\\\\ \ds~~~~~~~~~~~~~
~=~\pi^{{n-1\over 2}-z} \sum_{k=0}^\infty (-1)^k {(2\pi|\xi|)^{2k}\over (2k)!}~{\Gamma\left(k+{1\over 2}\right)    \over    \Gamma\left(k+{n\over 2}+1-z\right)}\qquad\hbox{\small{by (\ref{Beta}).}}
\end{array}
\eeq

\section{Fourier transform of $\Lambda^\alpha$}
\setcounter{equation}{0}
We derive the formula of $\Hat{\Lambda}^\alpha(\xi,\tau)$ in (\ref{Lambda Transform})
by following the lines in p. $253-284$,
Chapter III of Gelfand and Shilov \cite{Gelfand-Shilov}.

$\Hat{\hbox{U}}^\alpha$, $\Hat{\hbox{V}}^\alpha$ are distributions defined by analytic continuation from
\bel{UV}
\begin{array}{lr}\ds
\Hat{\hbox{U}}^\alpha(\xi,\tau)~=~\Gamma^{-1}\left(1-\alpha\right)
  \left({1\over \tau^2-|\xi|^2}\right)^\alpha_-,
\\\\ \ds
 \Hat{\hbox{V}}^\alpha(\xi,\tau)~=~\Gamma^{-1}\left(1-\alpha\right)
  \left( {1\over \tau^2-|\xi|^2}\right)^\alpha_+,
\qquad
  \hbox{\small{$\Re\alpha<1$}}.  
\end{array}
\eeq
Denote  $\lambda(\alpha)={n+1\over 2}-\alpha,~ \alpha\in\Cx$. $\Pi^\alpha$, $\Lambda^\alpha$ are distributions defined by analytic continuation from
\bel{PiLambda}
\begin{array}{lr}\ds
\Pi^\alpha(x,t)~=~\Gamma^{-1}\left(1-\lambda(\alpha)\right)
  \left({1\over t^2-|x|^2}\right)^{\lambda(\alpha)}_-,
\\\\ \ds
 \Lambda^\alpha(x,t)~=~\Gamma^{-1}\left(1-\lambda(\alpha)\right)
  \left({1\over t^2-|x|^2}\right)^{\lambda(\alpha)}_+,
\qquad
 \hbox{\small{$\Re\lambda(\alpha)<1$}}.  
  \end{array}
\eeq 
Let $a>0, b>0$ and
\bel{rho xi,tau}
\begin{array}{cc}\ds
\rho(\xi,\tau,a,b)~=~\Bigg\{ \left[\tau^2-|\xi|^2\right]^2+\left[a \tau^2+b|\xi|^2\right]^2\Bigg\}^{1\over 2},
\\\\ \ds
 \cos\theta(\xi,\tau,a,b)~=~{\tau^2-|\xi|^2\over\rho(\xi,\tau,a,b)}.
 \end{array}
\eeq
For $\Re\alpha<{n+1\over 2}$ and $|\xi|\neq|\tau|$, we assert
 \bel{P,R,ab}
 \begin{array}{lr}\ds
\Hat{\hbox{P}}^{\alpha~a~b}(\xi,\tau)~=~\left[ \tau^2-|\xi|^2+\i a \tau^2+\i b |\xi|^2\right]^{-\alpha}
\\\\ \ds~~~~~~~~~~~~~~~~~~
~=~\rho(\xi,\tau,a,b)^{-\alpha} e^{-\i\alpha \theta(\xi,\tau,a,b)},
\\\\ \ds
\Hat{\hbox{R}}^{\alpha~a~b}(\xi,\tau)~=~\left[ \tau^2-|\xi|^2-\i a \tau^2-\i b |\xi|^2\right]^{-\alpha}
\\\\ \ds~~~~~~~~~~~~~~~~~~
~=~\rho(\xi,\tau,a,b)^{-\alpha} e^{\i\alpha \theta(\xi,\tau,a,b)}.
\end{array}
\eeq
$\Hat{\hbox{P}}^{\alpha~a~b}$ and $\Hat{\hbox{R}}^{\alpha~a~b}$ are distributions defined in $\R^{n+1}$ agree with (\ref{P,R,ab}) whenever $|\xi|\neq|\tau|$. Define
\bel{PR Limit}
 \Hat{\hbox{P}}^\alpha~=~\lim_{a\mt0,~b\mt0} \Hat{\hbox{P}}^{\alpha~a~b},\qquad \Hat{\hbox{R}}^\alpha~=~\lim_{a\mt0,~b\mt0} \Hat{\hbox{R}}^{\alpha~a~b},\qquad \hbox{\small{$\Re\alpha<{n+1\over 2}$}}.
 \eeq
For $\Re\lambda(\alpha)<{n+1\over 2}$ and $|x|\neq|t|$, we consider
\bel{Phi,Psi,ab}
 \begin{array}{lr}\ds
\Phi^{\alpha~a~b}(x,t)~=~\left[ t^2-|x|^2+\i a t^2+\i b |x|^2\right]^{-\lambda(\alpha)}
\\\\ \ds~~~~~~~~~~~~~~~~~~
~=~\rho(x,t,a,b)^{-\lambda(\alpha)} e^{-\i\lambda(\alpha) \theta(x,t,a,b)},
\\\\ \ds
\Psi^{\alpha~a~b}(x,t)~=~\left[ t^2-|x|^2-\i a t^2-\i b |x|^2\right]^{-\lambda(\alpha)}
\\\\ \ds~~~~~~~~~~~~~~~~~~
~=~\rho(x,t,a,b)^{-\lambda(\alpha)} e^{\i\lambda(\alpha) \theta(x,t,a,b)}.
\end{array}
\eeq
$\Phi^{\alpha~a~b}$ and $\Psi^{\alpha~a~b}$ are distributions defined in $\R^{n+1}$ agree with (\ref{Phi,Psi,ab}) whenever $|x|\neq|t|$. Define
\bel{PhiPsi Limit}
\Phi^\alpha~=~\lim_{a\mt0,~b\mt0} \Phi^{\alpha~a~b},\qquad \Psi^\alpha~=~\lim_{a\mt0,~b\mt0} \Psi^{\alpha~a~b},\qquad\hbox{\small{$\Re\lambda(\alpha)<{n+1\over 2}$}}.
\eeq
\begin{remark}
$\Hat{\hbox{P}}^\alpha$, $\Hat{\hbox{R}}^\alpha$ are analytic for 
$\Re\alpha<{n+1\over 2}$. 
$\Phi^\alpha$, $\Psi^\alpha$ are analytic for $\Re\lambda(\alpha)<{n+1\over 2}$.

Regarding details  can be found in $2.2$, 
 Chapter III of  Gelfand and Shilov \cite{Gelfand-Shilov}.
 \end{remark}
Define $\Hat{\hbox{P}}^\alpha$ and $\Hat{\hbox{R}}^\alpha$ by analytic continuation from (\ref{PR Limit}).
Recall (\ref{UV}) and (\ref{rho xi,tau})-(\ref{PR Limit}). We find
\bel{PR}
\begin{array}{lr}\ds
 \Hat{\hbox{P}}^\alpha
 ~=~\Gamma\left(1-\alpha\right)\left[e^{-\i\pi \alpha}\Hat{\hbox{U}}^\alpha+ \Hat{\hbox{V}}^\alpha\right],
\qquad
\Hat{\hbox{R}}^\alpha~=~\Gamma\left(1-\alpha\right)\left[e^{\i\pi\alpha}\Hat{\hbox{U}}^\alpha+ \Hat{\hbox{V}}^\alpha\right].
\end{array}
\eeq
Define $\Phi^\alpha$ and $\Psi^\alpha$ by analytic continuation from (\ref{PhiPsi Limit}).
Recall (\ref{PiLambda}), (\ref{rho xi,tau}) and (\ref{Phi,Psi,ab})-(\ref{PhiPsi Limit}).
We have
\bel{PhiPsi} 
\begin{array}{lr}\ds
\Phi^\alpha
~=~\Gamma\left(1-\lambda(\alpha)\right)\left[e^{-\i\pi \lambda(\alpha)} \Pi^\alpha+\Lambda^\alpha\right],
\\\\ \ds
\Psi^\alpha
~=~\Gamma\left(1-\lambda(\alpha)\right)\left[e^{\i\pi \lambda(\alpha)} \Pi^\alpha+\Lambda^\alpha\right].
\end{array}
\eeq
For $\Im z\neq0, \Im w\neq0$ and $0<\Re\lambda(\alpha)<{n+1\over 2}$, we assert
\[
Q^{\alpha~z~w}(x,t)~=~\left\{{1\over z |x|^2+ wt^2}\right\}^{\lambda(\alpha)},\qquad \hbox{\small{$(x,t)\neq(0,0)$}}.
\]
Let $a>0, b>0$. From  direct computation, we find
\bel{Q ab Transform}
\begin{array}{lr}\ds
\Hat{Q}^{\alpha~\i a~\i b}(\xi,\tau)~=~(-\i)^{\lambda(\alpha)}\iint_{\R^{n+1}} e^{-2\pi\i \big[x\cdot\xi+t\tau\big]} \left\{{1\over a |x|^2+ b t^2}\right\}^{\lambda(\alpha)} dxdt
\\\\ \ds~~~~~~~~~~~~~~~~~~~
~=~(-\i)^{\lambda(\alpha)}{1\over (\sqrt{a})^n\sqrt{b}}\iint_{\R^{n+1}} e^{-2\pi\i \big[x\cdot\xi/\sqrt{a}+t\tau/\sqrt{b}\big]} \left({1\over  |x|^2+  t^2}\right)^{\lambda(\alpha)} dxdt
\\ \ds~~~~~~~~~~~~~~~~~~~~~~~~~~~~~~~~~~~~~~~~~~ ~~~~~~~~~~~~~~~~~~~~~~~~~~~~~~~~~~~~~~~~~~~~
\hbox{\small{$x\mt x/\sqrt{a}$,~~$t\mt t/\sqrt{b}$}}
\\\\ \ds~~~~~~~~~~~~~~~~~~~
~=~(-\i)^{\lambda(\alpha)}{\pi^{-{n+1\over 2}+2\lambda(\alpha)}\over (\sqrt{a})^n\sqrt{b}}{\Gamma\left({n+1\over 2}-\lambda(\alpha)\right)\over \Gamma\left(\lambda(\alpha)\right)}\left\{{1\over  |\xi|^2/a+  \tau^2/b}\right\}^{{n+1\over 2}-\lambda(\alpha)},
\\ \ds~~~~~~~~~~~~~~~~~~~~~~~~~~~~~~~~~~~~~~~~~~~~~~~~~~~~~~~~~~~~~~~~~~~~~~~~~~~~~~~~~~~~~~~~~~~~~~~~~
\hbox{\small{$0<\Re\lambda(\alpha)<{n+1\over 2}$}}
\end{array}
\eeq
whenever $(\xi,\tau)\neq(0,0)$. 

\begin{remark}
The last equality in (\ref{Q ab Transform}) is derived by using 
\[
\int_{\R^\N}e^{-2\pi\i x\cdot\xi} |x|^{\gamma-\N}dx~=~{\pi^{\N-\gamma\over 2}\Gamma\left({\gamma\over 2}\right)\over \pi^{\gamma\over 2}\Gamma\left({\N-\gamma\over 2}\right)} |\xi|^{-\gamma},\qquad \hbox{\small{$\xi\neq0$}},
\]
for $0<\Re\gamma<\N$.
\end{remark}
Replace $a=-\i z$ and $b=-\i w$ inside (\ref{Q ab Transform}). We have
\bel{Q zv Transform}
\begin{array}{lr}\ds
\Hat{Q}^{\alpha~z~w}(\xi,\tau)~=~\iint_{\R^{n+1}} e^{-2\pi\i \left(x\cdot\xi+t\tau\right)} \left\{{1\over z |x|^2+ w t^2}\right\}^{\lambda(\alpha)} dxdt
\\\\ \ds~~~~~~~~~~~~~~~~~~
~=~{\pi^{\lambda(\alpha)-\alpha}\over (\sqrt{ z})^n\sqrt{ w}}{\Gamma\left(\alpha\right)\over \Gamma\left(\lambda(\alpha)\right)}\left\{{1\over  |\xi|^2/z+  \tau^2/w}\right\}^\alpha,\qquad \hbox{\small{$(\xi,\tau)\neq(0,0)$}}.
\end{array}
\eeq
\begin{remark}
(\ref{Q zv Transform}) is true  for every $\Im z>0$ and $\Im w>0$.
\end{remark}
Let $w=\i b$ fixed. 
Both sides of (\ref{Q zv Transform}) are analytic whenever $\Im z>0$.  Moreover, they are equal at $z=\i a$ for every $a>0$.  By analytic continuation, this equality holds for $\Im z>0$.
Given $z,~\Im z>0$,  a vice versa argument shows that (\ref{Q zv Transform}) is true for  $\Im w>0$. 

On the other hand, we find
\bel{Q ab Transform -}
\begin{array}{lr}\ds
\Hat{Q}^{\alpha~-\i a~-\i b}(\xi,\tau)~=~(-\i)^{-\lambda(\alpha)}{\pi^{-{n+1\over 2}+2\lambda(\alpha)}\over (\sqrt{a})^n\sqrt{b}}{\Gamma\left({n+1\over 2}-\lambda(\alpha)\right)\over \Gamma\left(\lambda(\alpha)\right)}\left\{{1\over  |\xi|^2/a+  \tau^2/b}\right\}^{{n+1\over 2}-\lambda(\alpha)},
\\\\ \ds~~~~~~~~~~~~~~~~~~~~~~~~~~~~~~~~~~~~~~~~~~~~~~~~~~~~~~~~~~~~~~~~~~~~~~~~\hbox{\small{$(\xi,\tau)\neq(0,0)$}},\qquad
\hbox{\small{$0<\Re\lambda(\alpha)<{n+1\over 2}$}}
\end{array}
\eeq
by carrying out the same estimate in (\ref{Q ab Transform}).

Replace $a=\i z$ and $b=\i w$ inside (\ref{Q ab Transform -}),  we obtain (\ref{Q zv Transform}) again. Furthermore, an analogue argument below {\bf Remark B.1} shows (\ref{Q zv Transform}) hold for $\Im z<0, \Im w<0$.

 Let
\bel{rho a}
\rho(a)~=~\sqrt{1+a^2},\qquad \cos\theta(a)~=~-1/\rho(a),\qquad \sin\theta(a)~=~a/\rho(a).
\eeq
Note that 
\[
(-1\pm\i a)^{n\over 2}~=~\rho(a)^{n\over 2} e^{\pm\i \theta \left({n\over 2}\right)}~\mt~ e^{\pm\i\pi \left({n\over 2}\right)}\qquad \hbox{\small{as}} \qquad a\mt0.
\]
Consider $z=-1+\i a$ and $w=1+\i b$ inside (\ref{Q zv Transform}).
Let $0<\Re\lambda(\alpha)<{n+1\over 2}$. We have
\bel{Q zv Transform +}
\begin{array}{lr}\ds
\Hat{Q}^{\alpha~z~w}(\xi,\tau)~\doteq~\Hat{Q}^{\alpha~a~b}(\xi,\tau)
\\\\ \ds~~~~~~~~~~~~~~~~~~~
~=~\iint_{\R^{n+1}} e^{-2\pi\i \left(x\cdot\xi+t\tau\right)} \left\{{1\over (-1+\i a) |x|^2+ (1+\i b) t^2}\right\}^{\lambda(\alpha)} dxdt
\\\\ \ds~~~~~~~~~~~~~~~~~~~
~=~{\pi^{\lambda(\alpha)-\alpha}\over (-1+\i a)^{n\over2}(1+\i b)^{1\over 2}}{\Gamma\left(\alpha\right)\over \Gamma\left(\lambda(\alpha)\right)}\left\{{(-1+\i a)(1+\i b)\over  (1+\i b)|\xi|^2+ (-1+\i a) \tau^2}\right\}^\alpha
\end{array}
\eeq
for $(\xi,\tau)\neq(0,0)$.

Similarly, take $z=-1-\i a$ and $w=1-\i b$ inside (\ref{Q zv Transform}).
We find
\bel{Q zv Transform -}
\begin{array}{lr}\ds
\Hat{Q}^{\alpha~z~w}(\xi,\tau)~\doteq~\Hat{Q}^{\alpha~a~b}(\xi,\tau)
\\\\ \ds~~~~~~~~~~~~~~~~~~~
~=~{\pi^{\lambda(\alpha)-\alpha}\over (-1-\i a)^{n\over2}(1-\i b)^{1\over 2}}{\Gamma\left(\alpha\right)\over \Gamma\left(\lambda(\alpha)\right)}\left\{{(-1-\i a)(1-\i b)\over  (1-\i b)|\xi|^2+ (-1-\i a) \tau^2}\right\}^\alpha,\quad \hbox{\small{$(\xi,\tau)\neq(0,0)$}}.
\end{array}
\eeq

Define $\Hat{Q}^{\alpha~a~b}$, $a>0, b>0,  0<\Re\lambda(\alpha)<{n+1\over 2}$ as a distribution in $\R^{n+1}$ agree with (\ref{Q zv Transform +})-(\ref{Q zv Transform -}) whenever $|\xi|\neq|\tau|$. 

Recall (\ref{P,R,ab})-(\ref{PhiPsi Limit}). 
By  taking $a\mt0, b\mt0$, we find
\bel{Phi, R, Psi, P, Transform}
\begin{array}{lr}\ds
\Hat{\Phi}^\alpha~=~\pi^{\lambda(\alpha)-\alpha}e^{-\i\pi\left({ n\over 2}\right)}~{\Gamma\left(\alpha\right)\over \Gamma\left(\lambda(\alpha)\right)}~\Hat{\hbox{R}}^\alpha,
\\\\ \ds
\Hat{\Psi}^\alpha~=~\pi^{\lambda(\alpha)-\alpha}e^{\i\pi\left({ n\over 2}\right)}~{\Gamma\left(\alpha\right)\over \Gamma\left(\lambda(\alpha)\right)}~\Hat{\hbox{P}}^\alpha,\qquad \hbox{\small{$0<\Re\lambda(\alpha)<{n+1\over 2}$}}.
\end{array}
\eeq
Moreover,  $\Hat{\hbox{P}}^\alpha, \Hat{\hbox{R}}^\alpha$ and $\Phi^\alpha, \Psi^\alpha$ can be given   in terms of $\Pi^\alpha, \Lambda^\alpha$ and $\hbox{U}^\alpha, \hbox{V}^\alpha$ as (\ref{PR}) and (\ref{PhiPsi}). A direct computation shows
\bel{Lambda computa}
\begin{array}{lr}\ds
\Hat{\Lambda}^\alpha\left[ e^{i\pi\lambda(\alpha)}-e^{-i\pi\lambda(\alpha)}\right]
~=~\Gamma^{-1}(1-\lambda(\alpha))\Bigg[ e^{\i\pi\lambda(\alpha)} \Hat{\Phi}^\alpha-e^{-\i\pi\lambda(\alpha)}\Hat{\Psi}^\alpha\Bigg]
\\\\ \ds~~~~~~~~~~~~~~~~~~~~~~~~~~~~~~~~~~~
~=~\Gamma^{-1}(1-\lambda(\alpha)) \Gamma^{-1}(\lambda(\alpha)) 
\\ \ds~~~~~~~~~~~~~~~~~~~~~~~~~~~~~~~~~~~~~~~~~
\pi^{\lambda(\alpha)-\alpha}\Gamma(\alpha) 
\Bigg[ e^{-\i\pi\big[{n\over 2}-\lambda(\alpha)\big]} \Hat{\hbox{R}}^\alpha-e^{\i\pi\big[{n\over 2}-\lambda(\alpha)\big]} \Hat{\hbox{P}}^\alpha\Bigg],
~~~~~~
 \hbox{\small{$0<\Re\lambda(\alpha)<{n+1\over 2}$}}.
\end{array}
\eeq

From (\ref{Lambda computa}), by using the identity 
$\Gamma^{-1}(1-z)\Gamma^{-1}(z)={ \sin \pi z\over \pi}$ for $z\in\Cx$, 
we obtain
\bel{Lambda Transform Formula}
\begin{array}{lr}\ds
\Hat{\Lambda}^\alpha
~=~\pi^{\lambda(\alpha)-\alpha-1} \Gamma\left(\alpha\right){1\over 2\i} \Bigg\{- e^{\i\pi \big[{n\over 2}-\lambda(\alpha)\big]} \Hat{\hbox{P}}^\alpha+e^{-\i\pi \big[{n\over 2}-\lambda(\alpha)\big] }\Hat{\hbox{R}}^\alpha\Bigg\}
\\\\ \ds~~~~~
~=~\pi^{\lambda(\alpha)-\alpha-1} \Gamma\left(\alpha\right)\Gamma(1-\alpha)
\\ \ds~~~~~~~~~~
{1\over 2\i} \Bigg\{- e^{\i\pi \big[{n\over 2}-\lambda(\alpha)\big]} \Big[e^{-\i\pi\alpha}\Hat{\hbox{U}}^\alpha+\Hat{\hbox{V}}^\alpha\Big]+e^{-\i\pi \big[{n\over 2}-\lambda(\alpha)\big] }\Big[e^{\i\pi\alpha}\Hat{\hbox{U}}^\alpha+\Hat{\hbox{V}}^\alpha\Big]\Bigg\}\qquad \hbox{\small{by (\ref{PR})}}
\\\\ \ds~~~~~
~=~\pi^{\lambda(\alpha)-\alpha-1} \Gamma\left(\alpha\right)\Gamma(1-\alpha)
\\ \ds~~~~~~~~~~
{1\over 2\i} \Bigg\{- e^{\i\pi \big[\alpha-{1\over2}\big]} \Big[e^{-\i\pi\alpha}\Hat{\hbox{U}}^\alpha+\Hat{\hbox{V}}^\alpha\Big]+e^{-\i\pi \big[\alpha-{1\over 2}\big] }\Big[e^{\i\pi\alpha}\Hat{\hbox{U}}^\alpha+\Hat{\hbox{V}}^\alpha\Big]\Bigg\}
\\\\ ~~~~~
~=~\pi^{\lambda(\alpha)-\alpha-1}\Gamma(\alpha)\Gamma\left(1-\alpha\right)
\Bigg\{ \Hat{\hbox{U}}^\alpha-\sin\pi\left(\alpha-{1\over 2}\right)\Hat{\hbox{V}}^\alpha\Bigg\} 
\end{array}
\eeq
for $\lambda(\alpha)={n+1\over 2}-\alpha\notin\Z$ and $0<\Re\lambda(\alpha)<{n+1\over 2}$. 

Recall {\bf Remark B.1}. $\Hat{\hbox{P}}^\alpha$, $\Hat{\hbox{R}}^\alpha$ are analytic for 
$\Re\alpha<{n+1\over 2}$. By the principle of analytic continuation, (\ref{Lambda Transform Formula}) holds for $0<\Re\lambda(\alpha)<{n+1\over 2}$.

Let $0<\Re\alpha<1$.  $\Hat{\hbox{U}}^\alpha$ and $\Hat{\hbox{V}}^\alpha$ agree with  $\Hat{\hbox{U}}^\alpha(\xi,\tau)$ and $\Hat{\hbox{V}}^\alpha(\xi,\tau)$ in (\ref{UV}) for $|\tau|\neq|\xi|$. We have
\[
\begin{array}{lr}
  \Hat{\Lambda}^\alpha(\xi,\tau)~=~
\pi^{\lambda(\alpha)-\alpha-1}\Gamma\left(\alpha\right)
 \left\{ ~{\ds \left({1\over \tau^2-|\xi|^2}\right)^\alpha_-} - \sin\pi\left(\alpha-{1\over 2}\right){\ds \left( {1\over \tau^2- |\xi|^2}\right)^\alpha_+}~\right\}.
  \end{array}
\]

{\bf Acknowledgement

~~~~~~I am deeply grateful to my advisor Elias M. Stein for those stimulating talks and unforgettable lectures.}

\small{School of Science, Westlake University, Hangzhou, 310030, China}\\      
\small{email: wangzipeng@westlake.edu.cn}

\end{document}